# Weak Error of Dean-Kawasaki Equation with Smooth Mean-Field Interactions[*]

by Ana Djurdjevac[a], Xiaohao Ji[a], Nicolas Perkowski[a]

a. Freie Universität Berlin

**Abstract**

We consider the weak-error rate of the SPDE approximation by regularized Dean-Kawasaki equation with Itô noise for particle systems with mean-field interactions both on the drift and the noise. The global existence and uniqueness of the corresponding SPDEs are established using the variational approach to SPDEs, and the weak-error rate is estimated using the technique of Kolmogorov equations on the space of probability measures.

In particular, the rate derived in this paper coincides with that is the previous work [DKP24], which considered free Brownian particles using Laplace duality.

**Keywords:** mean-field systems, Dean-Kawasaki equation, weak error rate, master equation, Lions calculus, variational approach to SPDE.
**A.M.S. subject classification:** 60H15, 60J60, 60K35

## Table of contents

## 1 Introduction

In this paper, we consider the stochastic partial differential equation (SPDE) approximation of the *mean-field* interacting particle system $X_t = (X_t^i)_{i=1}^N$ on $\mathbb{T}^d$, whose dynamic is given by the following distribution-dependent stochastic differential equations:

$$dX_t^i = b(X_t^i, \rho_{X_t})dt + \sqrt{2}\, \Sigma(X_t^i, \rho_{X_t}) \cdot dW_t^i, \qquad (1.1)$$

---
[*]. This document has been written using GNU T$_E$X$_{MACS}$ [H+98].

on $\mathbb{T}^d$, where $W = (W^i)_{i=1}^N$ is the $d \times N$-dimensional real Brownian motion, and

$$\mathbb{T}^d \times \mathcal{P}(\mathbb{T}^d) \ni (x, \mu) \mapsto b(x, \mu) \in \mathbb{R}^d, \quad (x, \mu) \mapsto \Sigma(x, \mu) \in \mathbb{R}^{d \times d}$$

are functions modelling the dependence of the dynamic of a single particle on the position of itself and the whole system through the empirical measure (in the mean-field sense)

$$\rho_t^N = \rho_{X_t} := \frac{1}{N} \sum_{i=1}^N \delta_{X_t^i}.$$

Following the seminal work of Kawasaki [Kaw73] and Dean [Dea96], a natural candidate of approximation is the conservative SPDE of the form

$$dm_t^N = \nabla^2 : [a(\cdot, m_t^N) m_t^N] dt + \nabla \cdot [b(\cdot, m_t^N) m_t^N] dt + \frac{2}{\sqrt{N}} \nabla \cdot [f^N(m_t^N) \Sigma(\cdot, m_t^N) \cdot dW_t^N],$$

where $m \mapsto \nabla^2 : [a(\cdot, m) m] + \nabla \cdot [b(\cdot, m) m]$ is the formal dual of the generator given by (1.1), i.e.

$$a(x, \mu) = \Sigma(x, \mu) \cdot \Sigma^t(x, \mu), \tag{1.2}$$

for $x \in \mathbb{T}^d$ and $\mu \in \mathcal{P}(\mathbb{T}^d)$, $f^N : \mathbb{R} \to \mathbb{R}^+ \cup \{0\}$ is a family of well-chosen approximation of the square root, and $(dW^N)_{N \leq 0}$ denotes a choice of family of noises colored in space and white in time which approximates the genuine space-time white noise. We refer the reader to the following sections for discussion on the derivation and assumptions, especially the mollifications of the white noise and of the square-root nonlinearity used in the original works of Kawasaki and Dean.

## 1.1 Main Results

The main results of this paper can be summarized as follows: Firstly, we show the existence and uniqueness of the probabilistically strong $L^2$ *solution* (to be specified at the beginning of Chapter 2) given any square-integrable probability density as initial datum for the SPDE (1.12), under the *stochastic parabolicity* condition (specified in (2.15)) and regularity assumptions 2.1 on both drift and the noise. The exact statement and the proof are given in Corollary 2.10.

Furthermore, we characterize the approximation of SPDEs (1.12) to the particle systems (1.8) by a *quantitative* bound on the *weak error*, that is if we assume that the coefficients in (1.1) and an arbitrarily given continuous functional $\mathcal{F}[\cdot]$ on $\mathcal{P}(\mathbb{T}^d)$ are all regular enough in the sense of Assumption 4.1, the weak error rate of the approximation of the interacting particle system (1.1) by the mollified Dean-Kawasaki equation (1.12) is upper bounded by

$$|\mathbb{E}\mathcal{F}(\rho_t^N) - \mathbb{E}\mathcal{F}(m_t^{N,\epsilon})| \lesssim N^{-1-2/(d+2)} \log N, \tag{1.3}$$

where the constant does not depend on the number of particles $N$. Especially, the rate coincides with that in [DKP24], where the Brownian particles evolve freely (i.e. $b[\cdot]$ and $\Sigma[\cdot]$ are independent of the empirical measures) and a duality argument in the spirit of [Eth00, KLvR19] can be applied. The dependence of the rate on the dimension of the ambient space is essentially due to the balance of an entropy bound (specifically Lemma 3.2) and the approximation of the covariance of the white noise by the colored ones, which can be well observed in Lemma 3.7, since morally the generators of the Dean-Kawasaki equations only converge to those of the particle systems when the mollification of the white noise is taken away as the number of particles go to infinity, while the entropy bound deteriorates simultaneously.

## 1.2 Literature Review and Notations

This paper is an extension of the previous work [DKP24] by H. Kremp and the first and third authors. Despite the widely recognised effectiveness of the Dean-Kawasaki equation as the approximation of various particle systems by conservative SPDE, a rigorous solution theory thereof is highly nontrivial, and has only been mostly resolved recently. The main reason is two-folded: On one hand, the equation is highly singular due to the irregularity of the white noise, which is even worsened by the derivative on the noise term. Precisely, the equation is scaling supercritical in the language of singular SPDEs (c.f. [Hai14, GIP15] and others). On the other hand, the square-root in the noise leads to degeneracy as the solution goes to zero, which makes the equation difficult to handle even if the white noise is replaced by mollified ones.

More surprisingly, it has been shown in [KLvR19, KLvR20] that for as weak a solution notion as the martingale solution formulated there, there exists only trivial solutions (those given by the empirical measures of the particles) to the original Dean-Kawasaki equation with suitable coupling constants. Furthermore, the solution is fragile when one changes the parameters, and as such, the equation is not suitable for a stable numerical approximation

Consequently, in order to achieve a non-trivial solution theory, one needs to approximate in a suitable way the noise for example, so that the equation is at least subcritical with regard to parabolic scaling. To deal with the degeneracy, one requires more involved techniques such as (stochastic) entropy solution and kinetic solution theory to tackle the square root, c.f. [FG19, DG20, GGK22]. Such a progress reaches its peak, as the well-posedness of the Dean-Kawasaki equation with truncated Stratonovich spacetime white noise and with the full square-root diffusion coefficient is shown in the seminal work [FG21]. The Stratonovich noise is essential, however, due to an essential cancellation in the a priori entropy estimates. The incorporation of interacting kernels is then considered in [WWZ22]. We refer to [FG21] and the references therein for more details on the well-posedness of related equations and large deviation principle under suitable scalings, where the latter is proved using the method of weak convergence (c.f. [DE11] and references there) and translating the problem into the well-posedness and stability (with respect to the control) of the skeleton equations.

Numerical schemes for Dean-Kawasaki-type equations were considered in [CS23, CF23]. [CS23] introduce a discontinuous Galerkin scheme for the regularized DK equation from [CSZ19, CSZ20]. In [CF23], finite element and finite difference approximations without interaction are considered and weak error estimates are proved with some weak distance that is parametrized by the Sobolev regularity of the test functions, and the rate measured in their distance can be arbitrarily high. On the other hand, the positivity for the approximations is not proved, and a strong assumption on the existence of lower and upper bounds for the solution of a discrete heat equation, which means that initially the particles must be fairly spread out, has to be imposed.

As mentioned above, one of the main contributions of this paper is the derivation of a weak error rate of the approximation of the weakly interacting diffusion with the corresponding SPDEs with thermal noise. The methodology is highly influenced by that in the literature of mean-field games, mean-field control, and propagation of chaos pioneered in [Lio, BLPR17, CDLL19, Tse21, CD22] and others, where the solution theory of the master equations plays a fundamental role. Especially, if one considers the mean-field limit of (1.1) with constant diffusion coefficient, that is

$$\partial_t \bar{m}_t = \Delta \bar{m}_t + \nabla \cdot [b(\cdot, \bar{m}_t) \bar{m}_t],$$

it is proven in [DT21] that uniformly in time the weak propagation of chaos holds, i.e.

$$\mathbf{E}\mathcal{F}[\rho_t^N] - \mathcal{F}[\bar{m}_t] = O(N^{-1}),$$

under regularity assumptions on test functional $\mathcal{F}[\cdot]$ and the drift. One can infer from the comparison between the bound above and (1.3) that the Dean-Kawasaki equation indeed capture the fluctuation of the mean-field system. The result is further extended for example to the kinetic setting and higher-order weak propagation of chaos by [BD24] implementing the Glauber Calculus [Due21]. Regarding the strong propagation of chaos and recent progress on the Gaussian fluctuations of the related models, we refer the readers to [BRTV98, BGV07, CGM08, FG15, Lac23], [TH81, Tan84, Szn84, BD24, MM22, GK24, MM24] and references therein for further details. Our derivation of the regularity of the solution to the Kolmogorov equation is closely related to that in [Tse21], where a recursive formula for the arbitrary high-order expansion of the McKean–Vlasov equation is derived. It turns out that for the problem here, a second-order linearization (with respect to the probability measure) suffices, while we need to exploit more spatial regularity.

To conclude the review of the literature, it is worth emphasizing that the Dean-Kawasaki equation could also play a canonical rule corresponding to the Wasserstein diffusion on the space of probability measures as crucially observed in [vRS09] (see also [DH24] for a different perspective). Due to the infinite-dimensional essence of the underlying space, the construction of the diffusion is highly nontrivial. Inspiring results in this direction have recently been achieved by [Del22, Del24] as a higher-dimensional extension of the modified massive Arratia flow constructed in [KvR19], where a recurrent diffusion process on the space of probability measures over an arbitrary closed Riemannian manifold of dimension $d \geq 2$ is constructed in [Del22]. An integration by part formula is achieved, which in turn gives the closability of the corresponding bilinear form. [Del24] further provides a unified framework to treat the massive particle systems, the Wasserstein Brownian motions and the Dean-Kawasaki equation. In particular, a rigorous interpretation of the Dean-Kawasaki equation

$$dm_t = \sum_{x \in m_t} \Delta \delta_x dt + \nabla \cdot (\sqrt{m_t} \cdot dW_t) + \text{interaction}$$

is established, where $dW$ is the space-time white noise and the interactions can be quite singular. We refer the readers to the original works for details. While the processes in [Del22, Del24] have likely captured more mesoscopic essence of Dean's original work [Dea96] and the structure of the underlying space, they might deviate from the SPDE modelling and are much harder to manipulate under, for example, the numerical consideration.

### 1.2.1 Notation

We use $\lesssim$ to denote less or equal than modulo a generic positive constant, while we use $C$ to denote some generic positive constant that we define recursively in the paper. We will fix some arbitrary $T \in (0, \infty)$ and consider solutions on the time interval $[0, T]$. $\mathbb{R}^d$ denotes the $d$-dimensional Euclidean space, while $\mathbb{T}^d$ the $d$-dimensional unit torus. $\mathbb{Z}^d$ is the $d$-dimensional lattice and $\mathbb{Z}_+^d$ denotes those nonzero elements in $\mathbb{Z}^d$ whose first nonzero entry is positive. We call a function kernel on $\mathbb{R}^d$ if it's smooth, non-negative, radially symmetric, supported in the unit ball of $\mathbb{R}^d$ and has finite moment up to the second order. For function spaces on $\mathbb{T}^d$, we use $W^{k,p}$ to denote the Sobolev spaces. We specify especially the inner product on $W^{-1,2}(\mathbb{T}^d)$ as

$$\langle u, v \rangle_{W^{-1,2}} = \langle (1-\Delta)^{-1} u, v \rangle_{L^2}. \tag{1.4}$$

We will denote $e_k(x) := e^{2\pi i k \cdot x}$ with $k \in \mathbb{Z}^d$ to be the orthonormal basis of $L^2(\mathbb{T}^d, \mathbb{C})$, while $(B_t)_{t \geq 0}$ (with possible superscript) denotes the $d$-dimensional complex-valued Brownian motions. For an aribtrary distribution $f$ on $\mathbb{T}^d$, we denotes its Fourier coefficients as $\hat{f}(k) = \mathfrak{F}f(k)$ for $k \in \mathbb{Z}^d$.

Slightly different from [KLvR19], in the proof we will only consider the calculus on the probability space instead of that on the Radon measures with finite total variation. We denote the space of probability measure on the torus (with Borel $\sigma$-algebra) as $\mathcal{P} := \mathcal{P}(\mathbb{T}^d)$. $\mathcal{P}(\mathbb{T}^d)$ also denotes the same set, endowed with the 2-Wasserstein distance given by

$$\mathcal{W}_2^2(\mu, \nu) := \inf_{\Pi \in \Pi(\mu, \nu)} \int_{\mathbb{T}^d \times \mathbb{T}^d} d_{\mathbb{T}^d}^2(x, y) d\Pi(x, y),$$

where $d_{\mathbb{T}^d}$ denotes the canonical distance on $\mathbb{T}^d$ and $\Pi(\mu, \nu)$ denotes the family of couplings of $\mu$ and $\nu$, that is we consider the weak topology on the probability measure spaces. We use $\rightharpoonup$ to indicate weak convergence on the space of probability measures. We name as usual functionals on $\mathcal{P}$ of the form $\mathcal{F} = f(\langle \varphi, \cdot \rangle)$ cylindrical, if $f: \mathbb{R}^m \to \mathbb{R}$ and $\varphi: \mathbb{T}^d \to \mathbb{R}^m$ are both smooth for some $m \in \mathbb{N}$, and $f$ has bounded derivatives of all orders. For differentiation, $\delta/\delta\mu$ is used to denote the extrinsic derivative of functionals on $\mathscr{P}(\mathbb{T}^d)$, while $D_\mu$ is used to denote the interior derivative (also Lions derivative in the literature). A detailed definition is given in the Appendix A.3.

For any $b \in \mathbb{R}^d$, we abbreviate the $\ell^2$ norm by $|b|^2 := \sum_{i=1}^d |b_i|^2$. For $d \times d$ matrices $A$ and $B$, we denote by $A^t$ the transpose of $A$ and by

$$A : B := \sum_{i,j=1}^d A_{ij} B_{ij},$$

the Frobenius product of $A$ and $B$. We denote the Frobenius norm of $d \times d$ matrix $A$ as $|A| := (\sum_{i,j=1}^{d} |A_{ij}|^2)^{1/2}$, for which we have

$$|A:B| \leq |A| \cdot |B|, \quad |A \cdot B| \leq |A| \cdot |B|, \quad |A b| \leq |A| \cdot |b|, \tag{1.5}$$

all following from the Cauchy-Schwarz inequality. For some normed function space $X$ on $\mathbb{T}^d$ to $\mathbb{R}$, we extend the definition to $d \times d$ matrix-valued functions on $\mathbb{T}^d$ by $\|F\|_X := |(\|F_{ij}\|_X)|$ for any such $F$. We also use $L_2$ to denote the Hilbert-Schmidt norm.

## 1.3 Formal Derivation

We can now specify the first-order weakly-interacting particle systems introduced in the introduction on the $d$-dimensional torus as

$$dX_t^i = b(X_t^i, \rho_{X_t}) dt + \sqrt{2}\, \Sigma(X_t^i, \rho_{X_t}) \cdot dW_t^i, \tag{1.6}$$

where $i \in \{1, ..., N\}$ is used as the particle index, $b: \mathbb{T}^d \times \mathcal{P}(\mathbb{T}^d) \to \mathbb{R}^d$, $\Sigma: \mathbb{T}^d \times \mathcal{P}(\mathbb{T}^d) \to \mathbb{R}^{d \times d}$ are coefficients of the SDEs, and $\{W^i\}_{i=1}^N$ are independent real Brownian motions on $\mathbb{T}^d$. All stochastic integrals are considered in Itô sense throughout the paper, and we recall that the empirical measure is defined by

$$\rho_{X_t} := \frac{1}{N} \sum_{i=1}^{N} \delta_{X_t^i} \in \mathcal{P}(\mathbb{T}^d). \tag{1.7}$$

Especially, if we assume that the coefficients are given by smooth interaction kernel $K$, $G$, s.t.

$$b(x, \mu) := K * \mu(x), \quad \Sigma(x, \mu) := \sigma(G * \mu(x)), \tag{1.8}$$

for $x \in \mathbb{T}^d$ and $\mu \in \mathcal{P}(\mathbb{T}^d)$ with appropriate $\sigma: \mathbb{R}^d \to \mathbb{R}^{d \times d}$, (1.1) reduces to the first-order mean-field diffusion system

$$dX_t^i = \frac{1}{N} \sum_{j=1}^{N} K(X_t^i - X_t^j) dt + \sqrt{2}\, \sigma\!\left(\frac{1}{N} \sum_{j=1}^{N} G(X_t^i - X_t^j)\right) \cdot dW_t^i, \quad i \in \{1, ..., N\}. \tag{1.9}$$

Since only regular enough kernels are considered in this work, it is unnecessary to exclude the particle itself in the interactions. In order to capture the fluctuation of the particle system (1.1), we follow the idea of Dean and Kawasaki [Kaw73, Dea96] by calculating the formal quadratic variation of the empirical measure: Using Itô's formula for $\langle f, \rho_{X_t} \rangle$ with test function $f \in C^\infty(\mathbb{T}^d)$, we have

$$d\langle f, \rho_{X_t}\rangle = \langle [b(\cdot, \rho_{X_t}) \cdot \nabla] f, \rho_{X_t}\rangle dt + \langle ([\{\Sigma(\cdot, \rho_{X_t}) \cdot \Sigma(\cdot, \rho_{X_t})^t\} : \nabla^2) f, \rho_{X_t}\rangle dt + dM_t,$$

where $(M_t)_{t \geq 0}$ is a martingale, and we recall that $:$ is the Frobenius product of matrices. The quadratic variation of $(M_t)_{t \geq 0}$ is further given by

$$\begin{aligned}
d[M]_t = d[\langle f, \rho_X \rangle]_t &= \frac{2}{N} \langle [\Sigma(\cdot, \rho_{X_t})^t \cdot \nabla] f, \rho_{X_t}\rangle dt \\
&= \frac{2}{N} \int_{\mathbb{T}^d} |\sqrt{\rho_{X_t}} [\Sigma(\cdot, \rho_{X_t})^t \cdot \nabla] f|^2 dx\, dt \\
&= \frac{2}{N} \sum_{k \in \mathbb{Z}^d} |\langle \sqrt{\rho_{X_t}} [\Sigma(\cdot, \rho_{X_t})^t \cdot \nabla] f, e_k \rangle_{L^2}|^2 dt,
\end{aligned}$$

where $\{e_k\}_{k\in\mathbb{Z}^d}$ denotes the canonical ortho-normal basis of $L^2(\mathbb{T}^d,\mathbb{C})$ by Plancherel identity. Recall from (1.2) that the diffusion coefficient is defined by

$$a(x,\mu) := \Sigma(x,\mu) \cdot \Sigma(x,\mu)^t,$$

and the above calculation suggests that we may consider the following SPDE:

$$d\rho_t^N = \nabla^2 : [a(\cdot,\rho_t^N)\rho_t^N]dt - \nabla \cdot [b(\cdot,\rho_t^N)\rho_t^N]dt + \frac{2}{\sqrt{N}}\nabla \cdot \left[(\rho_t^N)^{\frac{1}{2}}\Sigma(\cdot,\rho_t^N) \cdot dW_t\right], \quad (1.10)$$

where $(W_t)_{t\geq 0}$ is the $d$-dimensional real cylindrical Wiener process on $L^2(\mathbb{T}^d)$, since formally it captures the quadratic variation $\rho_X$ with $N$ particles.

However, as crucially observed in [KLvR19, KLvR20], the following triviality result prevents one from achieving a useful solution theory for at least a large class of (1.10). We note that although the differentiation notion here is slightly different from that in [KLvR20], where they assume that the functionals are defined on the larger space of Radon measures on $\mathbb{T}^d$, the result remains valid, since Itô's formula (c.f. Proposition 3.5 below) still holds and enables us to apply Girsanov's theorem. We refer the solution notion used in [KLvR19, KLvR20] as the martingale solution to the Dean-Kawasaki equation.

**Lemma 1.1. (Triviality [KLvR19, KLvR20])** *If we assume that the noise in (1.1) is additive, that is*

$$\Sigma \equiv \mathrm{Id}_{d\times d},$$

*and the drift is given as the interior derivative (whose definition is given in Appendix A) of a potential function, that is*

$$b(\cdot,\mu) := \nabla\frac{\delta\mathcal{K}}{\delta\mu}(\mu,\cdot) = D_\mu\mathcal{K}(\mu,\cdot),$$

*for a scalar-valued $\mathcal{K}$ on $\mathcal{P}(\mathbb{T}^d)$ twice differentiable with Lipschitz-continuous interior derivatives, then there exists unique martingale solution to (1.10), if and only if the initial data is assumed to be of the form*

$$\rho_0^N := \frac{1}{N}\sum_{j=1}^N \delta_{x_j^N},$$

*for some $\{x_j^N\}_{j=1}^N \subset \mathbb{T}^d$. In that case, the unique martingale solution to (1.10) is given by (1.7), the empirical measure of the interacting diffusion.*

In particular, if the energy functional $\mathcal{K}$ on $\mathcal{P}(\mathbb{T}^d)$ in consideration is defined to be

$$\mathcal{K}[\mu] := \frac{1}{2}\int_{\mathbb{T}^d\times\mathbb{T}^d} V_\mathcal{K}(x-y)d\mu(x)d\mu(y), \quad (1.11)$$

for some even potential $V_\mathcal{K} \in C^2(\mathbb{T}^d)$, the drift term in (1.10) is given by $b(\cdot,\mu) = \nabla V_\mathcal{K} * \mu(\cdot)$, which corresponds to

$$\frac{1}{N}\sum_{j=1}^N K(X_t^i - X_t^j) = \frac{1}{N}\sum_{j=1}^N \nabla V_\mathcal{K}(X_t^i - X_t^j)$$

in (1.8). Due to the above-mentioned triviality result, we consider instead the regularized equation

$$dm_t^N = \nabla \cdot [a(\cdot, m_t^N)\nabla m_t^N - \tilde{b}(\cdot, m_t^N) m_t^N] dt + \frac{2}{\sqrt{N}} \nabla \cdot [f^N(m_t^N) \Sigma(\cdot, m_t^N) \cdot dW_t^N], \quad (1.12)$$

following [FG21, DKP24], where the noise $(W_t^N)_{t \geq 0}$ is now colored in space and white in time approximating the cylindrical Wiener process, $\{f^N\}_{N \in \mathbb{N}}$ approximates the square root in an appropriate manner and the modified drift is defined by

$$\mathbb{T}^d \times \mathcal{P}(\mathbb{T}^d) \ni (x, \mu) \mapsto \tilde{b}(x, \mu) := b(x, \mu) - \sum_{j=1}^d \partial_j a_{\cdot, j}(x, \mu). \quad (1.13)$$

We specify the approximation of the square root by imposing the following condition on $f^\delta \in C^1(\mathbb{R}_{\geq 0}, \mathbb{R})$ s.t. $f^\delta(0) = 0$:

$$\|(f^\delta)'\|_{L^\infty} \lesssim \frac{1}{\sqrt{\delta}}, \quad (f^\delta)'(x) \lesssim \frac{1}{\sqrt{x}}, \quad |f^\delta(x)^2 - x| \lesssim \delta, \quad (1.14)$$

for all $x \geq 0$ and constants uniform in $\delta \in (0, 1)$. One choice of such $(f^\delta)$ is given by Example 2.2 in [DKP24]. We abbreviate $f^N := f^{\delta_N}$ for some sequence $(\delta_N)_{N \geq 0}$, i.e., to mollify the square root in such a way that it converges to the square-root when the number of particles approaches infinity.

**Remark 1.2.** The family of real-valued spatially colored noise can be represented with complex Brownian motions by

$$W_t^N := \sum_{k \in \mathbb{Z}^d} \theta_k^N e_k B_t^k, \quad (1.15)$$

as in [FG21, DKP24], where we abuse the notation to denote $\{B^k\}_{k \in \mathbb{Z}^d}$ as $d$-dimensional complex-valued Brownian motions s.t. $B^{-k} = \overline{B^k}$, and $B^j$ is independent of $B^k$ if and only if $j \neq \pm k$.

We further impose the following conditions on the covariance structure of the colored noise that we consider in (1.12):

**Assumption 1.3.** *We assume that the Fourier modes of the colored noise in (1.15) are given by*

$$\theta_k^N := \hat{\eta}^{\frac{1}{2}}\left(\frac{k}{L_N}\right) = \mathfrak{F}_{\mathbb{R}^d}[L_N^d \eta(L_N(\cdot))](k)^{\frac{1}{2}}, \quad (1.16)$$

*for some symmetric $\eta \in S(\mathbb{R}^d, \mathbb{R})$ s.t. $\hat{\eta} \in C_c^\infty(\mathbb{R}^d, [0, \infty))$ and $\hat{\eta}(0) = 1$, while $\{L_N\}_{N \geq 0}$ is chosen such that $0 \leq L_N \to \infty$ as $N \to \infty$.*

**Remark 1.4.** Since $\theta^N$ is apparently symmetric for any $N \in \mathbb{N}$, i.e. $\theta_k^N = \theta_{-k}^N$ for all $k \in \mathbb{Z}^d$ by construction, the representation (1.15) is real. It also follows directly from Assumption 1.3 that

$$\sum_{k \in \mathbb{Z}^d} |k|^2 |\theta_k^N|^2 < \infty, \forall N \in \mathbb{N}, \quad \lim_{N \to \infty} \theta_k^N = 1, \forall k \in \mathbb{Z}^d, \quad (1.17)$$

where the second property indicates that the family of colored noises approximates the singular space-time white noise, as the covariances approximate dirac. Moreover, we have the following bound for $\psi^N := \sum_{k \in \mathbb{Z}^d} |\theta_k^N|^2 e_k$:

$$\int_{\mathbb{T}^d} |x|^2 |\psi^N(x)| \mathrm{d}x \lesssim_\eta L_N^{-2}. \tag{1.18}$$

Indeed, it follows from Poisson summation that $\psi^N = \sum_{z \in \mathbb{Z}^d} L_N^d \eta[L_N(\cdot + z)]$, since $\eta$ is a Schwartz function, and therefore

$$\int_{[-\frac{1}{2}, \frac{1}{2}]^d} |x|^2 |\psi^N|(x) \mathrm{d}x \leq \sum_{z \in \mathbb{Z}^d} \int_{([-\frac{1}{2}, \frac{1}{2}]^d + z) L_N} \left|\frac{x}{L_N} - z\right|^2 |\eta(x)| \mathrm{d}x$$

$$\leq \sum_{z \in \mathbb{Z}^d} \int_{([-\frac{1}{2}, \frac{1}{2}]^d + z) L_N} \frac{|x|^2}{L_N^2} |\eta(x)| \mathrm{d}x = L_N^{-2} \int_{\mathbb{R}^d} |x| |\eta(x)| \mathrm{d}x.$$

**Remark 1.5.** Under our assumptions on $a$, $b$ and $\theta^N$ it may be possible to extend the stochastic kinetic solution theory of [DG20] to prove the well-posedness of equation (1.12) with $f^N$ replaced by the non-truncated square root, provided that Stratonovich noise is used instead of Itô noise. As in [WZ22], the large deviations of this equation for $N \to \infty$ may be related to that of the corresponding particle system. But for controlling the weak error at finite $N$, Itô noise is the natural choice. Simply replacing it by Stratonovich noise introduces an error of order $N^{-1}$ in equation (1.3). And with Itô noise and square root, the equation may lose its parabolic character, and we do not expect it to be solvable for generic initial data.

The rest of the paper is organized as follows: In Chapter 2, we prove the existence and uniqueness of the equation (1.12) in Corollary 2.10 under the regularity assumptions 2.1. The proof combines the variational solutions provided by [RSZ24] in a Gelfand tripe suitable for the gradient multiplicative structure of the noise, a taming procedure (2.13) and the $L^1$ conservation law for non-negative initial data in Remark 2.8.

In Chapter 3, we derive the Itô's formula for the dynamics in concern using the approximation arguments provided by A.3.2 following [CD22] and then quantify the comparison of the generators by Lemma 3.7 using the key technical Lemma A.3. Fisher information naturally appears and can be controlled by the entropy estimate Lemma 3.2 mimicking [FG21]. In Chapter 4, we derive the weak error estimate using the Kolmogorov equation together with Lemma 3.7, whose conditions (essentially the $W^{2,\infty}$-bound in (3.16)) are translated into the spatial regularity of the solution to the Kolmogorov equations, that is (4.19) in Corollary 4.6. The main result is then stated and proved in Theorem 4.7.

Several assumptions are included as the paper proceeds. For the existence and uniqueness of the SPDE (1.12), Assumption 2.1 and the stochastic parabolicity bound (2.15) are crucially, while for the weak-error estimate, Assumption 1.3 provides the quantitative approximation of the covariance structure of the cylindrical Wiener process, and Assumption 4.1 assures sufficient regularity on the test functionals and coefficients of the equations. The main Theorem 4.7 combines them all.

# 2 Existence and Uniqueness

In this section, we focus on the global existence and uniqueness of the probabilistically strong solution to (1.12) under sufficient conditions on the drift and noise. Given an SPDE with coefficients in a Gelfand triple $V \subset H \subset V^*$, we define the $H$ solution to the SPDE to be the probabilistically strong solution in $C([0,T], H) \cap L^2([0,T], V)$ for any $0 < T < \infty$.

Due to our later need for both notions of $W^{-1,2}$ and $L^2$ solutions to (1.12), we assume that $b[\cdot]$ and $\Sigma[\cdot]$ admit extension to functionals on

$$X := W^{-1,2}(\mathbb{T}^d) \cup L^1(\mathbb{T}^d),$$

where we recall that the inner product on $W^{-1,2}(\mathbb{T}^d)$ is defined in (1.4). From now on we denote for brevity

$$b[\mu](\cdot) := b(\cdot, \mu),$$

and similarly for the distribution-dependent functions $\Sigma[\cdot]$ and $\tilde{b}[\cdot]$. We assume the following strong conditions on the drift and noise coefficients, which are anyway necessary in the chapter of weak-error estimates.

**Assumption 2.1.** *We assume that $b[\cdot]$ defined on $X \times \mathbb{T}^d$ to $\mathbb{R}^d$ and $\Sigma[\cdot]$, $a[\cdot]$ defined on $X \times \mathbb{T}^d$ to $\mathbb{R}^{d \times d}$ satisfy the continuity bounds on $W^{-1,2}$*

$$\|b[\mu] - b[\nu]\|_{L^\infty} \lesssim \|\mu - \nu\|_{W^{-1,2}}, \tag{2.1}$$

$$\|\Sigma[\mu] - \Sigma[\nu]\|_{W^{1,\infty}} \lesssim (1 \vee \|\nu\|_{W^{-1,2}}) \|\mu - \nu\|_{W^{-1,2}}, \tag{2.2}$$

$$\|a[\mu] - a[\nu]\|_{W^{1,\infty}} \lesssim (1 \vee \|\nu\|_{W^{-1,2}}) \|\mu - \nu\|_{W^{-1,2}}, \tag{2.3}$$

*and further for some arbitrarily small $\varepsilon > 0$ that*

$$\|a[\mu]\|_{W^{1+\varepsilon,\infty}} \lesssim 1 \vee \|\mu\|_{W^{-1,2}}. \tag{2.4}$$

*We further assume the uniform ellipticity and boundedness conditions on diffusion matrices, i.e. there exists $a_{\min} > 0$, $C_\Sigma < \infty$ s.t. uniformly in $\mu \in X$,*

$$a[\mu] \geq a_{\min} \mathrm{Id}_{d \times d}, \quad \|\Sigma[\mu]\|_{L^\infty} \leq C_\Sigma, \tag{2.5}$$

*and lastly the growth bound on $L^1$ that*

$$\|b[\mu]\|_{L^\infty} \vee \|\Sigma[\mu]\|_{W^{1,\infty}} \lesssim 1 \vee \|\mu\|_{L^1}. \tag{2.6}$$

**Remark 2.2.** It follows directly from the continuity conditions (2.1-2.3) that

$$\|b[\mu]\|_{L^\infty} \lesssim 1 \vee \|\mu\|_{W^{-1,2}}, \tag{2.7}$$

$$\|\Sigma[\mu]\|_{W^{1,\infty}} \vee \|a[\mu]\|_{W^{1,\infty}} \lesssim 1 \vee \|\mu\|_{W^{-1,2}}, \tag{2.8}$$

while the stronger condition (2.4) provides a rather simplified proof of the uniform-ellipticity bound in the $L^2 \subset W^{-1,2} \subset (L^2)^*$ Gelfand triple in Lemma 2.4 below, where we denote $(L^2(\mathbb{T}^d))^*$ to be the space of continuous functionals on $L^2(\mathbb{T}^d)$ with the subspace topology inherited from $W^{-1,2}$.

We also notice that the continuity condition on $a[\cdot]$ is stronger than the other conditions would imply, since then we would get

$$\|a[\mu] - a[\nu]\|_{W^{1,\infty}} \lesssim (\|\Sigma[\mu]\|_{W^{1,\infty}} + \|\Sigma[\mu]\|_{W^{1,\infty}})\|\Sigma[\mu] - \Sigma[\nu]\|_{W^{1,\infty}}$$
$$\lesssim (1 + \|\nu\|_{W^{-1,2}}^2 + \|\mu\|_{W^{-1,2}}^2)\|\mu - \nu\|_{W^{-1,2}}. \tag{2.9}$$

The dependence on $\|\nu\|_{W^{-1,2}}^2$ and $\|\mu\|_{W^{-1,2}}^2$ in (2.8) would prevent the direct application of Theorem 3.2 in [RSZ24] (specifically the local monotonicity condition).

Meanwhile, the assumptions (2.1) - (2.6) can be verified with the canonical example we have in mind.

**Example 2.3.** Let the drift be given by the interaction kernel (1.11) and the noise by

$$\Sigma[\mu] = \sigma(G * \mu), \quad \sigma \in C_b^2(\mathbb{R}^l, \mathbb{R}^{d \times d}). \tag{2.10}$$

under the condition that $K \in W^{2,\infty}(\mathbb{T}^d, \mathbb{R}^d)$, $G \in W^{2+\varepsilon,\infty}(\mathbb{T}^d, \mathbb{R}^l)$ for some $l \in \mathbb{N}$ and

$$\sigma\sigma^t \geq a_{\min} \mathrm{Id}_{d \times d}.$$

The assumptions (2.1) - (2.6) hold, since we have for example

$$\|\sigma(G * \mu) - \sigma(G * \nu)\|_{W^{1,\infty}} = \|\sigma(G * \mu) - \sigma(G * \nu)\|_{L^\infty}$$
$$+ \|\nabla(G * \mu) \cdot (\nabla\sigma)(G * \mu) - \nabla(G * \nu) \cdot (\nabla\sigma)(G * \nu)\|_{L^\infty},$$

where under Einstein summation convention

$$(\nabla(G * \mu) \cdot (\nabla\sigma)(G * \mu))_i := \partial^j \sigma(G * \mu) \partial_i(G_j * \mu),$$

we can bound separately

$$\|\sigma(G * \mu) - \sigma(G * \nu)\|_{L^\infty} \leq \|\sigma\|_{C_b^1} \|G * (\mu - \nu)\|_{L^\infty} \leq \|\sigma\|_{C_b^1} \|G\|_{W^{1,2}} \|\mu - \nu\|_{W^{-1,2}},$$

and similarly

$$\|\nabla(G * \mu) \cdot (\nabla\sigma)(G * \mu) - \nabla(G * \nu) \cdot (\nabla\sigma)(G * \nu)\|_{L^\infty}$$
$$\leq \|\sigma\|_{C_b^1} \|\nabla(G * \mu) - \nabla(G * \nu)\|_{L^\infty} + \|\nabla(G * \nu)\|_{L^\infty} \|\sigma\|_{C_b^2} \|G * (\mu - \nu)\|_{L^\infty}$$
$$\leq \|\sigma\|_{C_b^1} \|G\|_{W^{2,2}} \|\mu - \nu\|_{W^{-1,2}} + \|G\|_{W^{2,2}}^2 \|\nu\|_{W^{-1,2}} \|\sigma\|_{C_b^2} \|\mu - \nu\|_{W^{-1,2}}.$$

We thus conclude that (2.2) holds, with

$$\|\Sigma[\mu] - \Sigma[\nu]\|_{W^{1,\infty}} \lesssim_{G,\sigma} (\|\nu\|_{W^{-1,2}} \vee 1)\|\mu - \nu\|_{W^{-1,2}}.$$

The similar bounds (2.1) for $b[\cdot]$ and (2.3) for $a[\cdot]$ can be proven with identical arguments, since $\|\sigma \cdot \sigma^t\|_{C_b^2} \lesssim \|\sigma\|_{C_b^2}^2$, while (2.4) follows from the definition of $W^{1+\varepsilon,\infty}$, since

$$\|\nabla[a(G * \mu)]\|_{L^\infty} \leq \|a\|_{C_b^1} \|G\|_{W^{2,\infty}} \|\mu\|_{W^{-1,2}},$$

and via similar argument as above,

$$\sup_{x \neq y} \left\{ \frac{|\nabla[a(G * \mu)](x) - \nabla[a(G * \mu)](y)|}{|x - y|^\varepsilon} \right\} \lesssim \|a\|_{C_b^2} (1 + \|G\|_{W^{2+\varepsilon,\infty}}^2) \|\mu\|_{W^{-1,2}}.$$

(2.5) and (2.6) are straightforward to check.

To lighten the notation below, we decouple the drift part by denoting

$$A[\nu, u] := \nabla \cdot (a[\nu]\nabla u) + \nabla \cdot (\tilde{b}[\nu]u), \tag{2.11}$$

where $\tilde{b}$ is the modified drift given in (1.13). The following lemma on the coercivity of the elliptic operators is proven in Appendix A. Here we fix some $v \in W^{-1,2}(\mathbb{T}^d)$, and we consider the Gelfand triple[2.1]

$$V \subset H \subset V^* \quad \text{with} \quad V := L^2(\mathbb{T}^d), \ H := W^{-1,2}(\mathbb{T}^d).$$

**Lemma 2.4.** *The continuous linear map $u \mapsto \nabla \cdot (a[v]\nabla u)$ from $W^{1,2}$ to $(L^2)^*$ extends to a continuous linear map from $L^2$ to $(L^2)^*$ such that*

$$-{}_{(L^2)^*}\langle \nabla \cdot (a[v]\nabla u), u\rangle_{L^2} + \alpha(\|a[v]\|_{W^{1+\varepsilon,\infty}}^2)\|u\|_{W^{-1,2}}^2 \geq \frac{a_{\min}}{4}\|u\|_{L^2}^2, \qquad (2.12)$$

*for $a_{\min}$ in (2.5) and some affine function $\alpha$ with positive coefficients.*

As a direct consequence of Lemma 2.12, we may consider the drift term $A[u, u]$ as defined from $L^2$ to $(L^2)^*$, where the continuity of the drift term follows directly from the regularity that we impose on $b[\cdot]$ and $\Sigma[\cdot]$ in Assumption 2.1.

## 2.1 Global Existence and Uniqueness of Tamed Equations

As mentioned above, it turns out that it is more transparent to investigate the equation as [DKP24] in the Gelfand triple $V \subset H \subset V^*$ due to the gradient-multiplicative structure of the noise. The proof of the global existence and uniqueness of (1.12) is then arranged as follows: We first introduce the $W^{-1,2}$ taming to the equation to compensate for the quadratic growth while keeping the noise term:

$$\begin{aligned}
A^\chi[u] &:= A[u, u] - \chi(\|u\|_{W^{-1,2}}^2)u, \\
B^N[u] \cdot dW_t &:= \frac{\sqrt{2}}{\sqrt{N}} \sum_{k \in \mathbb{Z}^d} \theta_k^N \nabla \cdot (f^N(u)\Sigma[u]e_k \cdot dB_t^k),
\end{aligned} \qquad (2.13)$$

where some smooth function $\chi$ is fixed to be the taming functions on $\mathbb{R}_{\geq 0}$ s.t.

$$\begin{cases} \chi(r) = 0, & r \leq M_\chi - 1, \\ \chi(r) = C_\chi(r + 1 - M_\chi), & r \geq M_\chi, \end{cases} \qquad (2.14)$$

and $\chi'$ is non-decreasing on $\mathbb{R}_{\geq 0}$, where $C_\chi > 1$ will be chosen to be large enough in Lemma 2.5. The following assumption is essential for the SPDE modelling:

$$\frac{\|\theta^N\|_{\ell^2}^2}{N\delta_N} \leq \frac{a_{\min}}{32 C_\Sigma^2}, \qquad (2.15)$$

where the constant $C_\Sigma$ is defined in (2.5), and (2.15) will be determined in (2.23). Heuristically, such a bound indicates that the coercivity of the diffusion operator in (1.12) dominates the effect of the noise. We refer to (2.15) as the *stochastic parabolicity* condition of the SPDE.

The tamed version of (1.12) is well suited for the variational framework of [RSZ24, LR15], if we use the Gelfand triple $(L^2)^* \subset W^{-1,2} \subset L^2$. We then obtain a $C([0, T], L^2)$-valued solution by deriving (taming-dependent) $L^2$ energy estimates with the same argument as in [DKP24]. The uniqueness of the $L^2$ solution follows straightforwardly from that of the $W^{-1,2}$ solution, again for fixed taming.

---

2.1. In the Gelfand triple $W^{1,2} \subset L^2 \subset W^{-1,2}$ that we consider later, the coercivity bound directly holds with the multiplicaitive coefficient $a_{\min}$.

We further exploit the fact that the $L^1$ norm of the solution to the tamed equation is non-increasing since the equation without taming preserves probability measures, while the taming terms have a negative sign. This allows us to upgrade the $W^{-1,2}$ energy estimate of the tamed equation so that it is uniform in the taming by exploiting the fact that the interactions are controlled by the $L^1$ norm of the solution according to (2.6). We use this estimate to remove the taming, which gives the global existence and uniqueness that we desire in Corollary 2.10.

We remark that, in light of the recent advances for SPDEs with locally monotone coefficients - [RSZ24, LR15], the well-posedness of (1.12) is not surprising. However, we are not aware of any result which is directly applicable in our setting. Incorporating the $W^{-1,2}$ taming into the *variational approach* seems to be the most compact way to derive the results we require.

**Lemma 2.5.** *Assume that the stochastic-parabolicity condition (2.15) holds true for $C_\Sigma$ small enough only depending on $\Sigma[\cdot]$. Given $\chi$ as in (2.14) with sufficiently large $C_\chi > 1$ only depending on $b[\cdot]$ and $\Sigma[\cdot]$ and any finite $M_\chi > 1$, together with any initial datum in $W^{-1,2}$, there exists a unique $W^{-1,2}$ solution to the SPDE given by (2.13), that is*

$$\begin{cases} du_t = A^\chi[u_t]dt + B^N[u_t] \cdot dW_t, \\ u(0,\cdot) \in W^{-1,2}. \end{cases} \tag{2.16}$$

We denote the unique solution to the tamed equation for any fixed $N \in \mathbb{N}$ and taming function $\chi$ as $m^{N,\chi}$.

**Proof.** Thanks to [RSZ24], Theorem 3.2, we only need to verify that the conditions (H1-H5) therein hold. The hemicontinuity condition (H1) is apparently true, and the coercivity bound (H2) holds for the drift term, in the sense that

$$_{V^*}\langle A^\chi[u], u\rangle_V + \frac{a_{\min}}{8}\|u\|_V^2 \lesssim_{b,\Sigma} (M_\chi C_\chi + 1)\|u\|_H^2. \tag{2.17}$$

Indeed, thanks to Lemma 2.4 and the conditions we impose on $\Sigma[\cdot]$ and $b[\cdot]$,

$$\begin{aligned}\langle A[u,u] - \chi(\|u\|_H^2)u, u\rangle &\leq -\frac{a_{\min}}{4}\|u\|_{L^2}^2 - \chi(\|u\|_{W^{-1,2}}^2)\|u\|_{W^{-1,2}}^2 \\ &\quad + \alpha(\|a[u]\|_{W^{1,\infty}}^2)\|u\|_{W^{-1,2}}^2 + \|\nabla \cdot (\tilde{b}[u]u)\|_{W^{-1,2}}\|u\|_{W^{-1,2}} \\ &\leq -\frac{a_{\min}}{4}\|u\|_{L^2}^2 + \|\tilde{b}[u]\|_{L^\infty}\|u\|_{L^2}\|u\|_{W^{-1,2}} \\ &\quad + \alpha(\|a[u]\|_{W^{1,\infty}}^2)\|u\|_{W^{-1,2}}^2 - \chi(\|u\|_{W^{-1,2}}^2)\|u\|_{W^{-1,2}}^2 \\ &\stackrel{(2.8)}{\leq} -\frac{a_{\min}}{8}\|u\|_V^2 + \left[\frac{2}{\sigma_{\text{ellip}}}\|\tilde{b}[u]\|_{L^\infty}^2 + \alpha\left(C'_\Sigma(\|u\|_H^2 + 1)\right) - \chi(\|u\|_H^2)\right]\|u\|_H^2,\end{aligned}$$

for any $\kappa > 0$ and some $C'_\Sigma > 0$. We've used the uniformly elliptic bound on $A[u,\cdot]$ in the first inequality, which is allowed since $u \in V = L^2$, and apply Young's product inequality for the last one. Since we have the control on the drift by (2.7) and (2.8) that

$$\|\tilde{b}[u]\|_{L^\infty} \lesssim \|b[u]\|_{L^\infty} + \|\Sigma[u]\|_{W^{1,\infty}} \lesssim \|u\|_H \vee 1, \tag{2.18}$$

where the constant depends only on $b[\cdot]$ and $\Sigma[\cdot]$, we arrive at (2.17) with $C_\chi > 0$ being chosen and fixed large enough only depending on $b[\cdot]$ and $\Sigma[\cdot]$, since the bound on

$$\frac{2}{a_{\min}}\|\tilde{b}[u]\|_{L^\infty}^2 + \alpha\left(C'_\Sigma(\|u\|_H^2 + 1)\right)$$

is affine with respect to $\|u\|_H^2$ depending only on $b[\cdot]$ and $\Sigma[\cdot]$, which are compensated by $\chi(\|u\|_H^2)$ when $\|u\|_H^2 \geq M_\chi$. More precisely, by distinguishing the case of $\|u\|_H^2 \geq M_\chi$ from $\|u\|_H^2 < M_\chi$, we bound

$$\alpha\left(C'_\Sigma(\|u\|_H^2 + 1)\right) + \frac{2}{a_{\min}}\|\tilde{b}[u]\|_{L^\infty}^2 - \chi(\|u\|_H^2) \lesssim M_\chi C_\chi + 1$$

due to the very construction (2.14) of $\chi$. From now on, we fix $C_\chi$ that we have chosen and an arbitrary $M_\chi > 1$. This establishes (2.17).

Furthermore, local monotonicity (H3) holds, since (with $\theta = 0$, $\beta = 2$, $\gamma = 0$ and $\lambda = 2$ in [RSZ24], Theorem 3.2)

$$2\,{}_{V^*}\langle A^\chi[u] - A^\chi[v], u - v\rangle_V + \|B^N[u] - B^N[v]\|_{L_2(\ell^2, H)}^2 \lesssim [\rho(u) + \eta(v)]\|u - v\|_H^2 \qquad (2.19)$$

where the constant depends on $b[\cdot]$, $\Sigma[\cdot]$, $\chi$ and $N$, and $\rho$, $\eta$ are measurable on $V$ s.t.

$$|\rho(u)| := \|u\|_H^2,$$
$$|\eta(v)| := (1 + \|v\|_H^2)\|v\|_V^2 + 1 + \|v\|_H^2$$

following the notation of [RSZ24]. To prove the assertion, we first split the term $\langle A^\chi[u] - A^\chi[v], u - v\rangle$ as

$$\begin{aligned}(\mathrm{I}) + (\mathrm{II}) &:= {}_{V^*}\langle A[u, u - v] - \chi(1 + \|u\|_H^2)(u - v), u - v\rangle_V \\ &\quad + {}_{V^*}\langle A[u, v] - A[v, v] - [\chi(1 + \|u\|_H^2) - \chi(1 + \|v\|_H^2)]\,v, u - v\rangle_V,\end{aligned}$$

where we recall from (2.11) that $A[\cdot, \cdot]$ is only linear in the second argument. The first term is controlled with the same argument that we have used for (2.17) by

$$(\mathrm{I}) + \frac{a_{\min}}{8}\|u - v\|_V^2 \lesssim \|u - v\|_H^2. \qquad (2.20)$$

For the second term we have

$$(\mathrm{II}) \leq {}_{V^*}\langle A[u, v] - A[v, v], u - v\rangle_V + C_\chi \big|\|u\|_H^2 - \|v\|_H^2\big| \cdot \|u - v\|_H \|v\|_H,$$

where we use the bound $\|\chi'\|_{L^\infty} \leq C_\chi$, and we estimate the first term above by

$$\begin{aligned}\langle A[u, v] - A[v, v], u - v\rangle &\leq \|a[u] - a[v]\|_{W^{1,\infty}}\|v\|_{L^2}\|u - v\|_{L^2} \\ &\quad + \|\tilde{b}[u] - \tilde{b}[v]\|_{L^\infty}\|v\|_{L^2}\|u - v\|_{W^{-1,2}} \\ &\lesssim (1 + \|v\|_H)\|v\|_V\left(\|u - v\|_H \|u - v\|_V + \|u - v\|_H^2\right),\end{aligned}$$

according to Assumption 2.1, where we recall the definition of $\tilde{b}$ from (1.13), which yields

$$\|\tilde{b}[u] - \tilde{b}[v]\|_{L^\infty} \lesssim (1 + \|v\|_H)\|u - v\|_H$$

due to (2.1) and (2.3). We then apply Young's inequality and rearrange the terms to get

$$\langle A[u,v] - A[v,v], u-v\rangle - \frac{a_{\min}\|u-v\|_V^2}{16} \lesssim (1+\|v\|_H^2)\|v\|_V^2\|u-v\|_H^2 + (1+\|v\|_H)\|v\|_V\|u-v\|_H^2.$$

Together the estimates lead to

$$(\text{II}) - \frac{a_{\min}\|u-v\|_V^2}{16} \lesssim (1+\|v\|_H^2)\|v\|_V^2\|u-v\|_H^2$$
$$+ (1+\|v\|_H)\|v\|_V\|u-v\|_H^2 + (\|u\|_H + \|v\|_H)\|v\|_H\|u-v\|_H^2. \quad (2.21)$$

On the other hand, the Hilbert-Schmidt norm of the noise term is estimated by

$$\|B^N[u] - B^N[v]\|_{L_2}^2 - \frac{2C_\Sigma^2\|\theta^N\|_{\ell^2}^2}{N\delta_N}\|u-v\|_V^2 \lesssim \|v\|_V^2\|u-v\|_H^2, \quad (2.22)$$

since

$$\|f^N(u)\Sigma[u]e_k - f^N(v)\Sigma[v]e_k\|_{L^2} \leq \|f^N(u) - f^N(v)\|_{L^2}\|\Sigma[u]\|_{L^\infty} + \|f^N(v)\|_{L^2}\|\Sigma[u] - \Sigma[v]\|_{L^\infty},$$

and therefore

$$\|f^N(u)\Sigma[u]e_k - f^N(v)\Sigma[v]e_k\|_{L^2}^2 - \frac{2C_\Sigma^2}{\delta_N}\|u-v\|_V^2 \lesssim_\Sigma \|v\|_V\|u-v\|_H,$$

uniformly in $k\in\mathbb{Z}^d$. With all the above estimations (2.20) - (2.22), we therefore choose the coefficient $C_\Sigma$ in (2.15) as such that

$$(2.15) \Longrightarrow \frac{a_{\min}}{16} > \frac{2C_\Sigma^2\|\theta^N\|_{\ell^2}^2}{N\delta_N}, \quad (2.23)$$

and arrive at (2.19). Lastly, we have for (H5) that there exists

$$\|B^N(t,u)\|_{L_2}^2 \lesssim \|u\|_V^2, \quad (2.24)$$

with similar proof to (2.22). The growth condition (H4) is verified similarly. Since the coercivity condition (3.7) in [RSZ24], Theorem 3.2 is indicated directly by (2.15), we have obtained the global existence and uniqueness of the equation given by (2.13) in $H = W^{-1,2}(\mathbb{T}^d)$. □

We recall that we have found and fixed $C_\chi > 0$ in (2.14) according to the proof of Lemma 2.5. Our next goal is to improve the regularity of the solution constructed in Lemma 2.5, for which we follow the strategy from [DKP24, Proposition 1] of proving a priori estimates for Galerkin approximations. We introduce therefore the notation for the Fourier cutoff

$$\Pi_{\leq R} := 1_{|\cdot|\leq R}(D).$$

**Lemma 2.6.** *Assume the same conditions as in Lemma 2.5, and consider the Galerkin approximation of the tamed equation*

$$\begin{cases} du_t = A_R^\chi[u_t]dt + B_R^N[u_t]\cdot dW_t, \\ u(0,\cdot) \in \Pi_{\leq N}L^2, \end{cases} \quad (2.25)$$

*where*

$$A_R^\chi[u] := \Pi_{\leq R}A^\chi[\Pi_{\leq R}u],$$
$$B_R[u]\cdot dW_t^N := \Pi_{\leq R}(B_R^N[\Pi_{\leq R}u]\cdot dW_t^N).$$

*There exists globally a unique solution to the SDE (2.25), which we denote as $m^{N,\chi,R}$, and it satisfies the estimates*

$$\mathbf{E}\|m_t^{N,\chi,R}\|_{L^2}^2 + \mathbf{E}\int_0^t \|\nabla m_s^{N,\chi,R}\|_{L^2}^2 ds \lesssim \|m_0^{N,\chi,R}\|_{L^2}^2 \exp(C_{N,\chi} t), \qquad (2.26)$$

*for any $t \geq 0$, where the constants are independent of $R > 0$.*

**Proof.** Since for any finite $R > 0$, $\|\cdot\|_{W^{-1,2}}$ is equivalent to the Euclidean norm on $\Pi_{\leq R} L^2(\mathbb{T}^d)$, the local weak monotonicity and weak coercivity of (2.25) follows directly from the monotonicity and coercivity conditions that we have verified in Lemma 2.5. Theorem 3.1.1 in [LR15] is therefore in force, and the statement of Lemma 2.6 follows. □

Therefore, if we fix the initial condition $m_0^{N,\chi,R} = m_0$ in $L^2(\mathbb{T}^d)$, the solutions to the Galerkin-truncated equations are uniformly bounded in the reflexive space $L^2(\Omega \times [0,T], W^{1,2}(\mathbb{T}^d))$. As mentioned above, global existence and uniqueness of $L^2$ solution to the tamed equation, which constitutes the first part of the following corollary, follows then the same line as [DKP24], Theorem 3.10. We thus omit the details and refer the readers to the proof therein.

The solutions are, however, still taming-dependent, and we now invoke the $L^1$ conditions (2.6) in Assumption 2.1 and exploit further the fact that the solutions to any of the tamed equations have non-increasing $L^1$ norms. Intuitively, the noise is 'turned off' at the places where the density field is zero, and therefore the deterministic part of the equation dominates.

**Corollary 2.7.** *Under the same conditions as in Lemma 2.5, there exists for any $T > 0$ a unique $L^2$ solution $m^{N,\chi}$ to the tamed SPDE given by (2.13) in the space*

$$C([0,T], L^2(\mathbb{T}^d)) \cap L^2([0,T], W^{1,2}(\mathbb{T}^d)).$$

*Besides, the solution $(m_t^{N,\chi})_{t \in [0,T]}$ stays non-negative and the $L^1$ norm is non-increasing almost surely, given any non-negative initial data*

$$0 \leq m^{N,\chi}(0) = m_0^N \in L^2(\mathbb{T}^d).$$

*Consequently, we have the taming-independent energy bound for non-negative initial data: For any $t \in [0,T]$, there exists some constant $C_N > 0$ independent of $\chi$ and $t$, s.t.*

$$\mathbf{E}\sup_{s \in [0,t]} \|m_s^{N,\chi}\|_{L^2}^2 \leq \|m_0^N\|_{L^2}^2 \exp(C_N t). \qquad (2.27)$$

**Proof.** The first part of Corollary 2.7 has already been addressed above. We assume thereofre that the initial condition is non-negative in $L^2(\mathbb{T}^d)$ and denote

$$K^\zeta := \kappa_\zeta * (\cdot)_- \in C^\infty(\mathbb{T}^d, [0,\infty)), \quad \text{with} \quad \kappa_\zeta := \frac{1}{\zeta}\kappa(\frac{\cdot}{\zeta}) \quad \text{for} \quad \zeta > 0,$$

with a nonnegative and smooth convolution kernel $\kappa$, for which we compute

$$(K^\zeta)' = -\kappa_\zeta * 1_{\leq 0}, \quad (K^\zeta)'' = \kappa_\zeta,$$

in the distributional sense. We apply the Krylov-Itô formula ([Kry13], Chapter 2) to the process $t \mapsto \int_{\mathbb{T}^d} K^\zeta(m_t^{N,\chi})$, which approximates $t \mapsto \int_{\mathbb{T}^d} (m_t^{N,\chi})_-$, that is

$$\int_{\mathbb{T}^d} K^\zeta(m_t^{N,\chi}) + \int_0^t \left\{ \int_{\mathbb{T}^d} \kappa_\zeta(m_s^{N,\chi}) \nabla m_s^{N,\chi} \cdot (a[m_s^{N,\chi}] \nabla m_s^{N,\chi}) \right\} ds$$
$$- \int_0^t \left\{ \int_{\mathbb{T}^d} (\kappa_\zeta * 1_{\leq 0})(m_s^{N,\chi}) m_s^{N,\chi} \chi(\|m_s^{N,\chi}\|_{W^{-1,2}}^2) \right\} ds$$
$$= \int_{\mathbb{T}^d} K^\zeta(m_0^{N,\chi}) - \int_0^t \left\{ \int_{\mathbb{T}^d} \kappa_\zeta(m_s^{N,\chi}) \nabla m_s^{N,\chi} \cdot \tilde{b}[m_s^{N,\chi}] m_s^{N,\chi} \right\} ds$$
$$+ \frac{1}{N} \int_0^t \sum_{k \in \mathbb{Z}^d} |\theta_k^N|^2 \left\{ \int_{\mathbb{T}^d} \kappa_\zeta(m_s^{N,\chi}) |\nabla \cdot [f^N(m_s^{N,\chi}) \Sigma[m_s^{N,\chi}] e_k]|^2 \right\} ds + M_t^{\zeta,N,\chi},$$

where

$$M_t^{\zeta,N,\chi} := -\frac{2}{\sqrt{N}} \int_0^t \int_{\mathbb{T}^d} \kappa_\zeta(m_s^{N,\chi}) f^N(m_s^{N,\chi}) \nabla m_s^{N,\chi} \cdot (\Sigma[m_s^{N,\chi}] dW_s^N).$$

Since by assumption $(m_t^{N,\chi})_{t \in [0,T]}$ is continuous in $L^2(\mathbb{T}^d)$, the process $t \mapsto \|\Sigma[m_t^{N,\chi}]\|_{W^{1,\infty}}$ is almost surely locally bounded due to the assumption on $\Sigma[\cdot]$. This implies in particular that $M^{\zeta,N,\chi}$ in the Krylov-Itô formula above is indeed a martingale, since

$$\mathbf{E}[M^{\zeta,N,\chi}]_T \lesssim_{N,\zeta} \sum_{k \in \mathbb{Z}^d} |\theta_k^N|^2 \mathbf{E} \int_0^T \|\nabla m_t^{N,\chi}\|_{L^2}^2 \|m_t^{N,\chi}\|_{L^1}^2 dt \lesssim 1.$$

By further controlling the correction term according to (2.5) and (2.8) by the bound

$$\int_{\mathbb{T}^d} \kappa_\zeta(m_s^{N,\chi}) |\nabla \cdot [f^N(m_s^{N,\chi}) \Sigma[m_s^{N,\chi}] e_k]|^2$$
$$\lesssim \int_{\mathbb{T}^d} \kappa_\zeta(m_s^{N,\chi}) |\nabla f^N(m_s^{N,\chi})|^2 + (\|m_s^{N,\chi}\|_{L^2}^2 + |k|^2) \int_{\mathbb{T}^d} \kappa_\zeta(m_s^{N,\chi}) |f^N(m_s^{N,\chi})|^2$$
$$\overset{(1.14)}{\lesssim} \delta_N^{-1} \int_{\mathbb{T}^d} \kappa_\zeta(m_s^{N,\chi}) \|\nabla m_s^{N,\chi}\|^2 + (\|m_s^{N,\chi}\|_{L^2}^2 + |k|^2) \int_{\mathbb{T}^d} \kappa_\zeta(m_s^{N,\chi}) |m_s^{N,\chi}|^2,$$

and then using the same argument as in the proof of Lemma 2.5, the $L^\infty$ bound on $\tilde{b}[\cdot]$ and the parabolicity condition (2.15), we can absorb the first term above and bound (with constant depending on $N$, $b[\cdot]$ and $\Sigma[\cdot]$)

$$\mathbf{E} \int_{\mathbb{T}^d} K^\zeta(m_t^{N,\chi}) - \mathbf{E} \int_{\mathbb{T}^d} K^\zeta(m_0^{N,\chi}) + \frac{a_{\min}}{2} \mathbf{E} \int_0^t \left\{ \int_{\mathbb{T}^d} \kappa_\zeta(m_s^{N,\chi}) |\nabla m_s^{N,\chi}|^2 \right\} ds$$
$$- \mathbf{E} \int_0^t \left\{ \int_{\mathbb{T}^d} \kappa_\zeta * 1_{\leq 0}(m_s^{N,\chi}) m_s^{N,\chi} \chi(\|m_s^{N,\chi}\|_{W^{-1,2}}^2) \right\} ds$$
$$\lesssim \mathbf{E} \int_0^t (\|m_s^{N,\chi}\|_{L^2}^2 + 1) \left\{ \int_{\mathbb{T}^d} \kappa_\zeta(m_s^{N,\chi}) |m_s^{N,\chi}|^2 + |m_s^{N,\chi} \nabla m_s^{N,\chi}| \right\} ds. \tag{2.28}$$

By sending $\zeta \to 0^+$, the R.H.S. of (2.27) vanishes according to Lemma A.2 and dominated convergence, since

$$\int_0^t (\|m_s^{N,\chi}\|_{L^2}^2 + 1) \left\{ \int_{\mathbb{T}^d} \kappa_\zeta(m_s^{N,\chi}) (|m_s^{N,\chi}|^2 + |m_s^{N,\chi} \nabla m_s^{N,\chi}|) \right\} ds$$
$$\lesssim \int_0^t (\|m_s^{N,\chi}\|_{L^2}^2 + 1) \left\{ \int_{\mathbb{T}^d} |m_s^{N,\chi}| + |\nabla m_s^{N,\chi}| \right\} ds \lesssim \sup_{t \in [0,T]} \|m_s^{N,\chi}\|_{L^2}^4 + \int_0^T \|\nabla m_s^{N,\chi}\|_{L^2}^2 ds,$$

where we have used (A.3) for the first inequality, and the integrability of the upper bound follows from (A.10) and (2.26). Besides, dominated convergence also indicates that

$$\lim_{\zeta \to 0+} \mathbf{E} \int_0^t \left\{ \int_{\mathbb{T}^d} \kappa_\zeta * 1_{\leq 0}(m_s^{N,\chi}) m_s^{N,\chi} \chi(\|m_s^{N,\chi}\|_{W^{-1,2}}^2) \right\} ds$$
$$= \mathbf{E} \int_0^t \chi(\|m_s^{N,\chi}\|_{W^{-1,2}}^2) \int_{\mathbb{T}^d} -(m_s^{N,\chi})_- ds \leq 0. \tag{2.29}$$

We have therefore

$$\mathbf{E}\int_{\mathbb{T}^d}(m_t^{N,\chi})_- = \mathbf{E}\int_{\mathbb{T}^d}(m_t^{N,\chi})_- - \int_{\mathbb{T}^d}(m_0^N)_- \le 0, \qquad (2.30)$$

which shows that the solution stays non-negative for non-negative $L^2$ initial condition almost surely. By testing the tamed equation with the constant-one function on torus, we have almost surely

$$\int_{\mathbb{T}^d}|m_t^{N,\chi}| = \int_{\mathbb{T}^d}m_t^{N,\chi} \le \int_{\mathbb{T}^d}m_0^N, \qquad (2.31)$$

for all $t\in[0,T]$, since the taming term has the correct negative sign while the other terms vanish. We can now apply again the Krylov-Itô formula to the equation, now with functional $u \mapsto \frac{1}{2}\|u\|_{L^2}^2$, and derive

$$\begin{aligned}
& \frac{1}{2}\int_{\mathbb{T}^d}|m_t^N|^2 + \int_0^t\Big\{\int_{\mathbb{T}^d}\nabla m_s^N \cdot (a[m_s^N]\nabla m_s^N) + \int_{\mathbb{T}^d}|m_s^N|^2\chi(\|u\|_{W^{-1,2}}^2)\Big\}\mathrm{d}s \\
& = \frac{1}{2}\int_{\mathbb{T}^d}|m_0^N|^2 + M_t^N - \int_0^t\Big\{\int_{\mathbb{T}^d}\nabla m_s^N \cdot \tilde{b}[m_s^N]m_s^N\Big\}\mathrm{d}s \\
& \quad + \frac{1}{N}\int_0^t\sum_{k\in\mathbb{Z}^d}|\theta_k^N|^2\Big\{\int_{\mathbb{T}^d}|\nabla\cdot[f^N(m_s^N)\Sigma[m_s^N]e_k]|^2\Big\}\mathrm{d}s,
\end{aligned}$$

where thanks to (2.31) and the conditions (2.6) on the coefficients, we can bound

$$\sup_{t\in[0,T]}(\|\tilde{b}[m_t^N]\|_{L^\infty}\vee\|\Sigma[m_t^N]\|_{W^{1,\infty}}) \lesssim \|m_0^N\|_{L^1}. \qquad (2.32)$$

It follows that the martingale part can be controlled by

$$\begin{aligned}
\mathbf{E}\sup_{\tau\in[0,t]}|M_\tau^N| & = \mathbf{E}\sup_{\tau\in[0,t]}\Big|\int_0^\tau\int_{\mathbb{T}^d}f^N(m_s^N)\nabla m_s^N \cdot \Sigma[m_s^N]\mathrm{d}W_s^N\Big| \\
& \lesssim_N \mathbf{E}\Big(\int_0^t\|f^N(m_s^N)\nabla m_s^N\|_{L^1}^2\Big)^{\frac{1}{2}} \le C + \frac{a_{\min}}{4}\mathbf{E}\int_0^t\|\nabla m_s^N\|_{L^2}^2,
\end{aligned}$$

for some $C = C(N,\Sigma[\cdot]) > 0$, where the BDG inequality is applied. The absolutely continuous part can be treated easily with the upgraded estimates (2.32) on noise and drift by canonical argument, so that we can conclude that

$$\mathbf{E}\sup_{s\in[0,t]}\|m_s^N\|_{L^2}^2 + \frac{a_{\min}}{4}\mathbf{E}\int_0^t\|\nabla m_s^N\|_{L^2}^2 \lesssim 1 + \int_0^t\mathbf{E}\|m_s^N\|_{L^2}^2$$

using Young inequality and parabolicity condition (2.15), where the constant depends on $N$, $\|m_0\|_{L^1}$ and the coefficients of the equation, and (2.27) follows from Grönwall inequality. $\square$

**Remark 2.8.** If the taming term is "turned off", then in the above proof of the non-negativity of the solution the expression (2.28) vanishes, while the other parts remain the same, and thus we obtain the same conclusion (2.30) that the negative part vanishes. Since the untamed equation (1.12) preserves the integral on $\mathbb{T}^d$ of the solution against a constant test function (i.e. mass), we have the $L^1$ conservation

$$\|m_t^N\|_{L^1} \equiv \|m_0\|_{L^1}, \qquad (2.33)$$

for any $t\in[0,T]$ and any non-negative initial condition $m_0\in L^2(\mathbb{T}^d)$.

## 2.2 Global Existence and Uniqueness of (1.12)

Recall that

$$m^{N,\chi} \in C([0,T], L^2(\mathbb{T}^d)) \cap L^2([0,T], W^{1,2}(\mathbb{T}^d))$$

is the unique probabilistically strong solution to the equation given by (2.16) with initial data $0 \leq m_0^{N,\chi} \in L^2(\mathbb{T}^d)$. The last ingredient to prove the global existence and uniqueness of equation (1.12) is the following a priori bound, the proof of which follows directly from Remark 2.8 and the same argument we have used for (2.27).

**Lemma 2.9.** *Let $m^N$ be a (probabilistically strong) solution to (1.12) in*

$$C([0,T], L^2(\mathbb{T}^d)) \cap L^2([0,T], W^{1,2}(\mathbb{T}^d)),$$

*with initial probability $m_0^N$ with $L^2$ density. Then for any $t \in [0,T]$, there exists some $C_N > 0$ s.t.*

$$\mathbf{E} \sup_{s \in [0,t]} \|m_s^N\|_{L^2}^2 \leq \|m_0^N\|_{L^2}^2 \exp(C_N t). \tag{2.34}$$

We are ready to prove the main result of this section, i.e. the existence and uniqueness of the probabilistically strong solution of the regularized Dean-Kawasaki equation (1.12) in $C([0,T], L^2(\mathbb{T}^d)) \cap L^2([0,T], W^{1,2}(\mathbb{T}^d))$.

**Corollary 2.10.** *We assume that the stochastic parabolicity condition (2.15) holds true. Given any initial datum $m_0^N \in L^2(\mathbb{T}^d) \cap \mathcal{P}(\mathbb{T}^d)$, the regularized Dean-Kawasaki equation (1.12) admits a unique probabilistically strong solution in*

$$C([0,T], L^2(\mathbb{T}^d)) \cap L^2([0,T], W^{1,2}(\mathbb{T}^d)).$$

**Proof.** We fix the initial datum $m_0^N \in L^2 \subset W^{-1,2}$ and a sequence of $\{\chi_i\}_{i \geq 1}$ as in (2.14) s.t. $M_{\chi_i} \to \infty$ as $i \to \infty$. To prove the existence of the solution, we first glue the solutions $\{m^{N,\chi_i}\}_{i \geq 1}$ together by setting

$$m_t^N(\omega) := m_t^{N,\chi_i}(\omega), \quad \text{for } (t,\omega) \text{ s.t. } t < \tau^i(\omega), \tag{2.35}$$

where $\tau^i$ is the exit time (before $T$) of $m^{N,\chi_i}$ from the centered ball in $W^{-1,2}$ of radius $M_{\chi_i} - 1$. By the probabilistically strong uniqueness of the solutions to the tamed equations and the fact that all $m^{N,i}$ solve the same equation before they exceed their $W^{-1,2}$-norm thresholds, $m^N$ is well defined up to time $T$, as long as we have $\tau^i \to T$ as $i \to \infty$ a.s. This indeed is the case, since

$$\mathbf{P}(\tau^i < T) \leq \mathbf{P}\left(\max_{t \in [0,T]} \|m_t^{N,i}\|_{W^{-1,2}} > M_{\chi_i} - 1\right) \leq \frac{1}{M_{\chi_i} - 1} \mathbf{E}\left(\max_{t \in [0,T]} \|m_t^{N,i}\|_{W^{-1,2}}\right) \to 0,$$

as $i \to \infty$, due to the uniform-in-$\chi$ energy estimate (2.27). By construction, $(m_t^N)_{t \in [0,T]}$ is a probabilistically strong solution to (1.12) in

$$C([0,T], L^2(\mathbb{T}^d)) \cap L^2([0,T], W^{1,2}(\mathbb{T}^d)).$$

On the other hand, if we have another solution $\tilde{m}^N$ to (1.12), we can define $\tilde{\tau}^i$ as the minimum of the exit times of $\tilde{m}^N$ and the solution $m^N$ constructed above from the centered ball in $W^{-1,2}$ of radius $M_{\chi_i} - 1$. Due to the probablistically strong uniqueness of the solutions to the tamed equations, we have

$$m_t^N 1_{t < \tilde{\tau}^i} = \tilde{m}_t^N 1_{t < \tilde{\tau}^i},$$

and by the a priori $L^2$ energy estimate in Lemma 2.9, we have a.s. $\lim_{i \to \infty} \tilde{\tau}^i = T$, and therefore the two solutions coincide. □

Corollary 2.10 now provides us with the unique solution to the equation (1.12) with initial data $m_0^{N,\epsilon}$ as in (3.1) under the assumption (2.15).

# 3 Itô's formula and Kolmogorov equations

We now focus on the approximation of the diffusive interacting particle system (1.8) by the regularized Dean-Kawasaki equation (1.12) and specify the weak error to be the difference of the expectations of the two stochastic evolutions of probability measures evaluated by functionals on $\mathcal{P}(\mathbb{T}^d)$, which will be assumed to have sufficient regularity in Assumption 4.1. A concise review on the Lions calculus of functionals on probability measures on the torus is presented in Appendix A.3.

## 3.1 Initial Data and Entropy Bound

To prepare for the weak error estimate, we start by preparing suitable initial data for the equation. We denote

$$\mu_0^N := \frac{1}{N} \sum_{i=1}^{N} \delta_{x_i}, \quad m_0^{N,\epsilon} := \rho^\epsilon * \mu_0^N, \quad \text{where} \quad \rho^\epsilon := \sum_{z \in \mathbb{Z}^d} \frac{1}{\epsilon^d} \rho\left(\frac{\cdot + z}{\epsilon}\right), \tag{3.1}$$

for $\epsilon > 0$. $\rho^\epsilon$ is the rescaled periodization of the kernel $\rho$ with $\rho \in C_c^\infty(\mathbb{R}^d)$ being radially symmetric, normalized and positive.

**Lemma 3.1.** *The empirical measure $\mu_0^N$ is approximated by $m_0^{N,\epsilon}$ in the following sense:*

$$\mathcal{W}_2(\mu_0^N, m_0^{N,\epsilon}) \leq \epsilon \left( \int_{\mathbb{R}^d} |y|^2 \rho(y) \mathrm{d}y \right)^{\frac{1}{2}}. \tag{3.2}$$

*Further, we have the entropy bound*

$$\int_{\mathbb{T}^d} m_0^{N,\epsilon} \log(m_0^{N,\epsilon}) \leq d \log\left(\frac{1}{\epsilon}\right) + \log(\|\rho\|_{L^\infty}). \tag{3.3}$$

**Proof.** Let $(X, Y)$ be a realization of $\mu_0^N \otimes \rho^\epsilon$. Then $X + Y$ has the distribution of $m_0^{N,\epsilon}$ and therefore

$$\mathcal{W}_2^2(\mu_0^N, m_0^{N,\epsilon}) \leq \mathbf{E} \mathrm{d}_{\mathbb{T}^d}^2(X + Y, X) \leq \int_{\mathbb{R}^d} \frac{|x|^2}{|\epsilon|^d} \rho\left(\frac{x}{\epsilon}\right) \mathrm{d}x,$$

which establishes (3.2). The entropy bound is the same as in [DKP24], Lemma 4.3. □

We now prove an entropy estimate, which will be used to control the Fisher information that appears in the approximation of the formal generator in Lemma 3.7.

**Lemma 3.2.** *Given $(m_s^{N,\epsilon})_{s \in [0,T]}$ the solution in the sense of Corollary 2.10 to the regularized Dean-Kawasaki equation (1.12) and the initial datum $m_0^{N,\epsilon}$ as in (3.1), we have*

$$\mathbf{E} \int_0^T \int_{\mathbb{T}^d} \frac{|\nabla m_s^{N,\epsilon}|^2}{m_s^{N,\epsilon}} ds \lesssim 1 + \log\left(\frac{1}{\epsilon}\right) + \frac{\|(k\theta_k^N)\|_{\ell_k^2}^2}{N}, \qquad (3.4)$$

*for some constant independent of $N$ and $\epsilon$.*

**Proof.** For any $\zeta > 0$, we apply Krylov-Itô formula to $(m_s^{N,\epsilon})_{s \in [0,t]}$ and derive

$$\begin{aligned}
d \int_{\mathbb{T}^d} (\zeta + m_s^{N,\epsilon}) \log(\zeta + m_s^{N,\epsilon}) &= -\left[\int_{\mathbb{T}^d} \frac{(a[m_s^{N,\epsilon}]\nabla m_s^{N,\epsilon}) \cdot \nabla m_s^{N,\epsilon}}{\zeta + m_s^{N,\epsilon}}\right] ds \\
&\quad - \left[\int_{\mathbb{T}^d} \frac{m_s^{N,\epsilon}}{\zeta + m_s^{N,\epsilon}} \tilde{b}[m_s^{N,\epsilon}] \cdot \nabla m_s^{N,\epsilon}\right] ds \\
&\quad + \frac{1}{N} \sum_{k \in \mathbb{Z}^d} |\theta_k^N|^2 \left[\int_{\mathbb{T}^d} \frac{|\nabla \cdot (f^N(m_s^{N,\epsilon})\Sigma[m_s^{N,\epsilon}]e_k)|^2}{\zeta + m_s^{N,\epsilon}}\right] ds + dM_s^N,
\end{aligned}$$

where $M^N$ is again genuinely a martingale that vanishes after we take expectation. We absorb again the second term on the R.H.S. by the first one by

$$\left|\int_{\mathbb{T}^d} \frac{m_s^{N,\epsilon}}{\zeta + m_s^{N,\epsilon}} b[m_s^{N,\epsilon}] \cdot \nabla m_s^{N,\epsilon}\right| \leq \int_{\mathbb{T}^d} (m_s^{N,\epsilon})^{\frac{1}{2}} |\tilde{b}[m_s^{N,\epsilon}]| \frac{|\nabla m_s^{N,\epsilon}|}{(m_s^{N,\epsilon})^{\frac{1}{2}}}$$

$$\overset{2.1,(2.33)}{\lesssim} \left(\int_{\mathbb{T}^d} \frac{|\nabla m_s^{N,\epsilon}|^2}{m_s^{N,\epsilon}}\right)^{\frac{1}{2}} \leq C + \frac{a_{\min}}{2} \int_{\mathbb{T}^d} \frac{|\nabla m_s^{N,\epsilon}|^2}{m_s^{N,\epsilon}},$$

for some generic $C > 0$, while the third term on the right can be treated identically as in [DKP24], Proposition 4 by bounding

$$\int_{\mathbb{T}^d} \frac{|\nabla \cdot (f^N(m^{N,\epsilon})\Sigma[m_s^{N,\epsilon}]e_k)|^2}{\zeta + m_s^{N,\epsilon}} \leq 2C_\Sigma^2 |k|^2 + 2\|\Sigma[m_s^N]\|_{W^{1,\infty}}^2 + 2C_\Sigma^2 \|(f^N)'\|_{L^\infty}^2 \int_{\mathbb{T}^d} \frac{|\nabla m_s^{N,\epsilon}|^2}{\zeta + m_s^{N,\epsilon}},$$

where we recall that $C_\Sigma$ is defined in (2.5). Combining the estimates above with (2.15) and (2.23), we arrive at

$$\mathbf{E} \int_{\mathbb{T}^d} (\zeta + m_\cdot^{N,\epsilon}) \log(\zeta + m_\cdot^{N,\epsilon}) \Big|_0^T + \frac{a_{\min}}{4} \mathbf{E} \int_0^T \int_{\mathbb{T}^d} \frac{|\nabla m_s^{N,\epsilon}|^2}{m_s^{N,\epsilon}} ds \lesssim_\Sigma \frac{\|k\theta_k^N\|_{\ell_k^2}^2}{N},$$

and the desired bound (3.4) by using (3.3), Fatou Lemma and dominated convergence to send $\zeta$ to $0^+$. The detailed limit argument can be found in [DKP24]. $\square$

## 3.2 Itô's Formula

To start our investigation on the weak error, we derive the Itô's formulas for the interacting diffusions and the regularized SPDEs for cylindrical fucntionals and then extend them to general functionals by the approximation argument. Unlike in [KLvR19], the approximation is constructed with *Fourier multipliers* instead of Bernstein polynomials, due to the fact that the latter seems not well-suited for periodic boundary conditions. For this reason, we first recall the Fejér kernels $(\varphi^M)_{M \in \mathbb{N}}$ to be

$$\hat{\varphi}^M(k) := \begin{cases} \prod_{j=1}^d \left(1 - \frac{|k_j|}{M+1}\right), & |k|_\infty \leq M, \\ 0, & \text{else}, \end{cases} \qquad (3.5)$$

which is a positive approximation of the identity according to [Gra08]. We consider the following mollification for any given functional $\mathcal{F}[\cdot]$ on $\mathbb{T}^d$ by inner regularization

$$\mathcal{F}^M[\mu] := \mathcal{F}[\mu * \varphi^M]. \tag{3.6}$$

Notice that while it is relatively straightforward to derive the convergence of $\mathcal{F}^M[\cdot]$ to $\mathcal{F}[\cdot]$, the family $(\mathcal{F}^M[\,\cdot\,])_{M \in \mathbb{N}}$ does not consist of cylindrical functionals, since it seems difficult to derive the smoothness of the dependence of $\mathcal{F}[\mu * \varphi^M]$ on the Fourier coefficients of $\mu$ up to modes $|k|_\infty \leq M$.

**Lemma 3.3.** *Given a functional $\mathcal{F}: \mathcal{P}_2 \to \mathbb{R}$ with second-order interior derivatives in the sense of (A.15), s.t.*

$$D_\mu \mathcal{F}: \mathcal{P}(\mathbb{T}^d) \times \mathbb{T}^d \to R^d, \quad D_\mu^2 \mathcal{F}: \mathcal{P}(\mathbb{T}^d) \times \mathbb{T}^{2d} \to \mathbb{R}^{d \times d}$$

*are both jointly Lipschitz, we have that $\mathcal{F}^M$ approximates $\mathcal{F}$ in the sense that*

$$D_\mu^k \mathcal{F}^M \to D_\mu^k \mathcal{F}, \tag{3.7}$$

*uniformly on $\mathcal{P}(\mathbb{T}^d) \times \mathbb{T}^{d \times k}$, for $k = 0, 1, 2$.*

**Proof.** By the same argument as used in (3.2) and the Lipschitz continuity, we have

$$\sup_{\mu \in \mathcal{P}(\mathbb{T}^d)} |\mathcal{F}^M[\mu] - \mathcal{F}[\mu]| \leq \|\mathcal{F}\|_{\text{lip}} \sup_{\mu \in \mathcal{P}(\mathbb{T}^d)} \mathcal{W}_2(\mu, \mu * \varphi^M)$$

$$\leq \|\mathcal{F}\|_{\text{lip}} \int_{\mathbb{T}^d} d_{\mathbb{T}^d}^2(x, 0) \varphi^M(x) dx \to 0,$$

according to (A.16) as $M \to \infty$, since $\varphi^M \rightharpoonup \delta_0$.

For $k = 1$ in (3.7), it follows directly from the definitions (A.13) and (A.14) that

$$\frac{\delta \mathcal{F}^M}{\delta \mu}[\mu] = \frac{\delta \mathcal{F}}{\delta \mu}[\mu * \varphi^M] * \varphi^M,$$

and

$$D_\mu \mathcal{F}^M[\mu] = D_\mu \mathcal{F}^M[\mu * \varphi^M] * \varphi^M. \tag{3.8}$$

It thus follows from the triangle inequality and the same argument as above that

$$\sup_{\mu \in \mathcal{P}(\mathbb{T}^d), x \in \mathbb{T}^d} |D_\mu \mathcal{F}^M(\mu, x) - D_\mu \mathcal{F}(\mu, x)| \leq \|D_\mu \mathcal{F}\|_{\text{lip}} \sup_{\mu \in \mathcal{P}(\mathbb{T}^d)} \mathcal{W}_2(\mu, \mu * \varphi^M)$$

$$+ \sup_{\mu \in \mathcal{P}(\mathbb{T}^d)} \|D_\mu \mathcal{F}[\mu * \varphi^M] - D_\mu \mathcal{F}[\mu * \varphi^M] * \varphi^M\|_{L^\infty},$$

where

$$\sup_{\mu \in \mathcal{P}(\mathbb{T}^d)} \|D_\mu \mathcal{F}[\mu * \varphi^M] - D_\mu \mathcal{F}[\mu * \varphi^M] * \varphi^M\|_{L^\infty}$$

$$\leq \sup_{\mu \in \mathcal{P}(\mathbb{T}^d), x \in \mathbb{T}^d} \left| \int_{\mathbb{T}^d} \{D_\mu \mathcal{F}[\mu](x) - D_\mu \mathcal{F}[\mu](x-y)\} \varphi^M(y) dy \right|$$

$$\leq \sup_{\mu \in \mathcal{P}(\mathbb{T}^d)} \|D_\mu \mathcal{F}[\mu]\|_{W^{1,\infty}} \int_{\mathbb{T}^d} d_{\mathbb{T}^d}(x, 0) \varphi^M(y) dy \xrightarrow{(A.18)} 0,$$

as $M \to 0$, since $D_\mu \mathcal{F}$ is jointly Lipschitz continuous. The convergence of the second order derivatives is proved similarly. □

The following approximation due to [CD22][3.1] seems to be the most direct way to obtain Itô's formula in our case. The merit is that we can choose the functions appearing in the (non-unique) representation of the cylindrical functionals (explicitly $(F^M)_{M \geq 1}$ in (3.9) below) to be smooth.

**Proposition 3.4.** *For any Lipschitz continuous $\mathcal{F}[\cdot]$ on $\mathcal{P}(\mathbb{T}^d)$, there exists a sequence of functionals $(\bar{\mathcal{F}}^M[\,\cdot\,])_{M \geq 1}$ on $\mathcal{P}(\mathbb{T}^d)$, s.t. for each $M \geq 1$,*

$$\bar{\mathcal{F}}^M[\mu] = F^M(\{\langle \mu, \varphi \rangle : \varphi \in \mathscr{D}_M\}) \qquad (3.9)$$

*for some smooth function $F^M : \mathbb{R}^{\mathscr{D}_M} \to \mathbb{R}^d$ with $\mathscr{D}_M$ a finite subset of the orthonormal basis*

$$\{1\} \cup \{\sqrt{2}\cos(2\pi k \cdot x), \sqrt{2}\sin(2\pi k \cdot x)\}_{k \in \mathbb{Z}_+^d}$$

*of $L^2(\mathbb{T}^d, \mathbb{R})$, and $\bar{\mathcal{F}}^M[\cdot]$ approximates $\mathcal{F}[\cdot]$ uniformly on $\mathcal{P}(\mathbb{T}^d)$. If we further assume that $\mathcal{F}[\cdot]$ has Lipschitz derivatives up to second order, then (3.7) holds.*

We adjust the original proof in [CD22] slightly for our use in the Appendix A.3.2. We are now prepared to prove the following Itô formula.

**Proposition 3.5. (Itô's formula for Interacting Particles)** *For any $\mathcal{F}[\cdot]$ on $\mathcal{P}(\mathbb{T}^d)$ with Lipschitz-continuous interior derivatives up to the second order in the sense of (A.15), Krylov-Itô's formula holds for the empirical measure $\rho^N$ of the interacting particle system (1.8):*

$$\mathcal{F}[\rho_t^N] - \mathcal{F}[\rho_0^N] = \int_0^t \int_{\mathbb{T}^d} \{a[\rho_s^N](x) : \nabla D_\mu \mathcal{F}(\rho_s^N, x) + b[\rho_s^N](x) \cdot D_\mu \mathcal{F}(\rho_s^N, x)\} \mathrm{d}\rho_s^N(x) \mathrm{d}s$$

$$+ \frac{1}{N} \int_0^t \int_{\mathbb{T}^d} a[\rho_s^N](x) : D_\mu^2 \mathcal{F}(\rho_s^N, x, x) \mathrm{d}\rho_s^N(x) \mathrm{d}s + M_t^{N,\mathcal{F}}, \qquad (3.10)$$

*where $(M_t^{N,\mathcal{F}})_{t \geq 0}$ is a continuous martingale with quadratic variation*

$$[M^{N,\mathcal{F}}]_t = \frac{2}{N} \int_0^t \int_{\mathbb{T}^d} |\Sigma^t[\rho_s^N](x) \cdot D_\mu \mathcal{F}(\rho_s^N, x)|^2 \mathrm{d}\rho_s^N(x) \mathrm{d}s,$$

**Proof.** For any $\varphi = (\varphi_i)_{i=1}^I \in C^\infty(\mathbb{T}^d, \mathbb{R}^I)$ and $f \in C^\infty(\mathbb{R}^I, \mathbb{R})$ with an arbitrary $I \in \mathbb{N}$, due to Itô-Krylov formula we have for $\mathcal{F} = f(\langle \varphi, \cdot \rangle)$ that

$$f(\langle \varphi, \rho_t^N \rangle) - f(\langle \varphi, \rho_0^N \rangle) = \int_0^t Df(\langle \varphi, \rho_s^N \rangle) \{\langle (b[\rho_s^N] \cdot \nabla + a[\rho_s^N] : \nabla^2) \varphi, \rho_s^N \rangle\} \mathrm{d}s$$

$$+ \frac{1}{N} \int_0^t D^2 f(\langle \varphi, \rho_s^N \rangle) \{\langle |\Sigma[\rho_s^N]^t \cdot \nabla \varphi|^2, \rho_s^N \rangle\} \mathrm{d}s + M_t^{\varphi,f},$$

where we abbreviate under Einstein summation that

$$Df(\langle \varphi, \rho_s^N \rangle) \{\langle b[\rho_s^N] \cdot \nabla \varphi, \rho_s^N \rangle\} := \partial^i f(\langle \varphi, \rho_s^N \rangle) \langle b_l[\rho_s^N] \cdot \partial^l \varphi_i, \rho_s^N \rangle = \langle D_\mu \mathcal{F}[\rho_s^N] \cdot b[\rho_s^N], \rho_s^N \rangle,$$

$$Df(\langle \varphi, \rho_s^N \rangle) \cdot \langle (a[\rho_s^N] : \nabla^2) \varphi, \rho_s^N \rangle := \partial^i f(\langle \varphi, \rho_s^N \rangle) \cdot \langle (a[\rho_s^N] : \nabla^2) \varphi_i, \rho_s^N \rangle$$
$$= \langle a[\rho_s^N] : \nabla D_\mu \mathcal{F}[\rho_s^N], \rho_s^N \rangle,$$

---

[3.1]. c.f. Chapter 3 therein for details

and

$$D^2 f(\langle \varphi, \rho_s^N \rangle)\{\langle |\Sigma[\rho_s^N]^t \cdot \nabla \varphi|^2, \rho_s^N \rangle\} := \partial^{ij} f(\langle \varphi, \rho_s^N \rangle)\langle (\Sigma[\rho_s^N]^t \cdot \nabla \varphi_i) \cdot (\Sigma[\rho_s^N]^t \cdot \nabla \varphi_j), \rho_s^N \rangle$$
$$= \int_{\mathbb{T}^d} a[\rho_s^N](x) : D_\mu^2 \mathcal{F}[\rho_s^N](x,x) \mathrm{d}\rho_s^N(x),$$

for cylindrical functionals $\mathcal{F} = f(\langle \varphi, \cdot \rangle)$. It follows from similar calculation that the quadratic variation of the martingale is given by

$$[M^{\varphi,f}]_t = \frac{2}{N} \int_0^t \int_{\mathbb{T}^d} |\Sigma[\rho_s^N]^t \cdot D_\mu \mathcal{F}[\rho_s^N](x)|^2 \mathrm{d}\rho_s^N(x)\mathrm{d}s.$$

for $\mathcal{F} = f(\langle \varphi, \cdot \rangle)$. Comparing the corresponding terms, we have therefore derived (3.5) for cylindrical functionals, while Proposition 3.4 offers strong enough approximation such that Itô's formula for the general case follows as we pass $M$ to infinity. □

We abbreviate from now on

$$\mathcal{V}_\mathcal{F}[\mu] := \int_{\mathbb{T}^d} a[\mu](x) : D_\mu^2 \mathcal{F}(\mu, x, x) \mathrm{d}\mu(x), \quad \mu \in \mathcal{P}(\mathbb{T}^d). \tag{3.11}$$

It follows from the definition (with the Lipschitz constants as defined in (A.16)) that

$$|\mathcal{V}_\mathcal{F}[\mu] - \mathcal{V}_\mathcal{F}[\nu]| \lesssim \|a\|_{\mathrm{lip}} \|D_\mu^2 \mathcal{F}\|_{\mathrm{lip}} \mathcal{W}_2(\mu, \nu), \tag{3.12}$$

since

$$\begin{aligned}|\mathcal{V}_\mathcal{F}[\mu] - \mathcal{V}_\mathcal{F}[\nu]| &\leq \left|\int_{\mathbb{T}^d} a[\nu](x) : D_\mu^2 \mathcal{F}(\nu, x, x) \mathrm{d}(\mu-\nu)(x)\right| \\ &+ \int_{\mathbb{T}^d} |a[\mu](x) : D_\mu^2 \mathcal{F}(\mu, x, x) - a[\nu](x) : D_\mu^2 \mathcal{F}(\nu, x, x)| \mathrm{d}\mu(x) \\ &\lesssim \sup_{f \in W^{1,\infty}} \left|\int_{\mathbb{T}^d} f \mathrm{d}(\mu - \nu)\right| \cdot \sup_{\mu \in \mathcal{P}} \|D_\mu^2 \mathcal{F}[\mu] : a[\mu]\|_{W^{1,\infty}} + \|D_\mu^2 \mathcal{F}\|_{\mathrm{lip}} \|a\|_{\mathrm{lip}} \mathcal{W}_2(\mu, \nu)\end{aligned}$$

where we can then use the Kantorovich–Rubinstein duality to derive (3.12). We now turn to the Itô's formula for the approximating SPDE (1.12) and emphasize the intricate point that the dynamic of the mollified equation is supported on $L^2(\mathbb{T}^d) \cap \mathcal{P}(\mathbb{T}^d)$ in comparison to $\mathcal{V}_\mathcal{F}[\cdot]$.

**Proposition 3.6. (Itô's formula for SPDE)** *We fix a function $\mathcal{F}[\cdot]$ on $\mathcal{P}(\mathbb{T}^d)$, which is assumed to have Lipschitz-continuous interior derivatives up to second order in the sense of (A.15), and let $\mathcal{V}_\mathcal{F}^N[\cdot]$ defined on the the subspace $L^2(\mathbb{T}^d) \cap \mathcal{P}(\mathbb{T}^d)$ of $\mathcal{P}(\mathbb{T}^d)$ be*

$$\mathcal{V}_\mathcal{F}^N[\mu] := \sum_{k \in \mathbb{Z}^d} |\theta_k^N|^2 \mathfrak{F}_{\mathbb{T}^d \times \mathbb{T}^d} (D_\mu^2 \mathcal{F}[\mu] : \Lambda_\mu^N)(k, -k), \tag{3.13}$$

*with the abbreviation that*

$$\Lambda_\mu^N(x, y) := f^N(\mu(x)) f^N(\mu(y)) \Sigma[\mu](x) \cdot \Sigma[\mu]^t(y).$$

*The Itô's formula holds for the solution to the regularized Dean-Kawasaki equation (1.12) in the following form:*

$$\begin{aligned}\mathcal{F}[m_t^N] - \mathcal{F}[m_0^N] &= \int_0^t \int_{\mathbb{T}^d} \{(a[m_s^N] : \nabla D_\mu) \mathcal{F}(m_s^N, x) + b[m_s^N](x) \cdot D_\mu \mathcal{F}(m_s^N, x)\} \mathrm{d}m_s^N(x) \mathrm{d}s \\ &+ \frac{1}{N} \int_0^t \mathcal{V}_\mathcal{F}^N[m_s^N] \mathrm{d}s + \tilde{M}_t^{N,\mathcal{F}},\end{aligned} \tag{3.14}$$

where $(\tilde{M}_t^{N,\mathcal{F}})_{t\geq 0}$ is a continuous martingale with quadratic variation s.t.

$$\mathbf{E}[\tilde{M}^{N,\mathcal{F}}]_t \lesssim \frac{\|\theta^N\|_{\ell^2}^2}{N}\|m_0^N\|_{L^1}^2.$$

**Proof.** We notice that the Itô correction term for the SPDE with cylindrical functionals for $\mathcal{F}[\cdot] = f(\langle\varphi,\cdot\rangle)$ as in the proof of Proposition 3.5 reads (besides a $1/N$ factor)

$$\mathcal{V}_{\mathcal{F}}^N[m_s^N] = \sum_{k\in\mathbb{Z}^d}|\theta_k^N|^2 \int_{\mathbb{T}^{2d}} \overline{e_k(x)}e_k(y) D_\mu^2\mathcal{F}[\mu](x,y) : [(f^N(\mu)\Sigma[\mu](x))(f^N(\mu)\Sigma^t[\mu](y))] \mathrm{d}x\mathrm{d}y,$$

which corresponds exactly to (3.13). We have used the quadratic variation

$$[B^k, B^j]_t = \delta_{k=-j} t,$$

the symmetry of $(\theta^N)$ for fixed $N\in\mathbb{N}$ and (under Einstein summation again)

$$D_\mu^2\mathcal{F}[\mu](x,y) : (\Sigma[\mu](x)\Sigma^t[\mu](y))$$
$$= \partial^{i,l}f(\langle\varphi,\mu\rangle)\partial^j\varphi_i(x)\partial^k\varphi_l(x)\left(\sum_{m=1}^d \Sigma_{m,j}^t[\mu](x)\Sigma_{m,k}^t[\mu](y)\right)$$
$$= \sum_{m=1}^d \partial^{i,l}f(\langle\varphi,\mu\rangle)(\Sigma_{m,j}^t[\mu](x)\partial_j\varphi_i(x))(\partial_k\varphi_l(x)\Sigma_{m,k}^t[\mu](y)).$$

The functional in (3.13) is well-defined for $\mu\in L^2(\mathbb{T}^d)\cap\mathcal{P}(\mathbb{T}^d)$ and $\mathcal{F}[\cdot]$ with Lipschitz interior derivatives up to second order, since

$$\mathcal{V}_{\mathcal{F}}^N[\mu] \leq \|\theta^N\|_{\ell^2}^2 \sup_{k\in\mathbb{Z}^d}|\mathfrak{F}(D_\mu^2\mathcal{F}[\mu]:\Lambda_\mu^N)(k,-k)| \leq a\|\theta^N\|_{\ell^2}^2 \|D_\mu^2\mathcal{F}[\mu]:\Lambda_\mu^N\|_{L^1(\mathbb{T}^{2d})} < \infty, \quad (3.15)$$

where $D_\mu^2\mathcal{F}[\mu]\in L^\infty(\mathbb{T}^{2d})$, $\Sigma[\mu]\in L^\infty(\mathbb{T}^d)$ and $f^N(\mu)\in L^2(\mathbb{T}^d)$, and all bounds are uniform in $\mu\in\mathcal{P}\cap L^2$. Itô's formula for the general case can also be proved by the approximation introduced in Proposition 3.4. □

From what we have derived above, it is straightforward to see that the formal generators of two systems are almost identical, except for the Itô correction terms which are scaled by $N^{-1}$. The following Lemma further considers the difference of the correction terms, where it can be noticed that the Fisher information, i.e. the $L^2$ norm of the gradient of the square root, appears and will be closely related to an entropy bound for the solution to (1.12), which we will derive below.

Since the comparison also involves the approximation of functions on the diagonal due to (3.11), some regularity assumption on $D_\mu^2\mathcal{F}[\cdot]$ is expected. It turns out that uniform spatial $W^{2,\infty}$ regularity suffices for our purpose. We recall that $(\delta_N)_{N\geq 1}$ are the approximation parameters in (1.14) of the square root by $(f^N)_{N\geq 1}$, and $(\theta^N)_{N\geq 1}$ are the amplitudes of noise satisfying Assumption 1.3. In particular, $(\theta_k^N)_{N\geq 1}$ depends on $(L^N)_{N\geq 1}$ through (1.16), and

$$\int_{\mathbb{T}^d}|y|^2|\psi^N(y)|^2\mathrm{d}y \lesssim L_N^{-2}$$

for $\psi^N := \sum_{k \in \mathbb{Z}^d} |\theta_k^N|^2 e_k$ according to (1.18).

**Lemma 3.7.** *For $\mathcal{V}_\mathcal{F}^N$ and $\mathcal{V}_\mathcal{F}$ given above, we have for any $\mu \in \mathcal{P}(\mathbb{T}^d) \cap L^1(\mathbb{T}^d)$ that*

$$|\mathcal{V}_\mathcal{F}^N[\mu] - \mathcal{V}_\mathcal{F}[\mu]| \lesssim \|D_\mu^2 \mathcal{F}[\mu]\|_{W^{2,\infty}} \left[\delta_N + L_N^{-2}\left(1 + \int_{\mathbb{T}^d} \frac{|\nabla \mu|^2}{\mu}\right)\right], \qquad (3.16)$$

*where the constant is independent of $N \in \mathbb{N}$ and $\mu$.*

**Proof.** We assume without loss of generality that the R.H.S. of (3.16) is finite and split the estimate in (3.16) into two parts by

$$|\mathcal{V}_\mathcal{F}^N[\mu] - \mathcal{V}_\mathcal{F}[\mu]| \leq \left|\sum_{k \in \mathbb{Z}^d} |\theta_k^N|^2 \mathfrak{F}(\Lambda_\mu^N)(k, -k) - \int_{\mathbb{T}^d} \Lambda_\mu^N(x, x) \mathrm{d}x\right|$$
$$+ \left|\int_{\mathbb{T}^d} D_\mu^2 \mathcal{F}[\mu](x, x) : a[\mu](x) \left[f^N(\mu(x))^2 - \mu(x)\right] \mathrm{d}x\right| := (\mathrm{I}) + (\mathrm{II}), \qquad (3.17)$$

where we denote

$$\Lambda_\mu^N(x, y) := D_\mu^2 \mathcal{F}[\mu](x, y) : \left[(\Sigma[\mu] f^N(\mu))(x) \cdot (\Sigma^t[\mu] f^N(\mu))(y)\right].$$

According to the technical Lemma A.4 with $f = \Sigma[\mu] f^N(\mu)$ and $h = D_\mu^2 \mathcal{F}[\mu]$, the first term in (3.16) can be bounded by

$$(\mathrm{I}) \lesssim L_N^{-2} \|\Sigma[\mu] f^N(\mu)\|_{W^{1,2}}^2 \|D_\mu^2 \mathcal{F}[\mu]\|_{W^{2,\infty}}, \qquad (3.18)$$

where the constant is independent of $N$ and $\mu$. The second term in (3.16) is easily bounded by

$$(\mathrm{II}) \overset{(1.14)}{\lesssim} \delta_N \|D_\mu^2 \mathcal{F}[\mu] : a[\mu]\|_{L^\infty} \overset{(2.5)}{\lesssim} \delta_N \|D_\mu^2 \mathcal{F}[\mu]\|_{L^\infty}, \qquad (3.19)$$

It follows immediately from the estimates above and

$$\|\Sigma[\mu] f^N(\mu)\|_{W^{1,2}} \lesssim \|\Sigma[\mu]\|_{W^{1,\infty}} \|f^N(\mu)\|_{W^{1,2}} \overset{(1.14, 2.8)}{\lesssim} 1 + \int_{\mathbb{T}^d} \frac{|\nabla \mu|^2}{\mu},$$

where the constants are independent of $N$, that (3.16) holds. $\square$

## 3.3 Kolmogorov Equations on the Space of Probability Measures

From now on, we will fix a terminal time $T > 0$. Due to the approximation of the formal generators of the SPDEs to that of the particle systems as the particle number $N$ goes to infinity, we intend to solve the corresponding *Kolmogorov equation* for the derivation of the weak error rate. As it turns out, for Dean-Kawasaki equation, which contains the $N^{-1/2}$ scaling on the noise, it suffices to solve the corresponding Kolmogorov equations up to the first order *twice*, rather than solving the whole equations.

The results of solving such a first-order equation by characteristic line are well established in the literature of mean field equations, which we present briefly for completeness. More explicitly, we consider the following infinite-dimensional backward equation on $\mathcal{P}(\mathbb{T}^d)$:

$$\partial_t \mathcal{G}(t, \mu) + \int_{\mathbb{T}^d} \{a[\mu](x) : \nabla D_\mu \mathcal{G}(t, \mu, x) + b[\mu](x) \cdot D_\mu \mathcal{G}(t, \mu, x)\} \mathrm{d}\mu(x) = 0 \qquad (3.20)$$

with terminal condition $\mathcal{G}(T,\mu) = \mathcal{F}[\mu]$ for $T > 0$ that we will fix, where $a[\cdot]$ and $b[\cdot]$ are given by (1.2) and (1.8) respectively.

**Definition 3.8.** *For the fixed terminal time $T < \infty$, we define a classical solution to (3.20) to be some $\mathcal{G}:[0,T] \times \mathcal{P}(\mathbb{T}^d) \to \mathbb{R}$, which is continuously differentiable with respect to time, has Lipschitz-continuous first-order interior derivative and solves (3.20) pointwise.*

It is then established in [BLPR17] that (3.20) is in fact a transport equation on the space of probability measures, whose characteristic lines are the solutions to the corresponding McKean-Vlasov equation.

**Proposition 3.9.** ([BLPR17], Theorem 7.2) *Under the assumption that the coefficients $b[\cdot]$, $a[\cdot]$ and the test functional $\mathcal{F}[\cdot]$ all have jointly Lipschitz-continuous interior derivatives up to the second order, (3.20) admits a unique solution in the sense of Definition 3.8 with the terminal condition $\mathcal{G}(T,\mu) = \mathcal{F}[\mu]$, which is given by the characteristic line*

$$\mathcal{G}(s,\mu) = \mathcal{F}[\mathcal{L}(X_T^{s,\mu})], \quad s \in [0,T], \tag{3.21}$$

*where $(X_r^{s,\mu})_{r \in [s,T]}$ is the solution to the forward McKean–Vlasov SDE*

$$dX_r^{s,\mu} = b(X_r^{s,\mu}, \mathcal{L}(X_r^{s,\mu}))dr + \sqrt{2}\Sigma(X_r^{s,\mu}, \mathcal{L}(X_r^{s,\mu})) \cdot dB_r, \tag{3.22}$$

*with initial condition $\mathcal{L}(X_s^{s,\mu}) = \mu$. Especially, $\mathcal{G}(s,\cdot)$ has Lipschitz-continuous interior derivatives up to the second order for any time $s \in [0,T]$.*

Due to the uniform-Lipschitz condition which we imposed on the drift and the noise, the classical theory of McKean–Vlasov equations shows that there exists a unique probabilistically strong solution to (3.22). Besides, it follows directly from Itô's formula on $|X_\cdot^{s,\mu} - X_\cdot^{s,\nu}|^2$, the Grönwall inequality and the definition of the Wasserstein distance $\mathcal{W}_2$ that

$$|\mathcal{W}_2(\mathcal{L}(X_r^{s,\mu}), \mathcal{L}(X_r^{s,\mu}))| \lesssim \mathcal{W}_2(\mu, \nu), \tag{3.23}$$

which further implies that

$$|\mathcal{G}(s,\mu) - \mathcal{G}(s,\nu)| \lesssim \|\mathcal{F}\|_{\text{lip}} \mathcal{W}_2(\mu, \nu), \tag{3.24}$$

for any $s \in [0,T]$. We refer the reader to [CM18, HRW21] and references therein for a comprehensive review of the properties of McKean–Vlasov equations and the proof of (3.23).

# 4 Weak-Error Estimation

In this chapter, we derive a weak-error bound for the approximation of the Dean-Kawasaki equation. Although a solution theory to the Kolmogorov equation related to the SPDE (1.12) seems beyond the reach of current techniques, we only need to solve the equation up to the first order, i.e. considering a transport equation as explained above. On the other hand, we need both the test function and the coefficients of the equations to be sufficiently regular (in the sense of Assumption 4.1 below), so that the solution inherits this regularity.

## 4.1 Regularity of Solutions and Linearized Equations

In this section, we derive regularity estimates for the transport equation mentioned above. The well-posedness is provided by Proposition 3.9, which is due to [BLPR17] and which makes use of the deep link observed by Lions [Lio] between the interior derivatives of functionals on probability spaces with Fréchet derivatives of the lift of such functionals to the space of $L^2$ random variables.

These probabilistic arguments yield spatial $W^{1,\infty}$ regularity of the interior derivatives of the solutions to the Kolmogorov equation, while it seems less straightforward to derive higher spatial regularity of the solutions (see e.g. Chapter 6 in [BLPR17]). However, such spatial regularity is indispensable to our weak error estimates, see Lemma 3.7. Therefore, we will use analytic rather than probabilistic arguments, as in [Tse21, DT21].

The goal of this section is thus to prove that[4.1]

$$\sup_{\mu \in L^1 \cap \mathcal{P}} \sup_{t \in [0,T]} \|D_\mu^2 \{\mu \mapsto \mathcal{F}[\mathcal{L}(X_t^\mu)]\}\|_{W^{2,\infty}(\mathbb{T}^{2d})} < \infty, \tag{4.1}$$

if the functional $\mathcal{F}[\cdot]$ and the coefficients $b[\cdot]$ and $\Sigma[\cdot]$ are assumed to be sufficiently regular, in the sense that Assumption 2.1 holds. In particular, $\|a[\mu]\|_{W^{1,\infty}}$ and $\|b[\mu]\|_{L^\infty}$ are bounded uniformly in $\mu \in L^1 \cap \mathcal{P}$ according to (2.6), and $a_{\min} \leq a[\mu] \lesssim 1$ uniformly for $\mu \in L^1 \cap \mathcal{P}$ due to (2.5). The reason for considering (4.1) will be clear in (4.26) below.

We further impose the following regularity conditions in addition to Assumption 2.1:

**Assumption 4.1.** *The test functional $\mathcal{F}[\cdot]$ has Lipschitz-continuous interior derivatives up to the fourth order s.t.*

$$\max_{k \leq 4} \sup_{\mu \in L^1 \cap \mathcal{P}} \left\| \frac{\delta^k \mathcal{F}}{\delta \mu^k} \right\|_{W_{\mathbb{T}^{kd}}^{4,\infty}} < \infty, \tag{4.2}$$

*and the coefficients $b[\cdot]$ and $a[\cdot]$ have Lipschitz-continuous interior derivatives up to the fourth order, and furthermore, we assume they have spatial regularity*

$$\max_{k \leq 2} \sup_{\mu \in L^1 \cap \mathcal{P}} \left( \left\| \frac{\delta^k a}{\delta \mu^k} \right\|_{W_{\mathbb{T}^{(k+1)d}}^{6,\infty}} \vee \left\| \frac{\delta^k b}{\delta \mu^k} \right\|_{W_{\mathbb{T}^{(k+1)d}}^{5,\infty}} \right) < \infty. \tag{4.3}$$

The reason for requiring four derivatives with respect to the measure argument is that we need to solve the transport equation (3.23) effectively twice, with a nonlinear function of the first solution as terminal condition for the second equation. For Lemma 4.5 and Corollary 4.6, it suffices to have interior derivatives up to the second order. It would be possible to relax the assumption that all interior derivatives are in $W^{5,\infty}$, which we make for brevity.

---

4.1. We will denote

$$D_\mu^2 \{\mathcal{F}[\mathcal{L}(X_t^\mu)]\} = D_\mu^2 \{\mu \mapsto \mathcal{F}[\mathcal{L}(X_t^\mu)]\} : \mathcal{P}(\mathbb{T}^d) \times \mathbb{T}^{2d} \to \mathbb{R}$$

for $t \in [0, T]$ from now on.

**Example 4.2.** In the context of Example 2.3, the regularity assumptions above amount to the condition that $\sigma \in C_b^7(\mathbb{R})$ and $K, G \in W^{5,\infty}(\mathbb{T}^d)$, since we have

$$\frac{\delta^k}{\delta \mu^k}[\sigma(G * \mu(x))](y_1, \ldots, y_k) = \nabla^k \sigma(G * \mu(x))[G(x - y_1), \ldots, G(x - y_k)]$$

modulo some normalization constant.

We denote from now on the distribution of $X^\mu := X^{0,\mu}$ defined in (3.22) as

$$m_t^\mu := \mathcal{L}(X_t^\mu), \quad \text{for} \quad 0 \le t \le T,$$

which solves the following nonlinear Fokker-Planck equation

$$\begin{cases} \partial_t m_t^\mu = \nabla^2 : (a[m_t^\mu] m_t^\mu) - \nabla \cdot (b[m_t^\mu] m_t^\mu), \\ m_0^\mu = \mu, \end{cases} \tag{4.4}$$

and we will use throughout the almost trivial regularity estimate for the McKean–Vlasov equation (3.22) that for any initial datum $\mu \in L^1(\mathbb{T}^d) \cap \mathcal{P}(\mathbb{T}^d)$,

$$\sup_{t \in [0,T]} \|m_t^\mu\|_{L^1} = 1. \tag{4.5}$$

We also write $P_{\cdot,\cdot}^\mu$ as the solution map of the inhomogeneous decoupled forward equation

$$\begin{cases} \partial_t P_{s,t}^\mu \varphi = \nabla^2 : (a[m_t^\mu] P_{s,t}^\mu \varphi) - \nabla \cdot (b[m_t^\mu] P_{s,t}^\mu \varphi), \\ P_{s,s}^\mu \varphi = \varphi, \end{cases} \tag{4.6}$$

for any bounded measure $\varphi$ on $\mathbb{T}^d$ and $0 \le s \le t \le T$. Both the maps $[0, T] \times \mathcal{P}(\mathbb{T}^d) \ni (t, \mu) \mapsto m_t^\mu$ and $[0, T] \times \mathcal{P}(\mathbb{T}^d) \ni (t, \nu) \mapsto P_t^\mu \nu$ are flows of probability measures, and it is worth emphasizing that (4.6) decouples the dependence of $m^\mu$ of $\mu$ on the flow from that on the initial condition, which can also be seen later from (A.27) in the formal derivation.

Following the literature, we will denote from now on

$$W^{-k,1} := (W^{k,\infty})^*, \tag{4.7}$$

where $k > 0$. The following classical regularity estimates for (4.6) are needed.

**Lemma 4.3.** *We assume that the assumptions 2.1 and 4.1 hold and fix a $\mu \in \mathcal{P}(\mathbb{T}^d)$ with $L^1$ density. Considering $(P_{s,t}^\mu)_{0 \le s \le t}$ as the evolution system defined in (4.6), the following uniform bound follows:*

$$\sup_{\mu \in L^1 \cap \mathcal{P}} \left\| \int_0^\cdot P_{s,\cdot}^\mu (\nabla \cdot f_s) ds \right\|_{L^\infty([0,T], W^{-k,1})} \lesssim_{T,k} \|f\|_{L^\infty([0,T], W^{-k,1})}, \tag{4.8}$$

*for any $f \in L^\infty([0, T], W^{-k,1}(\mathbb{T}^d, \mathbb{R}^d))$ and $1 \le k \le 4$.*

**Proof.** We consider the Kolmogorov backward equation with respect to the dual of $P_{s,t}^\mu$, whose solution operator we denote by $(P_{s,t}^\mu)^*$. Duality here means that $\langle \varphi, (P_{s,t}^\mu)^* f \rangle = \langle P_{s,t}^\mu \varphi, f \rangle$. According to Theorem 3.2.6 and Corollary 2.4.5 in [Str08], we can estimate the regularity of the backward equation by

$$\|(P_{s,t}^\mu)^* (\nabla \cdot f_t)\|_{W^{k,\infty}} \lesssim \min \left( \|f_t\|_{W^{k+1,\infty}}, \frac{1}{1 \wedge |t-s|^{k/2}} \|f_t\|_{W^{1,\infty}}, \right)$$

for $0 \leq k' \leq 4$, which implies by Gagliardo–Nirenberg interpolation (c.f. [BM18]) that

$$\|(P_{s,t}^\mu)^*(\nabla \cdot f_t)\|_{W^{k,\infty}} \lesssim \frac{1}{1 \wedge |t-s|^{1/2}} \|f_t\|_{W^{k,\infty}}, \tag{4.9}$$

for $1 \leq k \leq 4$. By duality and recalling the definition (4.7), we thus have for the same range of $k$ that

$$\|P_{s,t}^\mu(\nabla \cdot f_s)\|_{W^{-k,1}} \lesssim \frac{1}{1 \wedge |t-s|^{1/2}} \|f_s\|_{W^{-k,1}}. \tag{4.10}$$

The multiplicative constants above only depend on $k$ and the regularity of the coefficients up to the fifth order, which we can further control by

$$\sup_{t \geq 0} \|a[m_t^\mu]\|_{W^{5,\infty}} \leq \sup_{\mu \in L^1 \cap \mathcal{P}} \|a[\mu]\|_{W^{5,\infty}} < \infty$$

and its counterpart for $b[\cdot]$ using (4.5) and the assumption (4.3). In particular, (4.10) is uniform in $\mu \in L^1 \cap \mathcal{P}(\mathbb{T}^d)$. We have therefore

$$\left\|\int_0^\cdot P_{s,\cdot}^\mu(\nabla \cdot f_s) ds\right\|_{L^\infty([0,T],W^{-k,1})} \lesssim \sup_{t \leq T} \int_0^t \frac{ds}{1 \wedge |t-s|^{1/2}} \|f\|_{L^\infty([0,T],W^{-k,1})} \lesssim_T \|f\|_{L^\infty([0,T],W^{-k,1})},$$

which is then also uniform in $\mu \in L^1 \cap \mathcal{P}(\mathbb{T}^d)$. □

In the previous work [Tse21] on the weak propagation of chaos, Lipschitz-continuity of the exterior derivatives (up to arbitrary degrees) of the nonlinear functionals with respect to the law of the McKean–Vlasov equations is derived in [Tse21]. The Lipschitz-continuity of the interior derivatives is then considered in [DT21], where the first-order linearization of the spatial variable is implemented. Besides those regularity bounds, we need furthermore a higher spatial-regularity control of the interior derivatives. We propose therefore the following *duality argument* using (4.5), which extends the Theorem (Main result) in [Tse21].

**Remark 4.4.** Although only the special case of constant diffusion coefficients is rigorously considered in [Tse21], similar results for the case considered in this paper[4.2] also hold true, since under the smoothness assumptions on the coefficients, the heat kernel bound essentially used in the proof in [Tse21] (for example (2.6) therein) also holds for the fundamental solution $[0,T] \times \mathbb{T}^{2d} \ni (t,x,y) \mapsto p^\mu(t,x,y) = (P_{0,t}^\mu \delta_x)(y)$. Specifically, we have for any $t > 0$, it holds that

$$\sup_{x,y \in \mathbb{T}^d} |\nabla p^\mu(t,x,y)| \lesssim t^{-\frac{1+d}{2}} g\left(\frac{|x-y|}{t^{1/2}}\right), \tag{4.11}$$

for some bounded $g$, and the constant is uniform in $\mu \in L^1 \cap \mathcal{P}(\mathbb{T}^d)$ due to the same reason as in (4.8). The proof can be found for example in Theorem 3.3.11 in [Str08]. The proofs of the main results in [Tse21] can therefore be carried out mutatis muntandis, if we take into acount the linearization of $a[m^\mu]$.

---

[4.2]. as considered in [BLPR17]

In the following technical lemma, we discuss the linearization of functionals of the McKean-Vlasov equation. The representation in (4.12) can be formally guessed by applying the (second order) chain rule to $\mathcal{F}[m^\mu]$ and then differentiating the evolution equation for $m^\mu$ with respect to $\mu$. We focus on the second interior derivative $D_\mu^2$, for which we need higher space regularity than in the references in order to apply Lemma 3.7. Later we will also need Lipschitz continuity of the derivatives $D_\mu^k$ for larger $k$, which is derived in [Tse21] and [DT21].

We recall the definition of $\tilde{b}[\cdot]$ in (1.13), which satisfies according to (4.3) that

$$\max_{k \leq 2} \sup_{\mu \in L^1 \cap \mathcal{P}} \left\| \frac{\delta^k \tilde{b}}{\delta \mu^k} \right\|_{W^{5,\infty}_{\mathbb{T}^{(k+1)d}}} < \infty.$$

**Lemma 4.5.** *We assume that the assumptions 2.1 and 4.1 are in force. For any $\phi, \psi \in C^\infty(\mathbb{T}^d)$ with zero mean and $\mu \in L^1(\mathbb{T}^d) \cap \mathcal{P}(\mathbb{T}^d)$, the following equality holds for a.e. $t \in [0, T]$:*

$$\left\langle \frac{\delta^2}{\delta \mu^2} \{\mathcal{F}[m_t^\mu]\}, \psi \otimes \phi \right\rangle = \left\langle \frac{\delta \mathcal{F}}{\delta \mu}[m_t^\mu], I_t^\mu(\psi, \phi) \right\rangle + \left\langle \frac{\delta^2 \mathcal{F}}{\delta \mu^2}[m_t^\mu], \mathcal{J}_t^\mu \psi \otimes \mathcal{J}_t^\mu \phi \right\rangle. \tag{4.12}$$

Here we denote $\mathcal{J}_t^\mu \psi := P_{0,t}^\mu \psi + \mathcal{J}_t^{\mu,1} \psi$ with the solution map $\mathcal{J}^{\mu,1}$ to the equation

$$\mathcal{J}_t^{\mu,1} \psi = \int_0^t P_{s,t}^\mu \nabla \cdot \left\{ \left[ \int_{\mathbb{T}^d} \frac{\delta a(\cdot, \mu)}{\delta \mu}[m_s^\mu](y) \mathrm{d}(\mathcal{J}_s^{\mu,1} \psi)(y) \right] \nabla m_s^\mu \right.$$
$$\left. - \left[ \int_{\mathbb{T}^d} \frac{\delta \tilde{b}(\cdot, \mu)}{\delta \mu}[m_s^\mu](y) \mathrm{d}(\mathcal{J}_s^{\mu,1} \psi)(y) \right] m_s^\mu + \mathcal{R}_s^\mu(\psi) \right\} \mathrm{d}s, \tag{4.13}$$

with the evolution system $(P_{\cdot,\cdot}^\mu)$ that is defined in (4.6) and a map $\mathcal{R}^\mu: [0, T] \times W^{-3,1}(\mathbb{T}^d) \to W^{-1,1}(\mathbb{T}^d)$ that we set to be[4.3]

$$\mathcal{R}_t^\mu(\psi) = \left( \int_{\mathbb{T}^d} \frac{\delta a(\mu, \cdot)}{\delta \mu}[m_t^\mu](y) \mathrm{d} P_{0,t}^\mu(\psi)(y) \right) \nabla m_t^\mu$$
$$- m_t^\mu \left( \int_{\mathbb{T}^d} \frac{\delta \tilde{b}(\mu, \cdot)}{\delta \mu}[m_t^\mu](y) \mathrm{d} P_{0,t}^\mu(\psi)(y) \right), \tag{4.14}$$

for any $t \in [0, T]$ and $\psi \in W^{-3,1}(\mathbb{T}^d)$, which satisfies

$$\sup_{\mu \in L^1 \cap \mathcal{P}} \sup_{0 \leq t \leq T} \|\mathcal{R}_t^\mu(\psi)\|_{W^{-1,1}} \lesssim \|\psi\|_{W^{-3,1}}.$$

*The bilinear map $I_t^\mu$ are defined as the sum of the following two bilinear maps $I_t^\mu = I_t^{\mu,1} + I_t^{\mu,2}$, where*

$$I_t^{\mu,1}(\psi, \phi) = \int_0^t P_{s,t}^\mu \nabla \cdot \left\{ \left[ \int_{\mathbb{T}^d} \frac{\delta a(\cdot, \mu)}{\delta \mu}[m_s^\mu](y) \mathrm{d}(\mathcal{J}_s^\mu \phi)(y) \right] \nabla P_{0,s}^\mu(\psi) \right.$$
$$\left. - \left[ \int_{\mathbb{T}^d} \frac{\delta \tilde{b}(\cdot, \mu)}{\delta \mu}[m_s^\mu](y) \mathrm{d}(\mathcal{J}_s^\mu \phi)(y) \right] P_{0,s}^\mu(\psi) \right\} \mathrm{d}s, \tag{4.15}$$

---

4.3. Note the divergence outside the curly bracket in (4.12).

and $I^{\mu,2}$ is the solution map to the parabolic non-local equation given by

$$I_t^{\mu,2}(\psi,\phi) = \int_0^t P_{s,t}^\mu \nabla \cdot \left\{ \left[ \int_{\mathbb{T}^d} \frac{\delta a(\mu,\cdot)}{\delta \mu}[m_s^\mu](y) \mathrm{d} I_s^{\mu,2}(\psi,\phi)(y) \right] \nabla m_s^\mu \right.$$
$$\left. - \left[ \int_{\mathbb{T}^d} \frac{\delta \tilde{b}(\mu,\cdot)}{\delta \mu}[m_s^\mu](y) \mathrm{d} I_s^{\mu,2}(\psi,\phi)(y) \right] m_s^\mu + \mathscr{T}_s^\mu(\psi,\phi) \right\} \mathrm{d}s, \qquad (4.16)$$

for a map $\mathscr{T}^\mu : [0,T] \times (W^{-3,1})^2(\mathbb{T}^d) \to W^{-2,1}(\mathbb{T}^d)$, such that

$$\sup_{\mu \in L^1 \cap \mathcal{P}} \sup_{0 \le t \le T} \|\mathscr{T}^\mu(\psi)\|_{W^{-2,1}} \lesssim \|\psi\|_{W^{-3,1}},$$

and it is also bilinear in the $(W^{-3,1})^2(\mathbb{T}^d)$-arguments[4.4]. Lastly, we have for any $\mu \in L^1(\mathbb{T}^d) \cap \mathcal{P}(\mathbb{T}^d)$ and $T < \infty$ fixed that the following estimates hold true:

$$\|\mathcal{J}^{\mu,1}\psi\|_{L^\infty([0,T],W^{-1,1})} \lesssim \|\psi\|_{W^{-3,1}}, \qquad (4.17)$$

and

$$\|I^{\mu,j}(\psi,\phi)\|_{L^\infty([0,T],W^{-4,1})} \lesssim \|\phi\|_{W^{-3,1}} \|\psi\|_{W^{-3,1}}, \quad j=1,2, \qquad (4.18)$$

where the constants are uniform in $\mu \in L^1(\mathbb{T}^d) \cap \mathcal{P}(\mathbb{T}^d)$.

The terms in the square bracket with $a[\cdot]$ used in (4.12) should be interpreted as

$$x \mapsto \int_{\mathbb{T}^d} \frac{\delta a(x,\mu)}{\delta \mu}[m_s^\mu](y) \mathrm{d}(\mathcal{J}_s^{\mu,1}\psi)(y),$$

and similarly for the others.

**Proof.** An outline essentially due to [Tse21] for the derivation of the formula (4.12), when $\phi, \psi$ are instead assumed to be of the forms

$$\phi \in \mathcal{P}(\mathbb{T}^d) - \mu, \quad \psi \in \mathcal{P}(\mathbb{T}^d) - \mu,$$

together with the derivation of the terms involved ($I^\mu$, $\mathcal{J}^\mu$, $\mathscr{R}^\mu$ and $\mathscr{T}^\mu$) are included in the Appendix A.4, where the differences with $\mu$ naturally appear due to the very definition of exterior derivatives (A.13). The proof of the boundedness of $\mathscr{R}^\mu$ is also postponed to Lemma A.8[4.5] directly after the derivation, and the existence and uniqueness of the linear nonlocal equations (4.12) and (4.15) given the regularity of the forcing terms are direct application of Leray-Schauder, whose proof can be found in Theorem 2.3 in [Tse21], in the case of constant diffusion coefficients. We recall the discussion on the general case of diffusions in Remark 4.4.

We now focus on the bounds (4.17). Using the equation (4.12) that $\mathcal{J}^{\mu,1}$ satisfies, the parabolic regularity estimate (4.8) of $(P_{\cdot,\cdot}^\mu)$ and the trivial embedding $W^{-1,1} \subset W^{-3,1}$, we have

$$\|\mathcal{J}_t^{\mu,1}\psi\|_{W^{-1,1}} \lesssim \int_0^t (t-s)^{-\frac{1}{2}} \|\mathcal{J}_s^{\mu,1}\psi\|_{W^{-1,1}} \mathrm{d}s + \sup_{t \le T} \|\mathscr{R}_t^\mu(\psi)\|_{W^{-1,1}},$$

---

[4.4]. The exact expression of $\mathscr{T}^\mu$ is postponed to (A.8) for brevity.

[4.5]. The proof of the boundedness of $\mathscr{R}^\mu$ and $\mathscr{T}^\mu$ is straightforward but lengthy, so we also postpone it; the bound for $\mathscr{T}^\mu$ depends on the estimate (4.17), which we derive below without relying on a bound for $\mathscr{T}^\mu$ and therefore the argument is not circular.

which implies by Young and Grönwall's inequalities (c.f. for example the proof of Theorem 5.1 in [EK86]) that

$$\|\mathcal{J}_t^{\mu,1}\psi\|_{W^{-1,1}} \lesssim \sup_{t\leq T}\|\mathcal{R}_t^\mu(\psi)\|_{W^{-1,1}}.$$

Uniform boundedness (in $\mu$) of $\mathcal{R}^\mu$ then leads to (4.17). Given the control of $\mathcal{J}^{\mu,1}$ via (4.17), it is then straightforward to show the boundedness of $\mathcal{T}^\mu$ using Lemma A.8. For the inequality (4.18) with $I^{\mu,1}$, we bound the first term in (4.14) by

$$\sup_{t\in[0,T]}\left\|\int_0^t P_{s,t}^\mu \nabla\cdot\left[\int_{\mathbb{T}^d}\frac{\delta a(\cdot,\mu)}{\delta\mu}[m_s^\mu](y)\mathrm{d}(\mathcal{J}_s^\mu\phi)(y)\nabla P_{0,s}^\mu(\psi)\right]\mathrm{d}s\right\|_{W^{-4,1}}$$

$$\overset{(4.8)}{\lesssim} \sup_{t\in[0,T]}\left\|\left[\int_{\mathbb{T}^d}\frac{\delta a(\cdot,\mu)}{\delta\mu}[m_t^\mu](y)\mathrm{d}(\mathcal{J}_t^\mu\phi)(y)\right]\nabla P_{0,t}^\mu(\psi)\right\|_{W^{-4,1}}$$

$$\lesssim \sup_{\mu\in L^1\cap\mathcal{P}}\left\|\frac{\delta a}{\delta\mu}[\mu]\right\|_{W^{4,\infty}}\sup_{t\in[0,T]}\|\mathcal{J}_t^\mu(\phi)\|_{W^{-3,1}}\sup_{t\in[0,T]}\|\nabla P_{0,t}^\mu(\psi)\|_{W^{-4,1}}$$

$$\lesssim \sup_{\mu\in L^1\cap\mathcal{P}}\left\|\frac{\delta a}{\delta\mu}[\mu]\right\|_{W^{4,\infty}}\|\phi\|_{W^{-3,1}}\|\psi\|_{W^{-3,1}},$$

where we use the fact that multiplication from $W^{-4,1}\times W^{4,\infty}$ to $W^{-4,1}$ is continuous since $W^{4,\infty}$ is an algebra for the second inequality and the estimate

$$\sup_{t\in[0,T]}\|P_{0,t}^\mu\psi\|_{W^{-3,1}} \lesssim \|\psi\|_{W^{-3,1}}$$

and the bound on $\mathcal{J}^\mu$ for the last one. A similar bound holds for $\tilde{b}[\cdot]$. The bound in (4.18)[4.6] with $I^{\mu,2}$ given by (4.15) is implied by (A.37) together with Young and Grönwall's inequalities as above.

Lastly, to show that (4.12) holds for smooth functions instead of differences of probability measures, we fix and decompose $\phi,\psi\in C^\infty(\mathbb{T}^d)$ with means zero as

$$\phi = \|\phi_+\|_{L^1}\left(\frac{\phi_+}{\|\phi_+\|_{L^1}}-\mu\right) - \|\phi_-\|_{L^1}\left(\frac{\phi_-}{\|\phi_-\|_{L^1}}-\mu\right)$$

(similarly for $\psi$) with the convention that $0/0=0$ and using the assumption that $\|\phi_+\|_{L^1}=\|\phi_-\|_{L^1}$. We derive therefore, by abbreviating $\phi_\pm^\mu := \frac{\phi_\pm}{\|\phi_\pm\|_{L^1}}-\mu\in\mathcal{P}(\mathbb{T}^d)-\mu$ and noticing that (4.12) holds for $(\phi_\pm,\psi_\pm)$ and both sides of (4.12) are bilinear, that the equality follows also for $\phi,\psi\in C^\infty(\mathbb{T}^d)$. □

Finally, we have the desired regularity estimate (4.1), which now follows from the technical Lemma 4.5 above.

**Corollary 4.6.** *Under the assumptions 2.1 and 4.1, we have*

$$\sup_{\mu\in L^1\cap\mathcal{P}}\sup_{t\in[0,T]}\|D_\mu^2\{\mathcal{F}[m_t^\mu]\}\|_{W^{2,\infty}} \lesssim \max_{k\leq 2}\sup_{\mu\in L^1\cap\mathcal{P}}\left\|\frac{\delta^k\mathcal{F}}{\delta\mu^k}\right\|_{W^{4,\infty}}. \tag{4.19}$$

---

4.6. $I^{\mu,2}$ is in fact a bit more regular, but we only care about the sum $I^\mu$.

**Proof.** It follows from the regularity estimates for the functionals defined in (4.12), (4.14) and (4.15) in Lemma 4.5 above that

$$\left|\left\langle \frac{\delta^2}{\delta\mu^2}\{\mathcal{F}[m_t^\mu]\}, \psi \otimes \phi \right\rangle\right| \lesssim \max_{k \leq 2} \sup_{\mu \in \mathcal{P}} \left\|\frac{\delta^k \mathcal{F}}{\delta\mu^k}\right\|_{W^{4,\infty}} \|\phi\|_{W^{-3,1}} \|\psi\|_{W^{-3,1}},$$

uniformly in $t \in [0, T]$. We have therefore

$$\left|\left\langle \nabla_1^j \nabla_2^j \frac{\delta^2}{\delta\mu^2}\{\mathcal{F}[m_t^\mu]\}, \psi \otimes \phi \right\rangle\right| = \left|\left\langle \frac{\delta^2}{\delta\mu^2}\{\mathcal{F}[m^\mu]\}, \nabla^j\psi \otimes \nabla^j\phi \right\rangle\right|$$

$$\lesssim \max_{k \leq 2} \sup_{\mu \in \mathcal{P}} \left\|\frac{\delta^k \mathcal{F}}{\delta\mu^k}\right\|_{L^\infty} \|\psi\|_{L^1} \|\phi\|_{L^1}, \quad (4.20)$$

for $1 \leq j \leq 3$, where the differentiation $\nabla_1$ falls on the first variable in $\mathbb{T}^d \times \mathbb{T}^d$ corresponding to $\psi$ and $\nabla_2$ on the second. We used in the inequality above that $\nabla^j : L^1 \to W^{-3,1}$ is continuous for $j \leq 3$, since $L^1 \subset (L^\infty)^*$. We now claim that (4.19) implies that the distributional derivatives of $\delta^2\{\mathcal{F}[m_t^\mu]\}/\delta\mu^2$ are in $L^\infty(\mathbb{T}^{2d})$.

Indeed, (4.19) holds for generic $\phi, \psi \in C^\infty(\mathbb{T}^d)$ without mean-zero condition since $j \geq 1$, while since the closed unit ball $B_{L^1}(1) := \{f \in L^1(\mathbb{T}^d) : \|f\|_{L^1} \leq 1\}$ is separable, we can find a dense countable family $\{\psi_j\}_{j \geq 0} \subset B_{L^1}(1)$ such that

$$\sup_{j \geq 0} |\langle \psi_j, g \rangle| = \sup_{\|\psi\|_{L^1} \leq 1} |\langle \psi, g \rangle| = \|g\|_{L^\infty},$$

for any $g \in L^\infty$, which then implies that for a.e. $y \in \mathbb{T}^d$,

$$\sup_{j \geq 0} \left|\nabla_1^k \nabla_2^k \frac{\delta^2}{\delta\mu^2}\{\mathcal{F}[m^\mu]\}(\psi_j, y)\right| \lesssim \max_{k \leq 2} \sup_{\mu \in \mathcal{P}} \left\|\frac{\delta^k \mathcal{F}}{\delta\mu^k}\right\|_{L^\infty}.$$

The claim then follows from the duality again and the density of $\{\psi_k\}_{k \geq 0} \subset B_{L^1}(1)$. □

## 4.2 Proof of Main Theorem

The main Theorem of this paper is now stated as follows:

**Theorem 4.7.** *We assume that the initial data to (1.12) are prepared as in (3.1), and that Assumptions 1.3, 2.1 and 4.1 together with the stochastic-parabolicity condition (2.15) are in force.*

*We have then the following weak-error rate estimate of the approximation of the interacting particle system (1.1) by the Dean-Kawasaki equation (1.12) holds for any $t \in [0, T]$:*

$$|\mathbf{E}\mathcal{F}[\rho_t^N] - \mathbf{E}\mathcal{F}[m_t^{N,\epsilon}]| \lesssim_T \epsilon + \frac{\delta_N}{N} + \frac{\log(\epsilon^{-1})}{NL_N^2}, \quad (4.21)$$

*where constant does not depend on $N$ and $\epsilon$.*

**Proof.** For given $\mathcal{F}$ as specified in the statement and $t \in (0, T)$ fixed, we denote $\mathcal{G}^{\mathcal{F},t} := \mathcal{G}_s^{\mathcal{F},t}[\mu]$ as the unique classical solution to the backward equation

$$\begin{cases} 0 = \partial_s \mathcal{G}^{\mathcal{F},t}(s,\mu) + \int_{\mathbb{T}^d} \{a[\mu](x) : \nabla D_\mu \mathcal{G}^{\mathcal{F},t}(s,\mu,x) + b[\mu](x) \cdot D_\mu \mathcal{G}^{\mathcal{F},t}(s,\mu,x)\} d\mu(x), \\ \mathcal{G}^{\mathcal{F},t}(t,\mu) = \mathcal{F}[\mu], \end{cases} \quad (4.22)$$

for $s \in [0, t]$. By Itô's formula in Proposition 3.5, we have therefore

$$\mathcal{G}^{\mathcal{F},t}(t, \rho_t^N) - \mathcal{G}^{\mathcal{F},t}(0, \rho_0^N) = \frac{1}{N}\int_0^t \mathcal{V}_{\mathcal{G}_s^{\mathcal{F},t}}[\rho_s^N]ds + M_t^{N,\mathcal{F}},$$

where $M^{N,\mathcal{F}}$ is genuinely a martingale, and the dynamics cancel each other by design. Together with the similar Itô formula for $m^{N,\epsilon}$ from Proposition 3.6, this leads to the first-order expansion

$$\begin{cases} \mathbf{E}\mathcal{F}[\rho_t^N] = \mathbf{E}\mathcal{G}^{\mathcal{F},t}(t, \rho_t^N) = \frac{1}{N}\int_0^t \mathbf{E}\mathcal{V}_{\mathcal{G}_s^{\mathcal{F},t}}[\rho_s^N]ds + \mathbf{E}\mathcal{G}^{\mathcal{F},t}(0, \rho_0^N), \\ \mathbf{E}\mathcal{F}[m_t^{N,\epsilon}] = \mathbf{E}\mathcal{G}^{\mathcal{F},t}(t, m_t^{N,\epsilon}) = \frac{1}{N}\int_0^t \mathbf{E}\mathcal{V}_{\mathcal{G}_s^{\mathcal{F},t}}^N[m_s^{N,\epsilon}]ds + \mathbf{E}\mathcal{G}^{\mathcal{F},t}(0, m_0^{N,\epsilon}). \end{cases} \quad (4.23)$$

At this stage we just replace $\mathcal{V}_{\mathcal{G}_s^{\mathcal{F},t}}^N[\cdot]$ by $\mathcal{V}_{\mathcal{G}_s^{\mathcal{F},t}}[\cdot]$ and bound the difference by Lemma 3.7. We further denote

$$\mathcal{G}^{\mathcal{V},s} = \mathcal{G}_\tau^{\mathcal{V},s}[\mu]$$

with $\tau \in [0, s]$ and $\mu \in \mathcal{P}(\mathbb{T}^d)$, as the classical solution to the Kolmogorov equation

$$\partial_\tau \mathcal{G}^{\mathcal{V},s}(\tau, \mu) + \int_{\mathbb{T}^d}\{a[\mu](x):\nabla D_\mu \mathcal{G}^{\mathcal{V},s}(\tau, \mu, x) + b(\mu, x)\cdot D_\mu \mathcal{G}^{\mathcal{V},s}(\tau, \mu, x)\}d\mu(x) = 0, \quad (4.24)$$

for $\tau \in [0, s]$ with terminal data

$$\mathcal{G}^{\mathcal{V},s}(s, \mu) = \mathcal{V}_{\mathcal{G}_s^{\mathcal{F},t}}[\mu],$$

which exists and is unique thanks to Proposition 3.9 and the regularity of $\mathcal{V}_{\mathcal{G}_s^{\mathcal{F},t}}$ provided in Lemma 4.9 below. By Itô's formula for the SPDE and the particle systems, that is Proposition 3.5 and Proposition 3.6, together with the terminal conditions again, we have

$$\begin{aligned} \mathbf{E}\mathcal{V}_{\mathcal{G}_s^{\mathcal{F},t}}[\rho_s^N] - \mathcal{G}^{\mathcal{V},s}(0, \rho_0^N) &= \mathbf{E}\mathcal{G}^{\mathcal{V},s}(s, \rho_s^N) - \mathcal{G}^{\mathcal{V},s}(0, \rho_0^N) \\ &= \frac{1}{N}\mathbf{E}\int_0^s \int_{\mathbb{T}^d} a[\rho_\tau^N](x):(D_\mu^2 \mathcal{G}_\tau^{\mathcal{V},s})(\rho_\tau^N, x, x)d\rho_\tau^N(x)d\tau \end{aligned} \quad (4.25)$$

and

$$\begin{aligned} \mathbf{E}\mathcal{V}_{\mathcal{G}_s^{\mathcal{F},t}}[m_s^{N,\epsilon}] - \mathcal{G}^{\mathcal{V},s}(0, \rho_0^N) &= \mathbf{E}\mathcal{G}^{\mathcal{V},s}(s, m_s^{N,\epsilon}) - \mathcal{G}^{\mathcal{V},s}(0, m_0^{N,\epsilon}) \\ &= \frac{1}{N}\mathbf{E}\int_0^s \sum_{k\in\mathbb{Z}^d}|\theta_k^N|^2 \mathfrak{F}(D_\mu^2 \mathcal{G}_\tau^{\mathcal{V},s}[m_\tau^{N,\epsilon}]:\Lambda_\tau^{N,\epsilon})(k, -k)d\tau, \end{aligned} \quad (4.26)$$

where we abbreviate

$$\Lambda_\tau^{N,\epsilon}(x, y) := (\Sigma[m_\tau^{N,\epsilon}](x)f^N(m_\tau^{N,\epsilon})(x)) \cdot (\Sigma^t[m_\tau^{N,\epsilon}](y)f^N(m_\tau^{N,\epsilon})(y)).$$

Due to the regularity bound we have for the transport equation, the differences in (4.24) and (4.25) are bounded by

$$|\mathbf{E}\mathcal{G}^{\mathcal{V},s}(s, \rho_s^N) - \mathcal{G}^{\mathcal{V},s}(0, \rho_0^N)| \lesssim \frac{s}{N}\sup_{\tau\in[0,s]}\sup_{\mu\in\mathcal{P}}\|D_\mu^2 \mathcal{G}_\tau^{\mathcal{V},s}[\mu]\|_{L^\infty},$$

and (brutally)

$$\begin{aligned} |\mathbf{E}\mathcal{G}^{\mathcal{V},s}(s, m_s^{N,\epsilon}) - \mathcal{G}^{\mathcal{V},s}(0, m_0^{N,\epsilon})| &\leq \frac{s\|\theta^N\|_{\ell^2}^2}{N}\sup_{\tau\in[0,s]}\|D_\mu^2 \mathcal{G}_\tau^{\mathcal{V},s}[m_\tau^{N,\epsilon}]:\Lambda_\tau^{N,\epsilon}\|_{L^1} \\ &\lesssim \frac{s\|\theta^N\|_{\ell^2}^2}{N}\sup_{\tau\in[0,s]}\sup_{\mu\in\mathcal{P}}\|D_\mu^2 \mathcal{G}_\tau^{\mathcal{V},s}[\mu]\|_{L^\infty}, \end{aligned}$$

for any $s \in [0, t]$, where the constants depend only on $\Sigma[\cdot]$, since

$$\|f^N(m^{N,\epsilon})\|_{L^2} \lesssim \|m^{N,\epsilon}\|_{L^1}^{\frac{1}{2}} \leq 1$$

due to the conservation law and (1.14). Combining (4.23), (4.24) and (4.25), we have for $t \in [0, T]$ fixed above that $|\mathbf{E}\mathcal{F}[\rho_t^N] - \mathbf{E}\mathcal{F}[m_t^{N,\epsilon}]| = |\mathbf{E}\mathcal{G}^{\mathcal{F},t}(t, \rho_t^N) - \mathbf{E}\mathcal{G}^{\mathcal{F},t}(t, m_t^{N,\epsilon})|$, which is further bounded by

$$\begin{aligned}
|\mathbf{E}\mathcal{G}^{\mathcal{F},t}(t, \rho_t^N) - \mathbf{E}\mathcal{G}^{\mathcal{F},t}(t, m_t^{N,\epsilon})| &\lesssim |\mathcal{G}^{\mathcal{F},t}(0, \rho_0^N) - \mathcal{G}^{\mathcal{F},t}(0, m_0^{N,\epsilon})| \\
&\quad + \frac{1}{N}\int_0^t |\mathcal{G}^{\mathcal{V},s}(0, m_0^{N,\epsilon}) - \mathcal{G}^{\mathcal{V},s}(0, \rho_0^N)|\mathrm{d}s \\
&\quad + \frac{1}{N}\left|\int_0^t \mathbf{E}\mathcal{V}_{\mathcal{G}_s^{\mathcal{F},t}}[m_s^{N,\epsilon}] - \mathbf{E}\mathcal{V}_{\mathcal{G}_s^{\mathcal{F},t}}[m_s^{N,\epsilon}]\mathrm{d}s\right| \\
&\quad + \frac{t^2}{N^2}(1 + \|\theta^N\|_{\ell^2}^2) \sup_{s \in [0,t]} \sup_{\tau \in [0,s]} \sup_{\mu \in \mathcal{P}} \|D_\mu^2 \mathcal{G}_\tau^{\mathcal{V},s}[\mu]\|_{L^\infty}. \quad (4.27)
\end{aligned}$$

For the first two terms, we deduce from the Lipschitz estimate (3.24) with respect to the $\mathcal{W}_2$ distance of the solution to the Kolmogorov equation that

$$|\mathcal{G}^{\mathcal{F},t}(0, \rho_0^N) - \mathcal{G}^{\mathcal{F},t}(0, m_0^{N,\epsilon})| \leq \|\mathcal{G}_0^{\mathcal{F},t}\|_{\mathrm{lip}} \mathcal{W}_2(\rho_0^N, m_0^{N,\epsilon}) \lesssim \|\mathcal{F}\|_{\mathrm{lip}} \epsilon,$$

where the constant relies on the coefficients of the SDE (3.22) but *not* on $\mathcal{F}[\cdot]$. Similarly, since $\mathcal{G}^{\mathcal{V},s}$ is the solution to the Kolmogorov equation of the McKean–Vlasov system with terminal condition $\mathcal{G}_s^{\mathcal{V},s}[\mu] = \mathcal{V}_{\mathcal{G}_s^{\mathcal{F},t}}[\mu]$,

$$|\mathcal{G}^{\mathcal{V},s}(0, \rho_0^N) - \mathcal{G}^{\mathcal{V},s}(0, m_0^{N,\epsilon})| \leq \|\mathcal{G}_0^{\mathcal{V},s}\|_{\mathrm{lip}} \mathcal{W}_2(\rho_0^N, m_0^{N,\epsilon}) \lesssim \epsilon \|\mathcal{G}_0^{\mathcal{V},s}\|_{\mathrm{lip}},$$

which we can further bound by

$$\ldots \lesssim \epsilon \|\mathcal{V}_{\mathcal{G}_s^{\mathcal{F},t}}\|_{\mathrm{lip}} \lesssim \epsilon \|D_\mu^2 \mathcal{G}_s^{\mathcal{F},t}\|_{\mathrm{lip}} \lesssim \epsilon \|D_\mu^2 \mathcal{F}\|_{\mathrm{lip}},$$

with the constant is again independent of $\mathcal{F}[\cdot]$, where the first inequality is due to the stability estimate (3.24) again, the second due to (3.12). The last term in (4.26) is simply bounded by

$$\frac{t^2}{N^2}(1 + \|\theta^N\|_{\ell^2}^2) \sup_{s \in [0,t]} \sup_{\tau \in [0,s]} \sup_{\mu \in \mathcal{P}} \|D_\mu^2 \mathcal{G}_\tau^{\mathcal{V},s}[\mu]\|_{L^\infty} \lesssim \frac{t^2}{N^2}(1 + \|\theta^N\|_{\ell^2}^2) \sup_{k \leq 4} \sup_{\mu \in \mathcal{P}} \|D_\mu^k \mathcal{F}[\mu]\|_{L^\infty}.$$

Finally, we can bound the second last term in (4.26), i.e. the cost of the replacement of quadratic variation, according to Lemma 3.7 by

$$\left|\int_0^t \mathbf{E}[\mathcal{V}_{\mathcal{G}_s^{\mathcal{F},t}}^N(m_s^{N,\epsilon}) - \mathcal{V}_{\mathcal{G}_s^{\mathcal{F},t}}(m_s^{N,\epsilon})]\mathrm{d}s\right|$$
$$\lesssim t \sup_{\substack{\mu \in \mathcal{P} \cap L^1 \\ s \in [0,t]}} (\|D_\mu^2 \mathcal{G}_s^{\mathcal{F},t}[\mu]\|_{W^{2,\infty}}) \left[\delta_N + L_N^{-2}\left(1 + \mathbf{E}\int_0^t \int_{\mathbb{T}^d} \frac{|\nabla m_s^{N,\epsilon}|^2}{m_s^{N,\epsilon}} \mathrm{d}s\right)\right].$$

Recall from Proposition 3.9 that the solutions to the infinite-dimensional transport equations are given by characteristic lines, i.e.

$$D_\mu^2 \mathcal{G}_s^{\mathcal{F},t}[\mu] = D_\mu^2\{\mathcal{F}[m_{t-s}^\mu]\}, \quad s \in [0, t].$$

Due to the regularity estimate (4.19) in Corollary 4.6, we can control the former term now in the space of desired spatial regularity by

$$\sup_{s\in[0,t]} \|D_\mu^2 \mathcal{G}_s^{\mathcal{F},t}[\mu]\|_{W^{2,\infty}} \lesssim 1.$$

Together with the preparation of the initial data, the entropy bound (3.4) and $\|(k\theta_k^N)\|_{\ell^2}^2 \lesssim L_N^2$ by construction, we bound

$$\text{RHS. of (4.26)} \lesssim_{T,\mathcal{F}} \varepsilon\left(1 + \frac{1}{N}\right) + \frac{\delta_N}{N} + \frac{1}{NL_N^2}\left[\log\left(\frac{1}{\varepsilon}\right) + \frac{\|(k\theta_k^N)\|_{\ell_k^2}^2}{N}\right] + \frac{1 + \|\theta^N\|_{\ell^2}^2}{N^2}$$

$$\lesssim \varepsilon + \frac{\delta_N}{N} + \frac{1}{NL_N^2}\log\left(\frac{1}{\varepsilon}\right),$$

since $\|\theta^N\|_{\ell^2}^2/N^2 \lesssim \delta_N/N$ due to the parabolicity condition (2.15). We have finished the proof of the main theorem. □

**Remark 4.8.** Optimizing (4.21) in $\varepsilon$, $\delta_N$ and $L_N$ subject to the constraint given by the parabolicity condition (2.15) leads to the choice

$$\varepsilon \sim \frac{1}{N^2}, \quad \delta_N \sim N^{-\frac{1}{1+d/2}}, \quad L_N \sim \delta_N^{-\frac{1}{2}}, \tag{4.28}$$

the optimal rate is given by $N^{-1-2/(d+2)}\log N$. We have therefore recovered the same rate for the case of smooth mean-field interactions as that of the free particles in [DKP24].

The following Lemma has been applied in the proof of the main Theorem 4.7:

**Lemma 4.9.** *We assume that the test functional $\mathcal{F}[\cdot]$ and also the coefficients $b[\cdot]$ and $a[\cdot]$ have Lipschitz interior derivatives up to the fourth order in the sense of (A.15). Then we have that $\mathcal{V}_{\mathcal{G}_s^{\mathcal{F},t}}[\cdot]$ for any $0 \leq s \leq t \leq T$, which is defined through (3.11) and the equation (4.22), also has Lipschitz interior derivatives up to second order, and*

$$\sup_{s\in[0,t]} \sup_{\substack{x,y\in\mathbb{T}^d \\ \mu\in\mathcal{P}(\mathbb{T}^d)}} |D_\mu^2 \mathcal{V}_{\mathcal{G}_s^{\mathcal{F},t}}(\mu,x,y)| < \infty. \tag{4.29}$$

In other word, we can choose $\mathcal{V}_{\mathcal{G}_s^{\mathcal{F},t}}$ as the terminal condition for the Kolmogorov equation (3.20) for any $s\in[0,t]$.

**Proof.** By verifying definition (A.13), we have

$$D_\mu \mathcal{V}_{\mathcal{G}_s^{\mathcal{F},t}}[\mu](x) = \int_{\mathbb{T}^d} D_\mu(a[\mu](w):D_\mu^2\mathcal{G}_s^{\mathcal{F},t}(\mu,w,w))(x)\mathrm{d}\mu(w) + a[\mu](x):D_\mu^2\mathcal{G}_s^{\mathcal{F},t}(\mu,x,x),$$

where

$$D_\mu(a[\mu](w):D_\mu^2\mathcal{G}_s^{\mathcal{F},t}(\mu,w,w))(x) = D_\mu a[\mu](w;x):D_\mu^2\mathcal{G}_s^{\mathcal{F},t}(\mu,w,w)$$
$$+ a[\mu](w;x):D_\mu^3\mathcal{G}_s^{\mathcal{F},t}(\mu,w,w,x).$$

The second-order derivative can be derived identically, so we omit it for brevity, while the uniformly Lipschitz continuity of $D_\mu^2 \mathcal{V}_{\mathcal{G}_s^{\mathcal{F},t}}$ follows from the assumption on $a[\cdot]$ and $\mathcal{F}[\cdot]$ in Assumption 4.1 and the fact that under such an assumption, the solution $\mathcal{G}_s^{\mathcal{F},t}$ to the equation (4.22) has uniformly Lipschitz interior derivatives up to the fourth order.

Such a regularity result can be achieved by further iterating the linearization procedure used in Section 4.1, and notice that since we don't require higher spatial regularity as we required in (4.19) but only Lipschitz continuity, it follows directly from the main Theorem in [Tse21], which presents the construction of arbitrarily high order of iterations and the proof of Lipschitz continuity of the higher-order interior derivatives of the solution. □

# A  Analytic Setup

## A.1  Analytic Estimates

The proof of Lemma 2.4 relies on the following commutator estimate. We will sometimes use square-brackets to denote commutators of operators and identify the multiplication operators of some functions with the functions themselves slightly abusing the notation.

**Lemma A.1.** *For any $f \in W^{1+\varepsilon,\infty}(\mathbb{T}^d)$ and $g \in L^2(\mathbb{T}^d)$ for an arbitrarily small $\varepsilon > 0$, we have*

$$\|[f,(1-\Delta)^{\frac{1}{2}}]g\|_{L^2} \lesssim_d \|f\|_{L_W^{1+\varepsilon,\infty}} \|g\|_{L^2}, \tag{A.1}$$

*where the constant is independent of $f$ and $g$.*

**Proof.** We have according to the kernel representation of the fractional Laplacian, which can be found for example in [Sil07], that

$$\begin{aligned}
[f,(1-\Delta)^{\frac{1}{2}}]g(x) &= \left| \text{pv.} \int_{\mathbb{T}^d} (f(x)-f(y))\, g(y)\, K_d(x-y)\mathrm{d}y \right| \\
&= \left| \text{pv.} \int_{\mathbb{T}^d} \left( \int_0^1 \nabla f(y+t(x-y))\cdot(x-y)\mathrm{d}t \right) g(y) K_d(x-y)\mathrm{d}y \right| \\
&\leq \left| \text{pv.} \int_{\mathbb{T}^d} (\nabla f(y)\cdot(x-y))\, g(y)\, K_d(x-y)\mathrm{d}y \right| \\
&\quad + \int_{\mathbb{T}^d}\int_0^1 |\nabla f(y+t(x-y))-\nabla f(y)||g(y)||(x-y)K_d(x-y)|\mathrm{d}t\mathrm{d}y,
\end{aligned}$$

where $K_d$ on $\mathbb{R}^d$ is periodic and

$$K_d(x) = \frac{c_d}{|x|^{d+1}} + K_d^s(x), \quad x \in [-\tfrac{1}{2},\tfrac{1}{2}]^d \setminus \{0\},$$

for some bounded $K_d^s$ on $[-1/2,1/2]^d$, and (A.1) then follows, since we have

$$\text{pv.} \int_{\mathbb{T}^d} \nabla f(y) \cdot \frac{x-y}{|x-y|^d} g(y)\mathrm{d}y = \sum_{j=1}^d R_j[(\partial_j f)\, g], \tag{A.2}$$

where $(R_j)_{j=1}^d$ denote the Riesz transforms which are bounded on $L^2(\mathbb{T}^d)$ (c.f. [Gra08], Chapter 5), while

$$\int_{\mathbb{T}^d}\int_0^1 \frac{|\nabla f(y+t(x-y))-\nabla f(y)|}{|x-y|^\varepsilon} |g(y)||x-y|^{1+\varepsilon}|K_d(x-y)|\mathrm{d}t\mathrm{d}y$$
$$\lesssim \int_0^1 t^\varepsilon \mathrm{d}t \, \|\nabla f\|_{W_{\mathbb{T}^d}^{\varepsilon,\infty}} \int_{\mathbb{T}^d} |g(y)||x-y|^{1+\varepsilon}|K_d(x-y)|\mathrm{d}y.$$

The claim follows from Young's inequality and the bound $|\cdot|^{1+\varepsilon}|K_d(\cdot)| \in L^1(\mathbb{T}^d)$. □

We are now ready to prove the coercivtiy bound in $W^{-1,2}(\mathbb{T}^d)$.

**Proof. (Lemma 2.4)** The first claim follows directly from the bound that

$$-{}_{(L^2)^*}\langle \nabla \cdot (a[v]\nabla u), w\rangle_{L^2} = \langle a[v]\nabla u, \nabla w\rangle_{W^{-1,2}} \leq \|a[v]\|_{W^{1,\infty}}\|u\|_{L^2}\|w\|_{L^2},$$

since multiplication is a well-defined continuous map from $W^{1,\infty} \times W^{-1,2}$ to $W^{-1,2}$. We now focus on the coercivity bound (2.12). By definition, we have

$$\begin{aligned}
-{}_{(L^2)^*}\langle \nabla \cdot (a[v]\nabla u), u\rangle_{L^2} &= -\left\langle (1-\Delta)^{-\frac{1}{2}}\nabla \cdot (a[v]\nabla u), (1-\Delta)^{-\frac{1}{2}}u\right\rangle_{L^2} \\
&= \left\langle a[v](1-\Delta)^{-\frac{1}{2}}\nabla u, (1-\Delta)^{-\frac{1}{2}}\nabla u\right\rangle_{L^2} \\
&\quad -\left\langle \left[(1-\Delta)^{\frac{1}{2}}, a[v]\right](1-\Delta)^{-\frac{1}{2}}\nabla u, (1-\Delta)^{-1}\nabla u\right\rangle_{L^2},
\end{aligned}$$

where the first term is lower bounded by

$$\begin{aligned}
\left\langle a[v](1-\Delta)^{-\frac{1}{2}}\nabla u, (1-\Delta)^{-\frac{1}{2}}\nabla u\right\rangle_{L^2} &\geq a_{\min}\left\|(1-\Delta)^{-\frac{1}{2}}\nabla u\right\|_{L^2}^2 \\
&= a_{\min}\sum_{k\in\mathbb{Z}^d}\frac{|2\pi k|^2}{1+|2\pi k|^2}|\hat{u}(k)|^2 \geq \frac{a_{\min}}{2}\sum_{k\in\mathbb{Z}^d\setminus\{0\}}|\hat{u}(k)|^2,
\end{aligned}$$

while for the second term we abbreviate the commutator

$$K[v] := \left[a[v], (1-\Delta)^{\frac{1}{2}}\right]$$

and derive

$$\begin{aligned}
\left|\left\langle K[v](1-\Delta)^{-\frac{1}{2}}\nabla u, (1-\Delta)^{-1}\nabla u\right\rangle_{L^2}\right| &\overset{(A.1)}{\leq} C\|a[v]\|_{W^{1+\epsilon,\infty}}\|u\|_{L^2}\|u\|_{W^{-1,2}} \\
&\leq \frac{a_{\min}}{4}\|u\|_{L^2}^2 + \frac{C^2\|a[v]\|_{W^{1+\epsilon,\infty}}^2}{a_{\min}}\|u\|_{W^{-1,2}}^2,
\end{aligned}$$

for some generic $C > 0$ depending only on $d$. Combing the estimates above, we arrive at (2.12), since $|\hat{u}(0)|^2 \leq \|u\|_{W^{-1,2}}^2$ by definition (1.4). □

**Lemma A.2.** *Given $\kappa_\zeta = \frac{1}{\zeta}\kappa(\frac{\cdot}{\zeta})$, for fixed kernel $\kappa \in C_c^\infty(\mathbb{R})$ and $\zeta > 0$ s.t. $\mathrm{supp}(\kappa) \in (-1, 1)$, we have for any $m \in W^{1,2}(\mathbb{T}^d)$ that*

$$\lim_{\zeta\to 0^+}\int_{\mathbb{T}^d}\kappa_\zeta(m(x))|m(x)|^2 dx = 0,$$
$$\lim_{\zeta\to 0^+}\int_{\mathbb{T}^d}\kappa_\zeta(m(x))|m(x)\nabla m(x)| dx = 0.$$

**Proof.** For the first limit, it suffices to notice that

$$\sup_{l\in\mathbb{R}}[l\kappa_\zeta(l)] = \sup_{l\in\mathbb{R}}[l\kappa(l)] \lesssim_\kappa 1, \tag{A.3}$$

which gives

$$\int_{\mathbb{T}^d}\kappa_\zeta(m(x))|m(x)|^2 dx \lesssim_\kappa \int_{\{x:|m(x)|\leq\zeta\}}|m(x)| dx \to 0, \quad \text{as}\quad \zeta \to 0,$$

since $\mathbb{T}^d$ is compact. For the second one, we have similarly

$$\int_{\mathbb{T}^d}\kappa_\zeta(m(x))|m(x)\nabla m(x)| dx \lesssim_\kappa \int_{\{x:|m(x)|\leq\zeta\}}|\nabla m(x)| dx \to \int_{\{x:m(x)=0\}}|\nabla m(x)| dx = 0,$$

as $\zeta \to 0$, where we have used the integrability of $\nabla m$ for the limit and Stampacchia's Lemma ([Eva22], Chapter 5, Exercises 18 (iii)) for the equality, which asserts that the distributional derivative $\nabla u$ vanishes a.e. on the set $\{m = 0\}$ for any $m \in W^{1,2}$. □

The following is the main technical lemma that we apply for the approximation of generators.

**Lemma A.3.** *Given functions $f, g \in C^\infty(\mathbb{T}^d, \mathbb{R})$ and $h \in C^\infty(\mathbb{T}^d \times \mathbb{T}^d, \mathbb{R})$, we have the following quantitative estimate*

$$\left| \sum_{k \in \mathbb{Z}^d} |\hat{\varphi}(k)|^2 \mathfrak{F}(f \otimes g\, h)(k, -k) - \int_{\mathbb{T}^d} (f \otimes g\, h)(x, x) \mathrm{d}x \right|$$
$$\lesssim \|f\|_{W^{1,2}(\mathbb{T}^d)} \|g\|_{W^{1,2}(\mathbb{T}^d)} \|h\|_{W^{2,\infty}(\mathbb{T}^{2d})} \int_{\mathbb{T}^d} |\xi|^2 |\psi(\xi)| \mathrm{d}\xi, \qquad (A.4)$$

*where $f \otimes g\, h$ on $\mathbb{T}^d \times \mathbb{T}^d$ is given by $(f \otimes g\, h)(x, y) := h(x, y) f(x) g(y)$, and $\varphi \in C^\infty(\mathbb{T}^d, \mathbb{R})$, s.t. $\hat{\varphi}(0) = 1$, and*

$$\psi(x) := \sum_{k \in \mathbb{Z}^d} |\hat{\varphi}(k)|^2 e^{-2\pi i k \cdot x}. \qquad (A.5)$$

**Proof.** We first transform the sum of Fourier coefficients on the left of (A.3) to the form of integration on the diagonal of $\mathbb{T}^d \times \mathbb{T}^d$. We notice that $\psi(x) = \psi(-x)$, since $\varphi$ is real, and therefore $\hat{\varphi}(-k) = \overline{\hat{\varphi}(k)}$. We further denote

$$\Phi(x, y) := \sum_{k \in \mathbb{Z}^d} |(\mathfrak{F}_{\mathbb{T}^d} \varphi)(k)|^2 e^{-2\pi i k \cdot (x-y)} = \psi(x - y), \qquad (A.6)$$

for which we have

$$\int_{\mathbb{T}^{2d}} \Phi = \mathfrak{F}_{\mathbb{T}^{2d}}(\Phi)(0) = |(\mathfrak{F}_{\mathbb{T}^d} \varphi)(0)|^2 = 1,$$

and as desired,

$$\int_{\mathbb{T}^d} \Phi * (f \otimes g\, h)(x, x) \mathrm{d}x = \sum_{k \in \mathbb{Z}^d} |(\mathfrak{F}_{\mathbb{T}^d} \varphi)(k)|^2 \mathfrak{F}_{\mathbb{T}^d \times \mathbb{T}^d}(f \otimes g\, h)(k, -k). \qquad (A.7)$$

For brevity, we omit the differentials in the integrals when all variables are integrated. By identifying functions $C^\infty(\mathbb{T}^d)$ with smooth periodic functions on $\mathbb{R}^d$, we can then rewrite and bound the L.H.S. of (A.3) by

$$\left| \int_{\mathbb{T}^d} \Phi * (f \otimes g\, h)(x, x) \mathrm{d}x - \int_{\mathbb{T}^d} (f \otimes g\, h)(x, x) \mathrm{d}x \right|$$
$$= \left| \int_{[-\frac{1}{2}, \frac{1}{2}]^{3d}} \Phi(y, z) \left[ (f \otimes g\, h)(x - y, x - z) - (f \otimes g\, h)(x, x) \right] \right|$$
$$= \left| \int_{[-\frac{1}{2}, \frac{1}{2}]^{3d}} \psi(y - z) \left[ (f \otimes g\, h)(x - y + z, x) - (f \otimes g\, h)(x, x) \right] \right|$$
$$= \left| \int_{[-\frac{1}{2}, \frac{1}{2}]^{2d}} \psi(y) \left[ (f \otimes g\, h)(x - y, x) - (f \otimes g\, h)(x, x) \right] \right|,$$

where we have used the invariance of the Lebesgue measure on $\mathbb{T}^{3d}$ under the push-forward of translation. We can split the difference by

$$(f \otimes g\, h)(x - y, x) - (f \otimes g\, h)(x, x) = \{[h(x - y, x) - h(x, x)] f(x)$$
$$+ [h(x - y, x) - h(x, x)][f(x - y) - f(x)] + h(x, x)[f(x - y) - f(x)]\} \cdot g(x), \qquad (A.8)$$

and consider them separately. For the first term in (A.7), we have

$$\left| \int_{[-\frac{1}{2},\frac{1}{2}]^{2d}} \psi(y)[h(x-y,x)-h(x,x)]f(x)g(x) \right|$$

$$= \left| \frac{1}{2} \int_{[-\frac{1}{2},\frac{1}{2}]^{2d}} \psi(y)[h(x-y,x)+h(x+y,x)-2h(x,x)]f(x)g(x) \right|$$

$$\lesssim \left( \int_{\mathbb{T}^d} |y|^2 |\psi(y)| \mathrm{d}y \right) \|h\|_{W^{2,\infty}} \|f\|_{L^2} \|g\|_{L^2},$$

where we use symmetry $\psi(-\cdot)=\psi$ and change of variable $y \mapsto -y$ for the equality, and

$$\sup_{x \in \mathbb{T}^d} |h(x-y,x)+h(x+y,x)-2h(x,x)| \lesssim \|D^2 h\|_{L^\infty} |y|^2$$

for the inequality. For the second term in (A.7), we have

$$\left| \int_{[-\frac{1}{2},\frac{1}{2}]^{2d}} \psi(y)[h(x-y,x)-h(x,x)][f(x-y)-f(x)]g(x) \right|$$

$$= \left| \int_{[-\frac{1}{2},\frac{1}{2}]^{2d}} \psi(y) g(x) \int_0^1 (y,0) \cdot Dh(x-sy,x) \mathrm{d}s \int_0^1 y \cdot Df(x-ty) \mathrm{d}t \right|$$

$$\lesssim \left( \int_{\mathbb{T}^d} |y|^2 |\psi(y)| \mathrm{d}y \right) \|h\|_{W^{1,\infty}} \|f\|_{W^{1,2}} \|g\|_{L^2},$$

while for the last term in (A.7), we abbreviate $\tilde{g}(x) := h(x,x)g(x)$ for $x \in \mathbb{T}^d$ and derive

$$\left| \int_{[-\frac{1}{2},\frac{1}{2}]^{2d}} \psi(y)[f(x-y)-f(x)] \tilde{g}(x) \right|$$

$$= \left| \int_{[-\frac{1}{2},\frac{1}{2}]^d} \psi(y) y \cdot \left( \int_0^1 \int_{[-\frac{1}{2},\frac{1}{2}]^d} \tilde{g}(x) Df(x-ty) \mathrm{d}x \mathrm{d}t \right) \mathrm{d}y \right|$$

$$= \left| \int_{[-\frac{1}{2},\frac{1}{2}]^d} \psi(y) y \cdot \left( \int_0^1 \int_{[-\frac{1}{2},\frac{1}{2}]^d_{-ty}} [\tilde{g}(x+ty)-\tilde{g}(x)] Df(x) \mathrm{d}x \mathrm{d}t \right) \mathrm{d}y \right|$$

$$\leq \left( \int_{\mathbb{T}^d} |y|^2 |\psi(y)| \mathrm{d}y \right) \|h\|_{W^{1,\infty}} \|g\|_{W^{1,2}} \|f\|_{W^{1,2}},$$

where we use for the second equality the change of variable $x \mapsto x-ty$ and

$$\int_{[-\frac{1}{2},\frac{1}{2}]^d} \psi(y) y \cdot \left( \int_0^1 \int_{[-\frac{1}{2},\frac{1}{2}]^d_{-ty}} \tilde{g}(x) Df(x) \right) = \left( \int_{[-\frac{1}{2},\frac{1}{2}]^d} \psi(y) y \right) \cdot \left( \int_{[-\frac{1}{2},\frac{1}{2}]^d} \tilde{g}(x) Df(x) \right) = 0,$$

due to the symmetry $\psi(-\cdot)=\psi$ again, in order to smuggle in the integral with $\tilde{g}(x)$. The last inequality follows from

$$\left| \int_0^1 \int_{[-\frac{1}{2},\frac{1}{2}]^d_{-ty}} [\tilde{g}(x+ty)-\tilde{g}(x)] Df(x) \right|$$

$$= \left| \int_{[0,1]^2} \int_{[-\frac{1}{2},\frac{1}{2}]^d_{-ty}} [ty \cdot D\tilde{g}(x+sty)] Df(x) \mathrm{d}x \mathrm{d}t \mathrm{d}s \right|$$

$$\leq \|Df\|_{L^2} \int_{[0,1]^2} t|y| \left( \int_{[-\frac{1}{2},\frac{1}{2}]^d_{-ty}} |D\tilde{g}(x+sty)|^2 \mathrm{d}x \right)^{\frac{1}{2}} \mathrm{d}t \mathrm{d}s,$$

where by periodicity again

$$\left( \int_{[-\frac{1}{2},\frac{1}{2}]^d_{-ty}} |D\tilde{g}(x+sty)|^2 \mathrm{d}x \right)^{\frac{1}{2}} = \|D\tilde{g}\|_{L^2} \lesssim \|g\|_{W^{1,2}} \|h\|_{W^{1,\infty}}.$$

Summarizing all the bounds we derive for the terms in (A.7), we arrive at (A.3). □

The trilinear estimate can be extended to the whole matrix-valued Sobolev spaces by a density argument, if we assume that $\hat{\varphi}$ is compactly supported. The statement can be sharpened by replacing the compact support condition by integrability of the Fourier coefficients, but the former suffices for our purpose.

**Corollary A.4.** *If $\hat{\varphi}$ is compactly supported, then it holds for all $f \in W^{1,2}(\mathbb{T}^d, \mathbb{R}^{d \times d})$ and $h \in W^{2,\infty}(\mathbb{T}^{2d}, \mathbb{R}^{d \times d})$ that*

$$\left| \sum_{k \in \mathbb{Z}^d} |\hat{\varphi}(k)|^2 \mathfrak{F}\Lambda(k, -k) - \int_{\mathbb{T}^d} \Lambda(x, x) \, dx \right| \lesssim \|f\|^2_{W^{1,2}(\mathbb{T}^d)} \|h\|_{W^{2,\infty}(\mathbb{T}^{2d})} \int_{\mathbb{T}^d} |\xi|^2 |\psi(\xi)| \, d\xi, \tag{A.9}$$

*where $\Lambda(x, y) := [f(x) \cdot f^t(y)] : h(x, y)$.*

**Proof.** If suffices to notice that for smooth approximation $f_\epsilon \to f$, $g_\epsilon \to g$ in $L^2(\mathbb{T}^d)$ and $h_\epsilon \to h$ in $L^\infty(\mathbb{T}^{2d})$, we have by Hölder inequality that in $L^1(\mathbb{T}^{2d})$, $f_\epsilon \otimes g_\epsilon h_\epsilon \to f \otimes g h$, which implies both the approximation of the finite Fourier coefficients in $\sum_{k \in \mathbb{Z}^d} |\hat{\varphi}(k)|^2 \mathfrak{F}(f \otimes g h)(k, -k)$, and

$$\int_{\mathbb{T}^{2d}} f_\epsilon \otimes g_\epsilon h_\epsilon \to \int_{\mathbb{T}^{2d}} f \otimes g h. \qquad \square$$

## A.2 Moment Bounds of (2.13)

In this section, we give a concise proof of the moment bounds of the $L^2$ solution $m^{N,\chi}$ to the tamed and regularized equation with drift and noise terms given by (2.13), that is for any $T > 0$ and $p \in [2, \infty)$,

$$\mathbf{E}\left( \sup_{t \leq T} \|m^{N,\chi}_t\|^{2p}_{L^2(\mathbb{T}^d)} \right) < \infty, \tag{A.10}$$

*assuming* that the conditions (H2) and (H5) hold for (2.13) in the Gelfand triple $W^{1,2} \subset L^2 \subset W^{-1,2}$.

**Remark A.5.** We emphasize that the reason for introducing the $(L^2)^* \subset W^{-1,2} \subset L^2$ Gelfand triple in the first place is mainly to tackle the local monontonicity condition (H3), while the other conditions in the $L^2$ Gelfand triple can be proven straightforwardly. Such moment bounds are used in the proof of the non-negativity of the solutions in Corollary 2.7. For the general case, we refer to Theorem 2.6 in [RSZ24] and [LR15].

We recall that $\|\cdot\|_{L_2}$ is used to denote the Hilbert-Schmidt norm of operators on $L^2(\mathbb{T}^d)$. We start with the Krylov-Itô formula applied on the solution which we abbreviate as $m := m^{N,\chi}$, i.e.

$$\begin{aligned}
\frac{1}{2}\|m_t\|^{2p}_{L^2} &= p \int_0^t \|m_s\|^{2(p-1)}_{L^2} \left( \langle A^\chi[m_s], m_s \rangle + \frac{1}{2}\|B^N[m_s]\|^2_{L_2} \right) ds \\
&\quad + p(p-1) \int_0^t \|m_s\|^{2(p-2)}_{L^2} [\langle m_\cdot, B^N[m_\cdot] \cdot dW_\cdot \rangle]_s ds \\
&\quad + \frac{1}{2}\|m_0\|^{2p}_{L^2} + p \int_0^t \|m_s\|^{2(p-1)}_{L^2} \langle m_s, B^N[m_s] \cdot dW_s \rangle \\
&\leq p \int_0^t \|m_s\|^{2(p-1)}_{L^2} \left( -\frac{a_{\min}}{2}\|m_s\|^2_{W^{1,2}} + C\|m_s\|^2_{L^2} \right) ds \\
&\quad + p(p-1) \int_0^t \|m_s\|^{2(p-2)}_{L^2} [\langle m_\cdot, B^N[m_\cdot] \cdot dW_\cdot \rangle]_s ds \\
&\quad + \frac{1}{2}\|m_0\|^{2p}_{L^2} + p \int_0^t \|m_s\|^{2(p-1)}_{L^2} \langle m_s, B^N[m_s] \cdot dW_s \rangle, \tag{A.11}
\end{aligned}$$

where the last equality in (A.5) is due to the coercivity condition of (2.13) in the Gelfand triple ($W^{1,2} \subset L^2 \subset W^{-1,2}$). Meanwhile, we have

$$[\langle m_\cdot, B^N[m_\cdot] \cdot dW_\cdot \rangle]_s \leq \|m_s\|_{L^2}^2 \|B^N[m_s]\|_{L_2}^2 ds,$$

by Itô Isometry, while *if* we assume that the last term in (A.5) is a martingale so that BDG inequality gives

$$\mathbf{E}\left(\sup_{t \leq T}\left|\int_0^t \|m_s\|_{L^2}^{2(p-1)} \langle m_s, B^N[m_s] \cdot dW_s \rangle\right|\right) \leq C\mathbf{E}\left(\int_0^T \|m_s\|_{L^2}^{4p-2} \|B^N[m_s]\|_{L_2}^2 ds\right)^{\frac{1}{2}}$$

$$\leq \frac{1}{2}\mathbf{E}\sup_{t \leq T}\|m_t\|_{L^2}^{2p} + C\mathbf{E}\left(\int_0^T \|m_s\|_{L^2}^{2p-2} \|B^N[m_s]\|_{L_2}^2 ds\right),$$

it follows from the (H5) condition that

$$\|B^N[m_s]\|_{L_2}^2 \lesssim 1 + \|m_s\|_{L^2}^2. \tag{A.12}$$

We can take expectation and rearrange (A.5) by

$$\mathbf{E}\left(\sup_{t \leq T} \|m_t\|_{L^2}^{2p}\right) \lesssim \mathbf{E}\int_0^T \|m_s\|_{L^2}^{2p} ds + \mathbf{E}\|m_0\|_{L^2}^{2p},$$

and then conclude that (A.10) holds with Grönwall inequality. The general case is tackled with a canonical stopping-time argument.

## A.3 Calculus on the Space of Probability Measures

In this section, we give a short review of the *Lions calculus* that we use on the probability measures to specify notions that we use. The results are more or less standard, and we refer to [Car10, BLPR17, CDLL19, CD22] for detailed proofs and discussions.

### A.3.1 Definitions

Throughout, we equip the space of probability measures on $\mathbb{T}^d$ with weak topology, or equivalently, with the Wasserstein distance $\mathcal{W}_2(\cdot,\cdot)$.

**Definition A.6.** *Given $\mathcal{F}: \mathcal{P}(\mathbb{T}^d) \to \mathbb{R}$ continuous, we define the directional derivative $\delta\mathcal{F}/\delta\mu: \mathcal{P}(\mathbb{T}^d) \times \mathbb{T}^d \to \mathbb{R}$ to be continuous and satisfy*

$$\lim_{\epsilon \to 0^+} \frac{1}{\epsilon}\{\mathcal{F}[\epsilon m' + (1-\epsilon)m] - \mathcal{F}[m]\} = \int_{\mathbb{T}^d} \frac{\delta\mathcal{F}}{\delta\mu}[m](x) d(m'-m)(x), \tag{A.13}$$

*for any $m, m' \in \mathcal{P}$, and*

$$\int_{\mathbb{T}^d} \frac{\delta\mathcal{F}}{\delta\mu}[m](x) dm(x) = 0. \tag{A.14}$$

*If $\delta\mathcal{F}/\delta\mu: \mathcal{P}(\mathbb{T}^d) \times \mathbb{T}^d \to \mathbb{R}$ is further continuously differentiable with respect to the spatial variable, we define the interior derivative to be*

$$D_\mu \mathcal{F}[m](x) := \nabla_x \frac{\delta\mathcal{F}}{\delta\mu}[m](x). \tag{A.15}$$

Since in the end we will only need the interior derivatives for the analysis, the normalization condition (A.14) is not essential. By fixing $m', m \in \mathcal{P}(\mathbb{T}^d)$, one can consider the map from $[0,1] \to \mathbb{R}$,

$$s \mapsto \mathcal{F}[sm' + (1-s)m] = \mathcal{F}[m + s(m'-m)],$$

which can be shown by elementary argument to be continuously differentiable at any point $s \in (0,1)$ following from (A.13). Furthermore,

$$\mathcal{F}[m'] - \mathcal{F}[m] = \int_0^1 \int_{\mathbb{T}^d} \frac{\delta \mathcal{F}}{\delta \mu}[sm' + (1-s)m](x) \mathrm{d}(m'-m)(x) \mathrm{d}s. \tag{A.16}$$

We assign any functional $\tilde{\mathcal{F}}[\cdot]$ defined on $\mathcal{P}(\mathbb{T}^d) \times \mathbb{T}^{d \times n}$ for $n \in \mathbb{N}$ the Lipschitz constant

$$\begin{aligned}\|\tilde{\mathcal{F}}\|_{\mathrm{lip}} &= \sup_{\mu \in \mathcal{P}} \sup_{x \in \mathbb{T}^{d \times n}} |\tilde{\mathcal{F}}[\mu](x)| \\ &\quad + \sup_{\mu \neq \nu \in \mathcal{P}} \sup_{x \neq y \in \mathbb{T}^{d \times n}} \left\{ \frac{|\tilde{\mathcal{F}}[\mu](x) - \tilde{\mathcal{F}}[\nu](y)|}{|x-y| + \mathcal{W}_2(\mu,\nu)} \right\}.\end{aligned} \tag{A.17}$$

It follows that

$$\sup_{\mu \in \mathcal{P}(\mathbb{T}^d)} \|\tilde{\mathcal{F}}[\mu]\|_{W^{1,\infty}} \lesssim \|\tilde{\mathcal{F}}\|_{\mathrm{lip}}. \tag{A.18}$$

We now briefly mention the fundamental observation by Lions [Lio], which bridges the probabilistic argument for example used in [BLPR17] and the PDE counterpart and also motivates heuristically the consideration of interior derivatives in our context. For our consideration, the PDE arguments provide a direct grasp on the spatial regularity of the solutions, but we believe the two perspectives should be more or less equivalent.

Given $\mathcal{F}: \mathcal{P}(\mathbb{T}^d) \to \mathbb{T}^d$ and a (atomless and Polish) probability space $(\Omega, \mathbf{P})$, we define the lift $\tilde{\mathcal{F}}$ of $\mathcal{F}$ onto $L^2(\Omega; \mathbb{R}^d)$ to be

$$L^2(\Omega, \mathbb{R}^d) \ni X \mapsto \tilde{\mathcal{F}}(X) := \mathcal{F}[\mathcal{L}(\tau_X(0))], \tag{A.19}$$

where $(\tau_x)_{x \in \mathbb{R}^d}$ is the group of translations on $\mathbb{T}^d$ and $x \mapsto \tau_x(0)$ is the canonical projection from $\mathbb{R}^d$ to $\mathbb{T}^d$. In other word, $\tau_X(0)$ is the periodization of $X$. For any $\mu \in \mathcal{P}(\mathbb{T}^d)$, if there exists $X_0 \in L^2(\Omega; \mathbb{R}^d)$ with

$$\mathcal{L}(\tau_{X_0}(0)) = \mu,$$

such that $\tilde{\mathcal{F}}$ is *Fréchét differentiable* at $X_0$, then it is *Fréchét differentiable* at any $X \in L^2(\Omega; \mathbb{R}^d)$ with $\mathcal{L}(\tau_X(0)) = \mu$ and the law of the derivative (identified as an element in $L^2(\Omega; \mathbb{R}^d)$) doesn't depend on the choice $X$. If we assume that the Fréchét derivative at any $X \in L^2(\Omega; \mathbb{R}^d)$, denoted as

$$D\tilde{\mathcal{F}}[X] \in L^2(\Omega; \mathbb{R}^d),$$

exists and is continuous with respect to $X$, then for any $\mu \in \mathcal{P}(\mathbb{T}^d)$, there exists a lift $\xi_\mu \in L^2_\mu(\mathbb{T}^d; \mathbb{R}^d)$ s.t.

$$D\tilde{\mathcal{F}}[X] = \xi_\mu(\tau_X(0)), \quad \mathbf{P}\text{-a.s.}$$

for any $X \in L^2(\Omega; \mathbb{R}^d)$ with $\mathcal{L}(\tau_X(0)) = \mu$. A detailed proof can be found in [Car10]. The Lions (interior) derivative now comes into play, since it is just the lift function given above, in the sense that for any continuous map $\mathcal{F}$ on $\mathcal{P}(\mathbb{T}^d)$ s.t. $D_\mu \mathcal{F}$ exists and is continuous, $\tilde{\mathcal{F}}$ defined in (A.19) is *Fréchét continuously differentiable* on $L^2(\Omega, \mathbb{R}^d)$, and we have

$$\xi_\mu = D_\mu \mathcal{F}[\mu], \tag{A.20}$$

Conversely, if $\tilde{\mathcal{F}}$ defined is *Fréchét continuously differentiable* at any $X \in L^2(\Omega, \mathbb{R}^d)$, $D_\mu \mathcal{F}$ then exists and (A.20) holds as well. The proof can be found in [CDLL19].

### A.3.2 Inner Mollification with Fourier Coefficients

In this section, we construct the mollification of functionals on $\mathcal{P}(\mathbb{T}^d)$ claimed in Proposition 3.4. We recall that $\varphi^M$ is the Fejér kernel of order $M$ as defined in (3.5) and we denote
$$\mathcal{P}_{\leq M}(\mathbb{T}^d) := \{\mu \in \mathcal{P}(\mathbb{T}^d) : \hat{\mu}(k) = 0, |k|_\infty > M\}.$$
There exists a bijection
$$\mathscr{I}^M : G_M \to \mathcal{P}_{\leq M}(\mathbb{T}^d), \quad (a_\varphi)_{\varphi \in \mathscr{D}_M} \mapsto 1 + \sum_{\varphi \in \mathscr{D}_M} a_\varphi \varphi,$$
where the zeroth Fourier mode is always one, since we are considering probability measures. $G_M$ is a closed subset of $\mathbb{R}^{\mathscr{D}_M}$ as a consequence of Bochner Theorem [Boc33], and
$$\mathscr{D}_M := \{\sqrt{2}\cos(2\pi k \cdot), \sqrt{2}\sin(2\pi k \cdot)\}_{|k|_\infty < M, k \in \mathbb{Z}^d_+},$$
where we recall that $\mathbb{Z}^d_+$ contains those non-zero lattice points whose components are non-negative, and
$$[(\mathscr{I}^M)^{-1}\mu]_\varphi = \int_{\mathbb{T}^d} \varphi \, d\mu, \quad \varphi \in \mathscr{D}_M.$$
We have apparently for any $k \in \mathbb{Z}^d_+$ that $\sqrt{2}\operatorname{Re}\hat{\mu}(k) = [(\mathscr{I}^M)^{-1}\mu]_\varphi$ for $\varphi = \sqrt{2}\cos(2\pi k \cdot)$, and $\sqrt{2}\operatorname{Im}\hat{\mu}(k) = [(\mathscr{I}^M)^{-1}\mu]_{\bar{\varphi}}$ for $\bar{\varphi} = \sqrt{2}\sin(2\pi k \cdot)$, where Re and Im are real and imaginary parts respectively. Due to the fact that[A.1]
$$\mathcal{O}_M := (\mathscr{I}^M)^{-1}\{\mu \in \mathcal{P}_{\leq M}(\mathbb{T}^d) : \mu > 0\}$$
is open in $\mathbb{R}^{\mathscr{D}_M}$, there exists a ball $B_M$ centered at zero (zero corresponding to uniform distribution on $\mathbb{T}^d$) of sufficiently small radius (depending on $M$) in $\mathbb{R}^{\mathscr{D}_M}$ that is contained in $\mathcal{O}_M \subset G_M$. We now construct the approximation of $\mathcal{F}$ on $\mathcal{P}(\mathbb{T}^d)$ as
$$\begin{aligned}\mathcal{F}^M[\mu] &= \int_{\mathbb{R}^{\mathscr{D}_M}} \mathcal{F}\left[(1-\frac{1}{M})\mu * \varphi^M + \frac{1}{M}\mathscr{I}_M(z)\right]\rho_M(z)\,dz \\ &= \int_{\mathbb{R}^{\mathscr{D}_M}} \mathcal{F} \circ \mathscr{I}^M\left((\mathscr{I}^M)^{-1}\left[(1-\frac{1}{M})\mu * \varphi^M + \frac{1}{M}\mathscr{I}_M(z)\right]\right)\rho_M(z)\,dz, \end{aligned} \quad (A.21)$$
where $(\rho_M)_{M \geq 1}$ is an approximate identity with smooth, radially symmetric and positive functions supported on $B_M$[A.2]. We used in (A.1) the fact that
$$(1-\frac{1}{M})\mu * \varphi^M + \frac{1}{M}\mathscr{I}_M(z) \in \mathcal{P}_{\leq M}(\mathbb{T}^d),$$
for any $z \in B_M$, and we can further write
$$\begin{aligned}\left[(\mathscr{I}^M)^{-1}\left(\frac{M-1}{M}\mu * \varphi^M + \frac{1}{M}\mathscr{I}_M(z)\right)\right]_\varphi &= (\frac{M-1}{M})\operatorname{Re}[\hat{\mu}(k)\hat{\varphi}^M(k)] + \frac{1}{M}z_\varphi \\ &= \frac{M-1}{\sqrt{2}M}\{[(\mathscr{I}^M)^{-1}\mu]_\varphi \operatorname{Re}[\hat{\varphi}^M(k)] - [(\mathscr{I}^M)^{-1}\mu]_{\bar{\varphi}}\operatorname{Im}[\hat{\varphi}^M(k)]\} + \frac{1}{M}z_\varphi,\end{aligned}$$

---
A.1. $\mu > 0$ in sense that $\mu$ has strictly positive density.

A.2. In fact, one could immediately extend $\mathcal{F} \circ \mathscr{I}^M$ to a continuous function on $\mathbb{R}^{\mathscr{D}_M}$, since $\mathcal{F} \circ \mathscr{I}^M$ is continuous on the closed subset $G_M$, as $|\mathcal{F} \circ \mathscr{I}^M(z) - \mathcal{F} \circ \mathscr{I}^M(w)| \lesssim_{\mathcal{F}[\cdot]} \mathcal{W}_2(\mathscr{I}^M(z), \mathscr{I}^M(w)) \lesssim_M |z-w|$, for any $z, w \in G_M$. The essential point of the construction here is the smoothness.

for any $k\in\mathbb{Z}_+^d$, $\varphi=\sqrt{2}\cos(2\pi k\cdot)$, $\bar{\varphi}=\sqrt{2}\sin(2\pi k\cdot)$ and $z\in B_M$. A similar formula holds for the $\sqrt{2}\sin(2\pi k\cdot)$ coordinate. It follows from the construction that

$$([(\mathcal{I}^M)^{-1}\mu]_\varphi)_{\varphi\in\mathcal{D}_M}\mapsto\text{R.H.S. of (A.1)}$$

is smooth on $\mathbb{R}^{\mathcal{D}_M}$, since the smooth function $\rho_M$ acts as a mollification, and we denote this function as $F^M$ in Proposition 3.4.

We now prove that $D_\mu\mathcal{F}^M$ approximates $D_\mu\mathcal{F}$ uniformly and omit the proof of the second order to avoid repetition. It follows from the definition of interior derivative (A.15) that

$$D_\mu\mathcal{F}^M[\mu](x) = \frac{M-1}{M}\int_{\mathbb{R}^{\mathcal{D}_M}}\left\{D_\mu\mathcal{F}\left[(1-\frac{1}{M})\mu*\varphi^M+\frac{1}{M}\mathcal{I}_M(z)\right]*\varphi^M\right\}(x)\rho_M(z)\mathrm{d}z,$$

which is apparently jointly Lipschitz in $(x,\mu)$ by the regulariy of $\mathcal{F}[\cdot]$. Uniform approximation follows from

$$|D_\mu\mathcal{F}^M[\mu](x)-D_\mu\mathcal{F}[\mu](x)| \leq \frac{1}{M}\sup_{\substack{\mu\in\mathcal{P}(\mathbb{T}^d)\\x\in\mathbb{T}^d}}|D_\mu\mathcal{F}[\mu](x)|$$
$$+\int_{\mathbb{R}^{\mathcal{D}_M}}|D_\mu\mathcal{F}[\mu_z^M]*\varphi^M(x)-D_\mu\mathcal{F}[\mu](x)|\rho_M(z)\mathrm{d}z,$$

where we use the fact that $(\rho_M)$ is an approximation identity and abbreviate

$$\mu_z^M:=(1-\frac{1}{M})\mu*\varphi^M+\frac{1}{M}\mathcal{I}_M(z).$$

The right side converges to zero uniformly with $\mu\in\mathcal{P}(\mathbb{T}^d)$ and $x\in\mathbb{T}^d$ due to the same argument we have applied in the proof of Lemma 3.3.

## A.4 Linearization of the McKean–Vlasov SDEs

In this section we outline the derivation of the *linearized equation* of the McKean–Vlasov SDEs (3.22), which is essentially due to the seminal works [BLPR17] and [CDLL19]. Since the rigorous derivation (even of higher-order interior derivatives) is already well established now majorly for its application in the quantative weak propagation of chaos, see e.g. [Tse21, DT21, DF22], we refer the readers to those papers for the comprehensive investigation, especially the proof that the formally derived equations below indeed characterize the differentiation with respect to the probability measures.

Firstly, we have the following chain rule due to the high regularity we assume on the functionals in concern:

$$\frac{\delta}{\delta\mu}\{\mathcal{F}[m_t^\mu]\}(y) \stackrel{(A.14)}{=} \left\langle\frac{\delta}{\delta\mu}\{\mathcal{F}[m_t^\mu]\},\delta_y-\mu\right\rangle = \left\langle\delta_{\mu,y}m_t^\mu,\frac{\delta\mathcal{F}}{\delta\mu}[m_t^\mu]\right\rangle, \quad (A.22)$$

for $t\in[0,T]$ as a special case of Theorem 2.6 in [Tse21] due to the normalization rule. We denote formally $\delta_{\mu,y}m^\mu(\cdot)$ (respectively the second-order derivative $\delta_{\mu,y,z}^2 m^\mu(\cdot)$ defined below) as the linearization (infinitesimal perturbation) of $m^\mu$ in the direction of $\delta_y$ (respectively in the direction of $\delta_y\otimes\delta_z$ for the second-order derivative), that is

$$\delta_{\mu,\nu}m_t^\mu = \int_{\mathbb{T}^d}(\delta_{\mu,y}m_t^\mu)(\nu-\mu)(\mathrm{d}y) = \lim_{\epsilon\to 0+}\frac{1}{\epsilon}(m_t^{\mu+\epsilon(\nu-\mu)}-m_t^\mu) \quad (A.23)$$

in the space of tempered distributions[A.3]. By chain rule again, we have the following second-order linearization assuming sufficient regularity on each term:

$$\frac{\delta^2}{\delta\mu^2}\{\mathcal{F}[m_t^\mu]\}(\mu,y,z) = \left\langle \delta^2_{\mu,y,z}m_t^\mu, \frac{\delta\mathcal{F}}{\delta\mu}[m_t^\mu]\right\rangle$$
$$+ \int_{\mathbb{T}^{2d}} \frac{\delta^2\mathcal{F}}{\delta\mu^2}[m_t^\mu](u,w)\delta_{\mu,y}m_t^\mu(du)\delta_{\mu,z}m_t^\mu(dw). \quad (A.24)$$

**Remark A.7.** The equation (A.23) above is just the equation (1.12) in [Tse21], and $\delta_{\mu,\nu}m^\mu$ is considered as the distributional solution to a linearised forward Kolmogorov equation, that is the equation (1.14) in [Tse21]. We emphasize, however, that to rigorously establish (A.22), we need to prove first that $\mu \mapsto \mathcal{F}[m_t^\mu]$ is indeed differentiable in the sense of (A.13). This is where the decoupling argument in [BLPR17] comes into play.

We denote therefore $\Phi_{\mu,\gamma}(t) := P_{0,t}^\mu(\gamma)$ for $\gamma \in \mathcal{P}(\mathbb{T}^d)$ with the evolution system defined in (4.6) to be the solution to the forward Kolmogorov equation corresponding to the decoupled SDEs

$$\begin{cases} dX_s^{\gamma,\mu} = b(X_s^{\gamma,\mu}, m_s^\mu)ds + \sqrt{2}\Sigma(X_s^{\gamma,\mu}, m_s^\mu) \cdot dB_s, & s \in [0,T], \\ \mathcal{L}(X_0^{\gamma,\mu}) = \gamma. \end{cases} \quad (A.25)$$

It follows that for any $\xi \in C(\mathbb{T}^d)$, we have

$$\int_{\mathbb{T}^d}\xi(x)m_t^\mu(x)dx = \mathbf{E}[\xi(X_t^\mu)] = \int_{\mathbb{T}^d}\mathbf{E}[\xi(X_t^{\delta_x,\mu})]d\mu(x) = \int_{\mathbb{T}^d}(P_{0,t}^\mu)^*\xi(x)d\mu(x), \quad (A.26)$$

where $t \mapsto (P_{0,t}^\mu)^*$ denotes the backward flow of the decoupled equation. This indicates that

$$\frac{d}{d\varepsilon}\bigg|_{\varepsilon=0^+}\int_{\mathbb{T}^d}\xi(x)m_t^{\mu+\varepsilon(\nu-\mu)}dx = \int_{\mathbb{T}^d}(P_{0,t}^\mu)^*\xi(x)d(\nu-\mu)(x)$$
$$+ \int_{\mathbb{T}^d}\frac{d}{d\varepsilon}|_{\varepsilon=0^+}(P_{0,t}^{\mu+\varepsilon(\nu-\mu)})^*\xi(x)d\mu(x).$$

The second term on the R.H.S. i.e. the linearization of the backward flow can also be expressed in terms of the signed measure $\nu - \mu$ as desired from (A.13), so one can effectively prove the differentiability of the flow of the McKean-Vlasov equation. We refer the readers to Theorem 2.6. in [Tse21] for the complete proof.

We now focus on the PDE side of the dynamics. According to the uniqueness (in law) of the McKean–Vlasov system, we have the equality $m^\mu = \Phi_{\mu,\mu}$, which implies the chain rules

$$\partial_{\mu,y}m^\mu = \partial_{\mu,y}\Phi_{\mu,\gamma}|_{\gamma=\mu} + \partial_{\gamma,y}\Phi_{\mu,\gamma}|_{\gamma=\mu}, \quad (A.27)$$

which is just (A.26) in its probabilistic formulation and

$$\partial^2_{\mu,y,z}m^\mu = \partial^2_{\mu,y,z}\Phi_{\mu,\gamma}|_{\gamma=\mu} + \partial_{\mu,z}\partial_{\gamma,y}\Phi_{\mu,\gamma}|_{\gamma=\mu}, \quad (A.28)$$

---

A.3. Rigorously speaking, one needs to define $\delta_{\mu,\nu}m^\mu$ as the solution to the corresponding linearized PDEs and then proves the limit in (A.23), as is done in Theorem 2.6 in [Tse21].

for any $y, z \in \mathbb{T}^d$, where we have used the fact for the second equality that

$$\partial^2_{\gamma, y, z} \Phi_{\mu, \gamma} = \partial_{\gamma, y} \partial_{\mu, z} \Phi_{\mu, \gamma} = 0,$$

since both linearizations should annihilate the initial data, and the linearized equations are all linear with respect to the initial conditions. As required by the proof of Lemma 4.5, we need to investigate the functions of the form

$$\partial_{\mu, \nu} \Phi_{\mu, \gamma} = \int_{\mathbb{T}^d} (\partial_{\mu, y} \Phi_{\mu, \gamma})(\nu - \mu)(\mathrm{d} y) = \lim_{\epsilon \to 0^+} \frac{1}{\epsilon} (\Phi_{\mu + \epsilon(\nu - \mu), \gamma} - \Phi_{\mu, \gamma})$$

for any $\nu \in \mathcal{P}$ regarding (A.3). We summarize the postulated equations appearing in the first two orders of linearization as follows:

1. The first-order linearization of the initial condition in the direction of $\nu \in \mathcal{P}$ satisfies

$$\partial_{\gamma, \nu} \Phi_{\mu, \gamma} \stackrel{(4.6)}{=} P^\mu_{0, t}(\nu - \gamma). \tag{A.29}$$

2. Together, we can decompose the first linearization of $m^\mu$ as

$$\partial_{\mu, \nu} m^\mu = \partial_{\mu, \nu} \Phi_{\mu, \gamma}|_{\gamma = \mu} + \partial_{\gamma, \nu} \Phi_{\mu, \gamma}|_{\gamma = \mu}, \tag{A.30}$$

where we first consider the first-order linearization (restricted for $\gamma = \mu$ for a closed equation):

$$\partial_{\mu, \nu} \Phi_{\mu, \gamma}|_{\gamma = \mu} := \int_{\mathbb{T}^d} (\partial_{\mu, y} \Phi_{\mu, \gamma})(\nu - \mu)(\mathrm{d} y)|_{\gamma = \mu}$$

of the flow (in the direction of $\nu \in \mathcal{P}$) is the solution to the equation (which is therefore independent of $\gamma \in \mathcal{P}$)

$$\partial_t (\partial_{\mu, \nu} \Phi_{\mu, \gamma})|_{\gamma = \mu} \stackrel{(A.30)}{=} \nabla \cdot \left( \left[ \int_{\mathbb{T}^d} \frac{\delta a(\mu, \cdot)}{\delta \mu} [m^\mu_t](y) (\partial_{\mu, \nu} \Phi_{\mu, \gamma})|_{\gamma = \mu}(t, \mathrm{d} y) \right] \nabla m^\mu_t \right)$$
$$- \nabla \cdot \left( \left[ \int_{\mathbb{T}^d} \frac{\delta \tilde{b}(\mu, \cdot)}{\delta \mu} [m^\mu_t](y) (\partial_{\mu, \nu} \Phi_{\mu, \gamma})|_{\gamma = \mu}(t, \mathrm{d} y) \right] m^\mu_t \right)$$
$$+ \nabla \cdot (a[m^\mu_t] \nabla (\partial_{\mu, \nu} \Phi_{\mu, \gamma})|_{\gamma = \mu} + \tilde{b}[m^\mu_t](\partial_{\mu, \nu} \Phi_{\mu, \gamma})|_{\gamma = \mu} + \mathcal{R}^\mu_t(\nu - \mu)), \tag{A.31}$$

for $t \in [0, T]$ with null initial condition, where the forcing term is given by

$$\mathcal{R}^\mu_t(\nu - \mu) = \left( \int_{\mathbb{T}^d} \frac{\delta a(\mu, \cdot)}{\delta \mu} [m^\mu_t](y) \mathrm{d} P^\mu_{0, t}(\nu - \mu)(y) \right) \nabla m^\mu_t$$
$$- m^\mu_t \left( \int_{\mathbb{T}^d} \frac{\delta \tilde{b}(\mu, \cdot)}{\delta \mu} [m^\mu_t](y) \mathrm{d} P^\mu_{0, t}(\nu - \mu)(y) \right).$$

(A.30) then establishes (4.12), after we identify

$$\mathcal{J}^{\mu, 1}(\nu - \mu) = (\partial_{\mu, \nu} \Phi_{\mu, \gamma})|_{\gamma = \mu}.$$

We have derived therefore $\partial_{\mu, \nu} m^\mu$ combining the consideration above. Plugging the solution to (A.30) back into (A.30), we derive the first-order linearization $\partial_{\mu, \nu} \Phi_{\mu, \gamma}$ of the flow as the solution to the equation

$$\partial_t \partial_{\mu, \nu} \Phi_{\mu, \gamma} - \nabla \cdot (a[m^\mu_t] \nabla \partial_{\mu, \nu} \Phi_{\mu, \gamma} + \tilde{b}[m^\mu_t] \partial_{\mu, \nu} \Phi_{\mu, \gamma}) = \mathcal{H}_t(\nu - \mu), \tag{A.32}$$

with null initial condition, and the forcing term defined by

$$\mathcal{H}(\nu-\mu) = \nabla \cdot \left( \left[ \int_{\mathbb{T}^d} \frac{\delta a(\mu,\cdot)}{\delta \mu}[m^\mu](y) \partial_{\mu,\nu} m^\mu(\mathrm{d}y) \right] \nabla \Phi_{\mu,\gamma} \right)$$
$$-\nabla \cdot \left( \left[ \int_{\mathbb{T}^d} \frac{\delta \tilde{b}(\mu,\cdot)}{\delta \mu}[m^\mu](y) \partial_{\mu,\nu} m^\mu(\mathrm{d}y) \right] \Phi_{\mu,\gamma} \right),$$

which is well-defined, since now $\partial_{\mu,\nu} m^\mu$ is known.

3. The second-order linearization (first on the initial condition in the direction of $\nu \in \mathcal{P}$, then on the flow in the direction of $\zeta \in \mathcal{P}$)

$$\partial_{\mu,\zeta} \partial_{\gamma,\nu} \Phi_{\mu,\gamma} := \int_{\mathbb{T}^{2d}} (\partial_{\mu,z} \partial_{\gamma,y} \Phi_{\mu,\gamma})(\nu-\gamma)(\mathrm{d}y)(\zeta-\mu)(\mathrm{d}z)$$

is the solution to the equation with null initial condition

$$\partial_t (\partial_{\mu,\zeta} \partial_{\gamma,\nu} \Phi_{\mu,\gamma}) = \nabla \cdot (a[m_t^\mu] \nabla (\partial_{\mu,\zeta} \partial_{\gamma,\nu} \Phi_{\mu,\gamma}))$$
$$-\nabla \cdot (\tilde{b}[m_t^\mu] (\partial_{\mu,\zeta} \partial_{\gamma,\nu} \Phi_{\mu,\gamma})) + \nabla \cdot \mathcal{J}_t(\nu-\gamma, \zeta-\mu), \qquad (A.33)$$

where the forcing term is given by

$$\mathcal{J}(\nu-\gamma, \zeta-\mu) \stackrel{(A.29)}{=} \left( \int_{\mathbb{T}^d} \frac{\delta a(\cdot, \mu)}{\delta \mu}[m^\mu](y) \mathrm{d}(\partial_{\mu,\zeta} m^\mu)(y) \right) \nabla P^\mu_{0,\cdot}(\nu-\gamma)$$
$$- \left( \int_{\mathbb{T}^d} \frac{\delta \tilde{b}(\cdot, \mu)}{\delta \mu}[m^\mu](y) \mathrm{d}(\partial_{\mu,\zeta} m^\mu)(y) \right) P^\mu_{0,\cdot}(\nu-\gamma).$$

We then derive (4.14) from (3) by identifying

$$I^{\mu,1}(\nu-\mu, \zeta-\mu) = (\partial_{\mu,\zeta} \partial_{\gamma,\nu} \Phi_{\mu,\gamma})|_{\gamma=\mu}.$$

4. Lastly, we notice that

$$\partial^2_{\mu,\nu,\zeta} m^\mu \stackrel{(A.28)}{=} (\partial^2_{\mu,\nu,\zeta} \Phi_{\mu,\gamma})|_{\gamma=\mu} + (\partial_{\mu,\zeta} \partial_{\gamma,\nu} \Phi_{\mu,\gamma})|_{\gamma=\mu}, \qquad (A.34)$$

as in (A.28), where we abbreviate $\partial^2_{\mu,\nu,\zeta} = \partial_{\mu,\zeta} \partial_{\mu,\nu}$, and therefore, the second-order linearization (both on the flow in the direction of $\nu, \zeta \in \mathcal{P}$)

$$(\partial^2_{\mu,\nu,\zeta} \Phi_{\mu,\gamma})|_{\gamma=\mu} = \int_{\mathbb{T}^{2d}} (\partial^2_{\mu,y,z} \Phi_{\mu,\gamma})|_{\gamma=\mu}(\nu-\mu)(\mathrm{d}y)(\zeta-\mu)(\mathrm{d}z)$$

is the solution to the equation with null initial condition

$$\partial_t (\partial^2_{\mu,\nu,\zeta} \Phi_{\mu,\gamma}|_{\gamma=\mu}) \stackrel{(A.32)}{=} \nabla \cdot ((a[m_t^\mu] \nabla - \tilde{b}[m_t^\mu])(\partial^2_{\mu,\nu,\zeta} \Phi_{\mu,\gamma}|_{\gamma=\mu}))$$
$$+ \nabla \cdot \left( \left[ \int_{\mathbb{T}^d} \frac{\delta a(\mu,\cdot)}{\delta \mu}[m_t^\mu](y) (\partial^2_{\mu,\nu,\zeta} \Phi_{\mu,\gamma}|_{\gamma=\mu})(\mathrm{d}y) \right] \nabla m_t^\mu \right)$$
$$- \nabla \cdot \left( \left[ \int_{\mathbb{T}^d} \frac{\delta \tilde{b}(\mu,\cdot)}{\delta \mu}[m_t^\mu](y) (\partial^2_{\mu,\nu,\zeta} \Phi_{\mu,\gamma}|_{\gamma=\mu})(\mathrm{d}y) \right] m_t^\mu \right)$$
$$+ \nabla \cdot (\mathcal{T}_t^\mu(\nu-\mu, \zeta-\mu)), \qquad (A.35)$$

where the nonlocal term comes from the decomposition of $\partial^2_{\mu,\nu,\zeta}m^\mu$ in $\mathscr{H}(\nu-\mu)$ using (A.34), and the forcing term is given by

$$\mathscr{T}^\mu(\nu-\mu, \zeta-\mu) = \left(\int_{\mathbb{T}^d} \frac{\delta a(\mu,\cdot)}{\delta \mu}[m^\mu](y)\partial_{\mu,\zeta}m^\mu(\mathrm{d}y)\right) \nabla(\partial_{\mu,\nu}\Phi_{\mu,\gamma}|_{\gamma=\mu})$$
$$-\left(\int_{\mathbb{T}^d} \frac{\delta \tilde{b}(\mu,\cdot)}{\delta \mu}[m^\mu](y)\partial_{\mu,\zeta}m^\mu(\mathrm{d}y)\right) \partial_{\mu,\nu}\Phi_{\mu,\gamma}|_{\gamma=\mu}$$
$$+\left(\int_{\mathbb{T}^d} \frac{\delta a(\mu,\cdot)}{\delta \mu}[m^\mu](y)\partial_{\mu,\nu}m^\mu(\mathrm{d}y)\right) \nabla(\partial_{\mu,\zeta}\Phi_{\mu,\gamma}|_{\gamma=\mu})$$
$$-\left(\int_{\mathbb{T}^d} \frac{\delta \tilde{b}(\mu,\cdot)}{\delta \mu}[m^\mu](y)\partial_{\mu,\nu}m^\mu(\mathrm{d}y)\right) \partial_{\mu,\zeta}\Phi_{\mu,\gamma}|_{\gamma=\mu}$$
$$+\left(\int_{\mathbb{T}^d} \frac{\delta a(\mu,\cdot)}{\delta \mu}[m^\mu](y)\partial_{\mu,\zeta}\partial_{\gamma,\nu}\Phi_{\mu,\gamma}|_{\gamma=\mu}(\mathrm{d}y)\right) \nabla m^\mu$$
$$-\left(\int_{\mathbb{T}^d} \frac{\delta \tilde{b}(\mu,\cdot)}{\delta \mu}[m^\mu](y)\partial_{\mu,\zeta}\partial_{\gamma,\nu}\Phi_{\mu,\gamma}|_{\gamma=\mu}(\mathrm{d}y)\right) m^\mu$$
$$+\left(\int_{\mathbb{T}^{2d}} \frac{\delta^2 a(\mu,\cdot)}{\delta \mu^2}[m^\mu](y,z)\partial_{\mu,\nu}m^\mu(\mathrm{d}y)\partial_{\mu,\zeta}m^\mu(\mathrm{d}z)\right) \nabla m^\mu$$
$$-\left(\int_{\mathbb{T}^{2d}} \frac{\delta^2 \tilde{b}(\mu,\cdot)}{\delta \mu^2}[m^\mu](y,z)\partial_{\mu,\nu}m^\mu(\mathrm{d}y)\partial_{\mu,\zeta}m^\mu(\mathrm{d}z)\right) m^\mu. \quad (\text{A.36})$$

Then by evaluating at $\gamma=\mu$, we can close the equation (4.15) with respect to $\partial^2_{\mu,\nu,\zeta}\Phi_{\mu,\gamma}|_{\gamma=\mu}$.

We can now conclude that the expression (4.12) holds using (A.3), (A.27), (A.28) and the corresponding equations, while $\phi, \psi$ in (4.12) are instead assumed to be of the forms $\phi \in \mathcal{P}(\mathbb{T}^d) - \mu$, $\psi \in \mathcal{P}(\mathbb{T}^d) - \mu$ according to the definition of exterior derivatives. We finally recall that we extended $R^\mu$ onto $[0,T] \times W^{-3,1}(\mathbb{T}^d)$ in (4.3), and we now prove that the same holds for $\mathscr{T}^\mu$.

**Lemma A.8.** *The map* $\mathscr{R}^\mu : [0,T] \times W^{-3,1}(\mathbb{T}^d) \to W^{-1,1}(\mathbb{T}^d)$ *satisfies*

$$\sup_{\mu \in L^1 \cap \mathcal{P}} \sup_{0 \le t \le T} \|\mathscr{R}^\mu_t(\psi)\|_{W^{-1,1}} \lesssim \|\psi\|_{W^{-3,1}}.$$

*If we further assume that (4.17) holds, then the function* $\mathscr{T}^\mu$ *defined in (A.35) naturally extends to a function* $\mathscr{T}^\mu : [0,T] \times (W^{-3,1})^{\otimes 2}(\mathbb{T}^d) \to W^{-2,1}(\mathbb{T}^d)$ *such that*

$$\sup_{\mu \in L^1 \cap \mathcal{P}} \sup_{0 \le t \le T} \|\mathscr{T}^\mu_t(\psi)\|_{W^{-2,1}} \lesssim \|\psi\|_{W^{-3,1}},$$

**Proof.** We recall that

$$\mathscr{R}^\mu_t(\psi) = \left(\int_{\mathbb{T}^d} \frac{\delta a(\mu,\cdot)}{\delta \mu}[m^\mu_t](y) \mathrm{d}P^\mu_{0,t}(\psi)(y)\right) \nabla m^\mu_t - m^\mu_t \left(\int_{\mathbb{T}^d} \frac{\delta \tilde{b}(\mu,\cdot)}{\delta \mu}[m^\mu_t](y) \mathrm{d}P^\mu_{0,t}(\psi)(y)\right),$$

and the boundedness of $\mathscr{R}^\mu$ follows from

$$\left\|\left(\int_{\mathbb{T}^d} \frac{\delta a(\mu,\cdot)}{\delta \mu}[m^\mu_t](y) P^\mu_{0,t}\psi(\mathrm{d}y)\right) \nabla m^\mu_t\right\|_{W^{-1,1}} \lesssim \left\|\frac{\delta a}{\delta \mu}[m^\mu_t]\right\|_{W^{3,\infty}} \|\psi\|_{W^{-3,1}} \|m^\mu_t\|_{L^1}$$
$$\overset{(4.5)}{\lesssim} \sup_{\mu \in \mathcal{P} \cap L^1} \left\|\frac{\delta a}{\delta \mu}\right\|[\mu]_{W^{3,\infty}} \|\psi\|_{W^{-3,1}},$$

which is uniform in $t \in [0,T]$. We have used the continuity of $\nabla: L^1 \to W^{-1,1}$ and the joint regularity of

$$\mathbb{T}^{2d} \ni (x,y) \mapsto \frac{\delta a}{\delta \mu}[m_t^\mu](x;y) := \frac{\delta a(\mu, x)}{\delta \mu}[m_t^\mu](y).$$

Similarly, we have the estimate with the drift term $b[\cdot]$ that

$$\left\| \left( \int_{\mathbb{T}^d} \frac{\delta \tilde{b}(\mu,\cdot)}{\delta \mu}[m_t^\mu](y) P_{0,t}^\mu \psi(\mathrm{d}y) \right) m_t^\mu \right\|_{L^1} \lesssim \sup_{\mu \in \mathcal{P} \cap L^1} \left\| \frac{\delta \tilde{b}}{\delta \mu}[\mu] \right\|_{W^{3,\infty}} \|\psi\|_{W^{-3,1}}.$$

We can now extend $\mathcal{T}^\mu$ in (A.35) by

$$\begin{aligned}
\mathcal{T}^\mu(\psi, \phi) &= \left( \int_{\mathbb{T}^d} \frac{\delta a(\mu,\cdot)}{\delta \mu}[m^\mu](y) \mathcal{J}^\mu(\phi)(\mathrm{d}y) \right) \nabla(\mathcal{J}^{\mu,1}(\psi)) \\
&\quad - \left( \int_{\mathbb{T}^d} \frac{\delta \tilde{b}(\mu,\cdot)}{\delta \mu}[m^\mu](y) \mathcal{J}^\mu(\phi)(\mathrm{d}y) \right) \mathcal{J}^{\mu,1}(\psi) \\
&\quad + \left( \int_{\mathbb{T}^d} \frac{\delta a(\mu,\cdot)}{\delta \mu}[m^\mu](y) \mathcal{J}^\mu(\psi)(\mathrm{d}y) \right) \nabla(\mathcal{J}^{\mu,1}(\phi)) \\
&\quad - \left( \int_{\mathbb{T}^d} \frac{\delta \tilde{b}(\mu,\cdot)}{\delta \mu}[m^\mu](y) \mathcal{J}^\mu(\psi)(\mathrm{d}y) \right) \mathcal{J}^{\mu,1}(\phi) \\
&\quad + \left( \int_{\mathbb{T}^d} \frac{\delta a(\mu,\cdot)}{\delta \mu}[m^\mu](y) I^{\mu,1}(\psi,\phi)(\mathrm{d}y) \right) \nabla m^\mu \\
&\quad - \left( \int_{\mathbb{T}^d} \frac{\delta \tilde{b}(\mu,\cdot)}{\delta \mu}[m^\mu](y) I^{\mu,1}(\psi,\phi)(\mathrm{d}y) \right) m^\mu \\
&\quad + \left( \int_{\mathbb{T}^{2d}} \frac{\delta^2 a(\mu,\cdot)}{\delta \mu^2}[m^\mu](y,z) \mathcal{J}^\mu(\psi)(\mathrm{d}y) \mathcal{J}^\mu(\phi)(\mathrm{d}z) \right) \nabla m^\mu \\
&\quad - \left( \int_{\mathbb{T}^{2d}} \frac{\delta^2 \tilde{b}(\mu,\cdot)}{\delta \mu^2}[m^\mu](y,z) \mathcal{J}^\mu(\psi)(\mathrm{d}y) \mathcal{J}^\mu(\phi)(\mathrm{d}z) \right) m^\mu.
\end{aligned} \qquad (A.37)$$

Therefore, we have

$$\begin{aligned}
\sup_{t \in [0,T]} \|\mathcal{T}_t^\mu\|_{W^{-2,1}} &\lesssim \max_{j=1,2} \sup_{\mu \in L^1 \cap \mathcal{P}} \left\| \frac{\delta^j a}{\delta \mu^j}[\mu] \right\|_{W^{4,\infty}} \|\psi\|_{W^{-3,1}} \|\phi\|_{W^{-3,1}} \\
&\quad + \max_{j=1,2} \sup_{\mu \in L^1 \cap \mathcal{P}} \left\| \frac{\delta^j \tilde{b}}{\delta \mu^j}[\mu] \right\|_{W^{4,\infty}} \|\psi\|_{W^{-3,1}} \|\phi\|_{W^{-3,1}},
\end{aligned} \qquad (A.38)$$

since for example the worst term is controlled by

$$\sup_{t \in [0,T]} \left\| \int_{\mathbb{T}^d} \frac{\delta a(\mu,\cdot)}{\delta \mu}[m_t^\mu](y) (\mathcal{J}_t^\mu \phi)(\mathrm{d}y) \nabla(\mathcal{J}_t^{\mu,1} \psi) \mathrm{d}s \right\|_{W^{-2,1}}$$

$$\lesssim \sup_{\mu \in L^1 \cap \mathcal{P}} \left\| \frac{\delta a}{\delta \mu}[\mu] \right\|_{W^{4,\infty}} \sup_{t \in [0,T]} \|\mathcal{J}_t^\mu \phi\|_{W^{-3,1}} \sup_{t \in [0,T]} \|\mathcal{J}_t^{\mu,1} \psi\|_{W^{-1,1}}$$

$$\overset{(4.17)}{\lesssim} \sup_{\mu \in L^1 \cap \mathcal{P}} \left\| \frac{\delta a}{\delta \mu}[\mu] \right\|_{W^{4,\infty}} \|\phi\|_{W^{-3,1}} \|\psi\|_{W^{-3,1}},$$

while the rest are controlled in a similar manner. $\square$

# Acknowledgement

We gratefully acknowledge funding by Deutsche Forschungsgemeinschaft (DFG) under Germany's Excellence Strategy – The Berlin Mathematics Research Center MATH+ (EXC-2046/1, project ID: 390685689) and through the grant CRC 1114 "Scaling Cascades in Complex Systems", Project Number 235221301, Project C10 "Numerical analysis for nonlinear SPDE models of particle systems".

The second author thanks François Delarue, Jin Feng, Zhengyan Wu and Wei Huang for inspiring discussions and Lorenzo Dello Schiavo for carefully reading a preliminary draft of the paper and providing valuable suggestions and insights.